\documentclass[12pt]{amsart}
\usepackage{amsthm,amssymb,amsmath,amscd}

\newtheorem{thm}{Theorem}[section]
\newtheorem{lm}[thm]{Lemma}
\newtheorem{cor}[thm]{Corollary}
\newtheorem{pro}[thm]{Proposition}
\newtheorem{pro-def}[thm]{Proposition and Definition}

\theoremstyle{definition}
\newtheorem{df}[thm]{Definition}

\input epsf

\def \int {\text{Int}}

\begin{document}

\title[The invariant of $n$-punctured ball tangles]
{The invariant of $n$-punctured ball tangles}
\author[Jae-Wook Chung]{Jae-Wook Chung}
\address{Department of Mathematics, Yeungnam University,
Kyongsan, Kyongbuk, Korea 712-749} \email{jwchung@ynu.ac.kr,
jwchung@alumni.ucr.edu}
\thanks{This work was supported by the Korea Research Foundation
Grant funded by the Korean Government (MOEHRD, Basic Research
Promotion Fund) (KRF-2007-355-C00006).}

\begin{abstract}{Based on the Kauffman bracket at $A=e^{i
\pi/4}$, we defined an invariant for a special type of
$n$-punctured ball tangles. The invariant $F^n$ takes values in
the set $PM_{2\times2^n}(\mathbb Z)$ of $2\times 2^n$ matrices
over $\mathbb Z$ modulo the scalar multiplication of $\pm1$. We
provide the formula to compute the invariant of the $k_1 + \cdots
+ k_n$-punctured ball tangle composed of given
$n,k_1,\dots,k_n$-punctured ball tangles. Also, we define the
horizontal and the vertical connect sums of punctured ball tangles
and provide the formulas for their invariants from those of given
punctured ball tangles. In addition, we introduce the elementary
operations on the class $\textbf{\textit{ST}}$ of $1$-punctured
ball tangles, called spherical tangles. The elementary operations
on $\textbf{\textit{ST}}$ induce the operations on
$PM_{2\times2}(\mathbb Z)$, also called the elementary operations.
We show that the group generated by the elementary operations on
$PM_{2\times2}(\mathbb Z)$ is isomorphic to a Coxeter group.}
\end{abstract}

\maketitle

\section{Introduction}

Throughout the paper, we work in either the smooth or the
piecewise linear category. For basic terminologies of knot theory,
see \cite{A,BZ}.

We introduced a general definition of an $n$-punctured ball tangle
and basic properties on them in \cite{C-L}. However, our main
interest still lies in a special type of $n$-punctured ball
tangles, each boundary component of which intersects with the
$1$-dimensional proper submanifold at exactly $4$ points. Hence,
we restrict our scope to only such punctured ball tangles
(Definition 2.1). In the case of $n=0$, it corresponds exactly to
Conway's notion of tangles in the 3-ball $B^3$ \cite{Conway}, and
we call them {\it ball tangles}. Using the Kauffman bracket at
$A=e^{i \pi/4}$, we defined an invariant for this special type of
$n$-punctured ball tangles \cite{C-L}. The invariant $F^n(T^n)$
for such an $n$-punctured ball tangle $T^n$ is an element of the
set $PM_{2\times2^n}(\mathbb Z)$ of $2\times 2^n$ matrices over
$\mathbb Z$ modulo the scalar multiplication of $\pm1$. Specially,
$F^0(T^0)$ is Krebes' invariant \cite{K}.

In this paper, we generalize the formula for the invariant of the
ball tangle induced by an $n$-punctured ball tangle and $n$ ball
tangles. The invariant $F^n$ behaves well under the operadic
composition of $n$-punctured ball tangles. As a punctured ball
tangle valued function, an $n$-punctured ball tangle $T^n$ has the
class of all punctured ball tangles as domain. When we put $n$
many punctured ball tangles at the $n$ holes of $T^n$, we have a
new punctured ball tangle. That is, given an $n$-punctured ball
tangle $T^n$ and $k_1,\dots,k_n$-punctured ball tangles
$T^{k_1(1)},\dots,T^{k_n(n)}$, respectively, we consider the $k_1
+ \cdots + k_n$-punctured ball tangle
$T^n(T^{k_1(1)},\dots,T^{k_n(n)})$, where $n \in \mathbb N$ and
$k_1,\dots,k_n \in \mathbb N \cup \{0\}$. In this case, we show
how to compute the invariant
$F^{k_1+\cdots+k_n}(T^n(T^{k_1(1)},\dots,T^{k_n(n)}))$ if
$F^n(T^n), F^{k_1}(T^{k_1(1)}),\dots,F^{k_n}(T^{k_n(n)})$ are
given (Theorem 3.2). Also, we consider the horizontal connect sum
$T^{k_1(1)} +_h T^{k_2(2)}$ and the vertical connect sum
$T^{k_1(1)} +_v T^{k_2(2)}$ of $k_1$ and $k_2$-punctured ball
tangles $T^{k_1(1)}$ and $T^{k_2(2)}$, respectively, and provide
the formulas for the invariants $F^{k_1+k_2}(T^{k_1(1)} +_h
T^{k_2(2)})$ and $F^{k_1+k_2}(T^{k_1(1)} +_v T^{k_2(2)})$ from
$F^{k_1}(T^{k_1(1)})$ and $F^{k_2}(T^{k_2(2)})$ (Theorem 3.3).

These generalizations can reduce much work when we try to compute
the invariant for rather complicated punctured ball tangles. In
order to compute the invariant for a given $n$-punctured ball
tangle, we may successfully decompose it appropriately by already
known ball tangles and punctured ball tangles in terms of
compositions and connect sums. Then we will get the invariant of
it by our formulas.

Finally, we introduce the elementary operations on the class
$\textbf{\textit{ST}}$ of $1$-punctured ball tangles, called
spherical tangles. The elementary operations on
$\textbf{\textit{ST}}$ induce the operations on
$PM_{2\times2}(\mathbb Z)$, which is also called the elementary
operations. We show that the group generated by the elementary
operations on $PM_{2\times2}(\mathbb Z)$ is isomorphic to a
Coxeter group (Theorem 4.6).

\section{Preliminaries}

In this section, we give a bunch of definitions and statements
required for our main theorems. All of them come from our previous
paper \cite{C-L}.

The notion of {\it tangles} was introduced by J. Conway
\cite{Conway} as the basic building blocks of links in the
3-dimensional sphere $S^3$. A tangle $T$ is defined by a pair
$(B^3,T)$, where $B^3$ is a 3-dimensional closed ball and $T$ is a
1-dimensional proper submanifold of $B^3$ with 2 non-circular
components. The points in $\partial T\subset\partial B^3$ will be
fixed all the time. Here, we considered holes inside the tangle
such that if they are filled up with any tangles, we have a new
tangle. In this sense, we define an {\it $n$-punctured ball
tangle} slightly modified that in \cite{C-L} to fit our purpose.

\begin{df} Let $n$ be a nonnegative integer, and let
$H_0$ be a closed 3-ball, and let $H_1,\dots,H_n$ be pairwise
disjoint closed 3-balls contained in the interior ${\rm Int}(H_0)$
of $H_0$. For each $k \in \{0,1,\dots,n\}$, take $4$ distinct
points $a_{k1},a_{k2},a_{k3},a_{k4}$ of $\partial H_k$. Then a
1-dimensional proper submanifold $T$ of $H_0 - \bigcup_{i=1}^n
{\rm Int}(H_i)$ is called an $n$-punctured ball tangle with
respect to $(H_k)_{0 \leq k \leq n}$ and
$((a_{k1},a_{k2},a_{k3},a_{k4}))_{0 \leq k \leq n}$ or, simply, an
$n$-punctured ball tangle if $\partial T = \bigcup_{k=0}^n
\{a_{k1},a_{k2},a_{k3},a_{k4}\}$. Hence, $\partial T \cap \partial
H_k = \{a_{k1},a_{k2},a_{k3},a_{k4}\}$ for each $k \in
\{0,1,\dots,n\}$.
\end{df}

Note that we can regard an $n$-punctured ball tangle $T$ with
respect to $(H_k)_{0 \leq k \leq n}$ and
$((a_{k1},a_{k2},a_{k3},a_{k4}))_{0 \leq k \leq n}$ as a $4$-tuple
$(n,(H_k)_{0 \leq k \leq n},((a_{k1},a_{k2},a_{k3},a_{k4}))_{0
\leq k \leq n},T)$.

Let $n \in \mathbb N \cup \{0\}$, and let $\textbf{\textit{nPBT}}$
be the class of all $n$-punctured ball tangles with respect to
$(H_k)_{0 \leq k \leq n}$ and $((a_{k1},a_{k2},a_{k3},a_{k4}))_{0
\leq k \leq n}$, and let $X = H_0 - \bigcup_{i=1}^n {\rm
Int}(H_i)$. Define $\cong$ on $\textbf{\textit{nPBT}}$ by $T_1
\cong T_2$ if and only if there is a homeomorphism $h:X
\rightarrow X$ such that $h|_{\partial X}={\rm Id}_X|_{\partial
X}$, $h(T_1)=T_2$, and $h$ is isotopic to ${\rm Id}_X$ relative to
the boundary $\partial X$ for all $T_1,T_2 \in
\textbf{\textit{nPBT}}$. Then $\cong$ is an equivalence relation
on $\textbf{\textit{nPBT}}$, where ${\rm Id}_X$ is the identity
map from $X$ to $X$. $n$-punctured ball tangles $T_1$ and $T_2$ in
$\textbf{\textit{nPBT}}$ are said to be equivalent or of the same
isotopy type if $T_1 \cong T_2$. Also, for each $n$-punctured ball
tangle $T$ in $\textbf{\textit{nPBT}}$, the equivalence class of
$T$ with respect to $\cong$ is denoted by $[T]$. Without any
confusion, we will also use $T$ for $[T]$.

Like link diagrams, to deal with diagrams of $n$-punctured ball
tangles in the same isotopy type, we need Reidemeister moves among
them. For link diagrams or ball tangle diagrams, we have 3 kinds
of Reidemeister moves. However, we need one and only one more kind
of moves which are called the Reidemeister moves of type IV.

\vskip 0.1in

\bigskip
\centerline{\epsfxsize=5 in \epsfbox{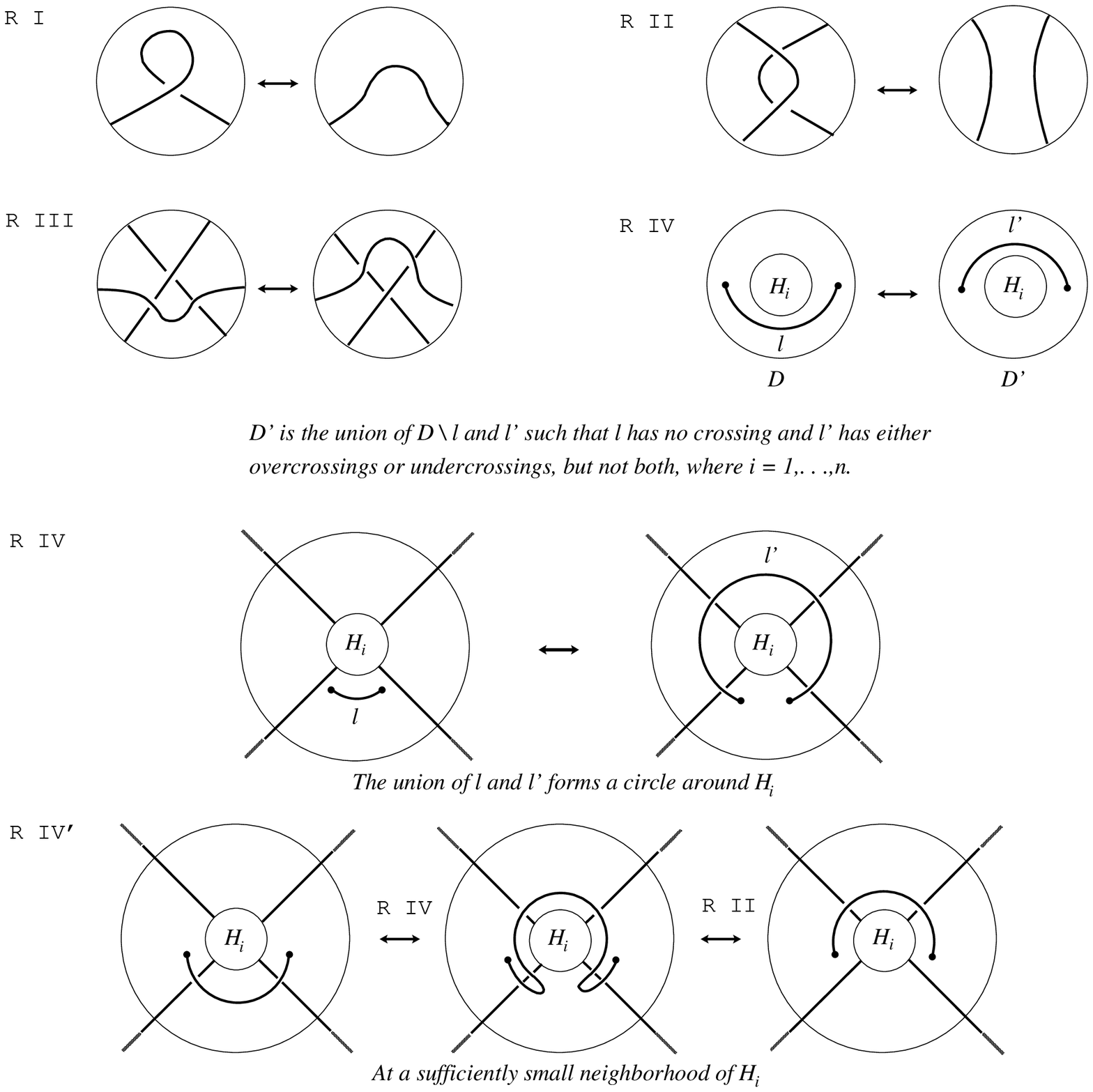}}
\medskip
\centerline{\small Figure 1. Tangle Reidemeister moves.}
\bigskip

The Reidemeister moves for diagrams of $n$-punctured ball tangles
are illustrated in Figure 1. Like link diagrams, tangle diagrams
also have Reidemeister Theorem involving the Reidemeister moves of
type IV. Let us call Reidemeister moves including type IV Tangle
Reidemeister moves.

\begin{thm} Let $n$ be a nonnegative integer, and let $D_1$ and $D_2$
be diagrams of $n$-punctured ball tangles. Then $D_1 \cong D_2$ if
and only if $D_2$ can be obtained from $D_1$ by a finite sequence
of Tangle Reidemeister moves.
\end{thm}

There are many models for a class of $n$-punctured ball tangles.
It is convenient to use normalized ones. One model for a class of
$n$-punctured ball tangles is described in \cite{C-L}.

Our invariant is based on the Kauffman bracket at $A=e^{i \pi/4}$.
Recall the Kauffman bracket is a regular isotopy invariant of link
diagrams. That is, it will not be changed under Reidemeister moves
of type II and III.

Note that a state $\sigma$ of a link diagram $L$ with $n$
crossings $c_1,\dots,c_n$ is regarded as a function
$\sigma:\{c_1,\dots,c_n\} \rightarrow \{A,B\}$, where $A$ and $B$
are the $A$-type and $B$-type splitting functions, respectively.
Therefore, a link diagram $L$ with $n$ crossings has exactly $2^n$
states of it. Apply a state $\sigma$ to $L$ in order to change $L$
to a diagram $L_\sigma$, called the resolution of $L$ by $\sigma$,
without any crossing.

\vskip 0.1in

\bigskip
\centerline{\epsfxsize=4 in \epsfbox{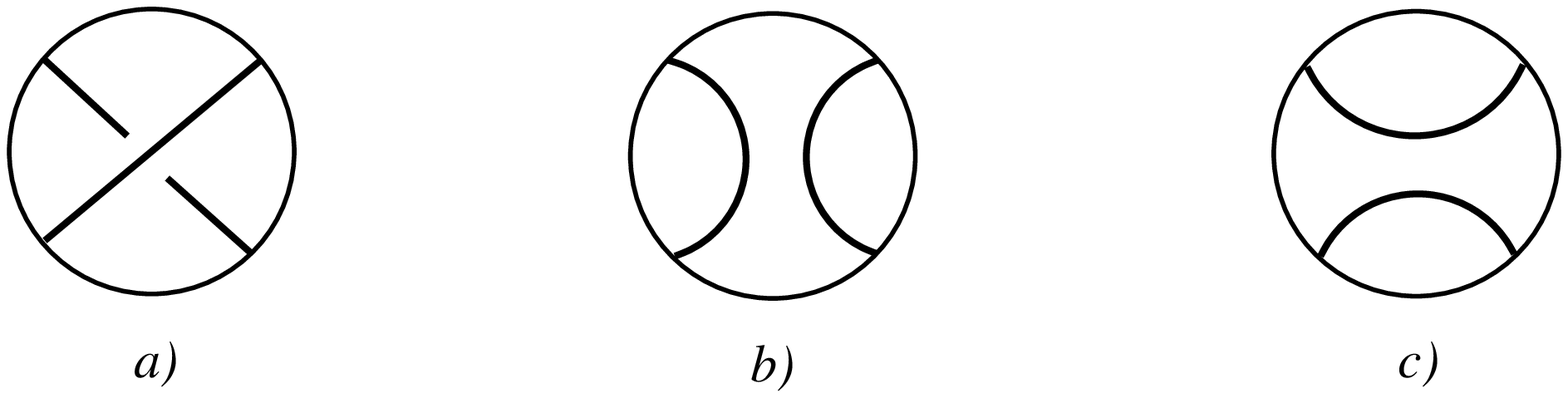}}
\medskip
\centerline{\small a) a crossing $c$ of $L$, b) the part of
$L_\sigma$ by $\sigma(c)=A$, c) the part of $L_\sigma$ by
$\sigma(c)=B$.}
\centerline{\small Figure 2. Two types of splitting
of a crossing of $L$.}
\bigskip

\begin{df} Let $L$ be a link diagram. Then the Kauffman bracket
$\langle L \rangle_A$, or simply, $\langle L \rangle$, is defined
by $$\langle L \rangle_A=\sum_{\sigma \in S}
A^{\alpha(\sigma)}(A^{-1})^{\beta(\sigma)}
(-A^2-A^{-2})^{d(\sigma)-1},$$ where $S$ is the set of all states
of $L$, $\alpha(\sigma)=|\sigma^{-1}(A)|$,
$\beta(\sigma)=|\sigma^{-1}(B)|$, and $d(\sigma)$ is the number of
circles in $L_\sigma$.
\end{df}

We have the following skein relation of the Kauffman bracket.

\begin{pro} Let $L$ be a link diagram, and let $c$ be a crossing
of $L$. Then if $L_A$ and $L_B$ are link diagrams obtained from
$L$ by $A$-type splitting and $B$-type splitting only at $c$,
respectively, then $\langle L \rangle = A\langle L_A \rangle +
A^{-1}\langle L_B \rangle$.
\end{pro}

\begin{proof} Suppose that $S$ is the set of all states
of $L$ and $S_A=\{\sigma \in S|\sigma(c)=A\}$ and $S_B=\{\tau \in
S|\tau(c)=B\}$. Then $\langle L \rangle = A\sum_{\sigma \in S_A}
A^{\alpha(\sigma)-1}(A^{-1})^{\beta(\sigma)}
(-A^2-A^{-2})^{d(\sigma)-1} + A^{-1}\sum_{\tau \in S_B}
A^{\alpha(\tau)}(A^{-1})^{\beta(\tau)-1} (-A^2-A^{-2})^{d(\tau)-1}
= A\langle L_A \rangle + A^{-1}\langle L_B \rangle$ because $S$ is
the disjoint union of $S_A$ and $S_B$. This proves the
proposition.
\end{proof}

Following \cite{K}, a state $\sigma$ of a link diagram $L$ is
called a monocyclic state of $L$ if $d(\sigma)=1$. That is, we
have only one circle when we remove all crossings of $L$ by
$\sigma$.

Also, it is proved in \cite{K} that monocyclic states $\sigma$ and
$\sigma'$ of $L$ differ at an even number of crossings. The
following lemma is a generalization of this statement.

\begin{lm} [J.-W. Chung and X.-S. Lin \cite{C-L}] Let $L$ be a link
diagram. Then states $\sigma$ and $\sigma'$ of $L$ are of the same
parity, i.e., $d(\sigma) \equiv d(\sigma')$ {\rm mod} 2, if and
only if $\sigma$ and $\sigma'$ differ at an even number of
crossings, where $d(\sigma)$ and $d(\sigma')$ are the numbers of
circles in $L_\sigma$ and $L_{\sigma'}$, respectively.
\end{lm}

\begin{proof} Let $\sigma$ be a state of a link diagram $L$ with
$n$ crossings $c_1,\dots,c_n$. Change the value of $\sigma$ at
only one crossing $c_i$ to get another state $\sigma_i$ and
observe what happens to $d(\sigma_i)$, where $1 \leq i \leq n$. We
claim that $\sigma$ and $\sigma_i$ have different parities, more
precisely, $d(\sigma)=d(\sigma_i) \pm1$. Hence, we will have
$d(\sigma) \equiv d(\sigma_i)+1$ {\rm mod} 2. Now, to consider
$\sigma(c_i)$ and $\sigma_i(c_i)$, take a sufficiently small
neighborhood $B_i$ at the projection of $c_i$ so that the
intersection of ${\rm Int}(B_i)$ and the set of all double points
of $L$ is the projection of $c_i$ and the intersection of
$\partial B_i$ and the projection of $L$ has exactly 4 points on
the projection plane of $L$ which are not double points of $L$.

{\it Case 1}. If these 4 points are on a circle in $L_\sigma$,
then
$$d(\sigma_i)=d(\sigma)+1.$$

{\it Case 2}. If two of 4 points are on a circle and the other
points are on another circle in $L_\sigma$, then
$$d(\sigma_i)=d(\sigma)-1.$$

\vskip 0.1in

\bigskip
\centerline{\epsfxsize=5 in \epsfbox{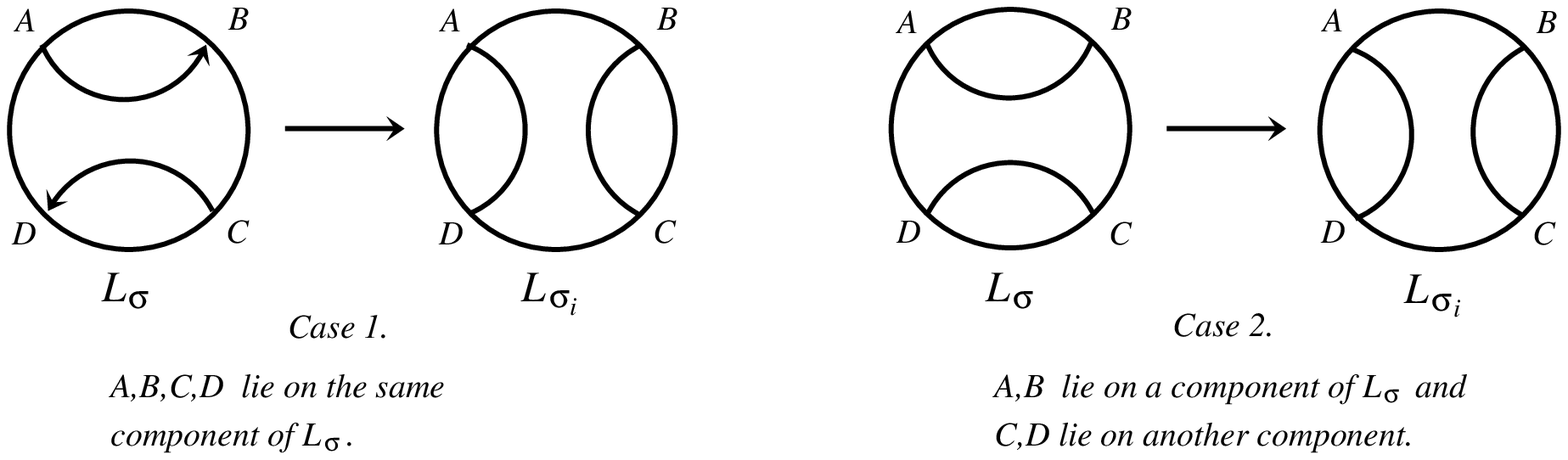}}
\medskip
\centerline{\small Figure 3. Proof of Lemma 2.5.}
\bigskip

Now, it is easy to show the lemma. Suppose that $\sigma$ and
$\sigma'$ are states of $L$ which differ at $k$ crossings of $L$
for some $k \in \{0,1,\dots,n\}$. Then $d(\sigma') \equiv
d(\sigma)+k$ {\rm mod} 2. If $d(\sigma) \equiv d(\sigma')$ {\rm
mod} 2, then $k$ is even. Conversely, if $d(\sigma) \equiv
d(\sigma')+1$ {\rm mod} 2, then $k+1$ is even, that is, $k$ is
odd. This proves the lemma.
\end{proof}

Suppose that $A=e^{i \pi/4}$. Then $-A^2-A^{-2}=0$. Therefore,
$$\langle L \rangle=\sum_{\sigma \in M}
A^{\alpha(\sigma)-\beta(\sigma)},$$ where $M$ is the set of all
monocyclic states of $L$.

From now on, we use only the Kauffman brackets at $A=e^{i \pi/4}$.
Note that, since $|A|=1$, the determinant $|\langle L \rangle|$ of
$L$ is an isotopy invariant.

\begin{lm} If $L$ is a link diagram, then there are $p \in \mathbb Z$ and
$u \in \mathbb C$ such that $u^8=1$ and $\langle L \rangle=pu$.
\end{lm}

The following notations throughout the rest of the paper:

\noindent$\bullet$ $\Phi = \{z \in \mathbb{C}\,|\,z^8 =1\} =
\{A^k\,|\,k \in\mathbb{Z}\}$ and
$\mathbb{Z}\Phi=\{kz\,|\,k\in\mathbb{Z},z\in\Phi\}$.

\noindent$\bullet$ $PM_{m\times n}(\mathbb{Z})$ is the quotient of
$M_{m\times n}(\mathbb{Z})$ under the scalar multiplication by
$\pm1$.

\noindent$\bullet$ $\textbf{\textit{BT}}$ is the class of diagrams
of 0-punctured ball tangles (i.e. ball tangles).

\noindent$\bullet$ $\textbf{\textit{ST}}$ is the class of diagrams
of 1-punctured ball tangles (they will be called spherical
tangles).

\begin{pro} If $a,b,k,l \in \mathbb{Z}$, then $aA^k + bA^l \in
\mathbb{Z}\Phi$ if and only if $ab=0$ or $k \equiv l$ {\rm mod} 4.
\end{pro}

Given a ball tangle diagram $B$, consider 2 kinds of closures as
in Figure 4. The link diagrams $B_1$ and $B_2$ are called the
numerator closure and the denominator closure of $B$,
respectively. A monocyclic state of $B_1$ is called a numerator
state of $B$ and that of $B_2$ is a denominator state of $B$.

\bigskip
\centerline{\epsfxsize=5 in \epsfbox{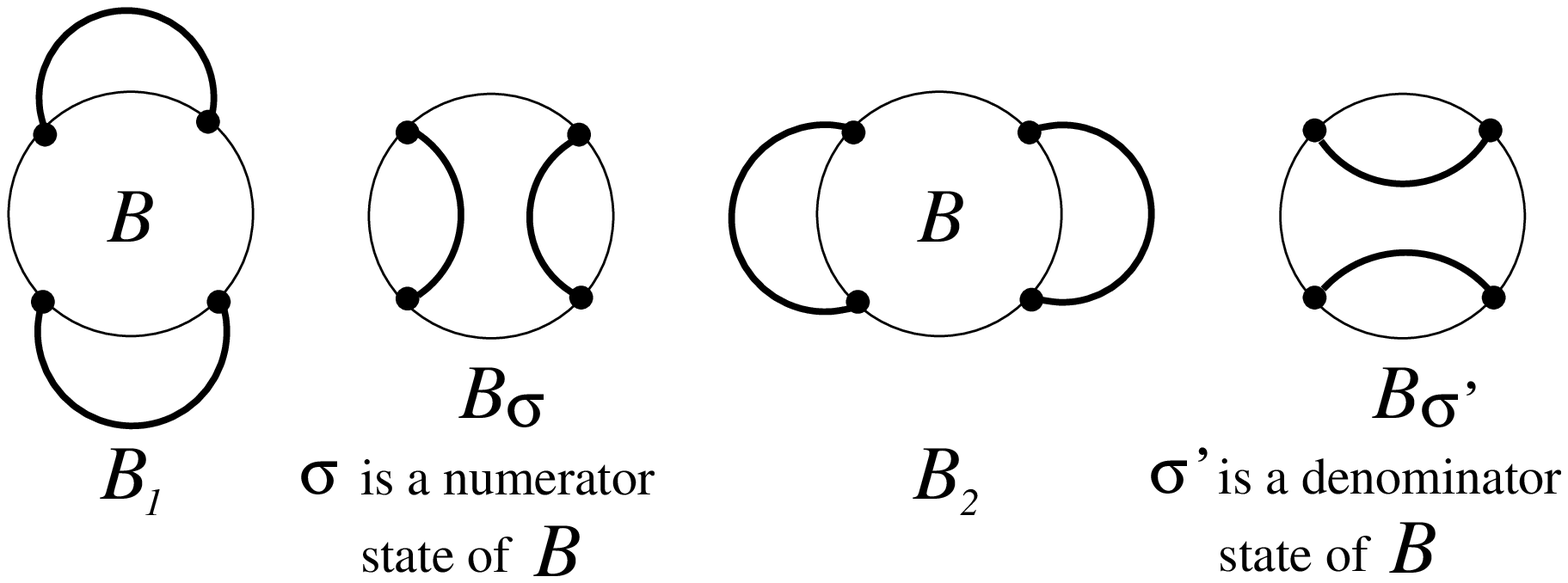}}
\medskip
\centerline{\small Figure 4. The numerator closure $B_1$ and the
denominator closure $B_2$.}
\bigskip

Notice that a numerator state $\sigma$ and a denominator state
$\sigma'$ of a ball tangle diagram $B$ differ at an odd number of
crossings. To see this, we think of a link diagram $L$ such that
$B$ embeds in $L$ and $L$ has one and only one more crossing $c$
at the outside of the ball containing $B$ and $L$ has no
self-twist at the outside of the ball. We have two monocyclic
states of $L$ from the numerator state $\sigma$ and the
denominator state $\sigma'$, respectively, which differ at $c$.
Hence, $\sigma$ and $\sigma'$ differ at an odd number of
crossings. Without loss of generality, we may assume that
$$\langle L\rangle=A\langle B_1\rangle+A^{-1}\langle B_2\rangle\in\mathbb{Z}
\Phi.$$ If $\langle B_1\rangle = pA^k$ and $\langle
B_2\rangle=qA^l$, by Proposition 2.7, we have $l\equiv k+2$ {\rm
mod} 4. Hence, there is a unique $(\alpha,\beta) \in \mathbb{Z}^2$
such that
$$\left\{\begin{pmatrix} z\langle B_1 \rangle \\ iz\langle B_2 \rangle
\end{pmatrix}\,|\,z \in \Phi\right\} \cap M_{2\times 1}(\mathbb{Z})=\left\{
\begin{pmatrix} \alpha \\ \beta \end{pmatrix}, \begin{pmatrix} -\alpha \\ -\beta
\end{pmatrix}\right\}:=\left[\begin{matrix} \alpha \\ \beta
\end{matrix}\right] \in PM_{2\times 1}(\mathbb{Z}).$$

\begin{df} (Krebes \cite{K}) Define $f:\textbf{\textit{BT}} \rightarrow PM_{2\times1}
(\mathbb{Z})$ by
$$f(B)=\left\{\begin{pmatrix} z\langle B_1 \rangle \\ iz\langle B_2
\rangle \end{pmatrix}\,|\,z \in \Phi\right\} \cap
M_{2\times1}(\mathbb{Z}) \in PM_{2\times 1}(\mathbb{Z})$$ for each
$B \in \textbf{\textit{BT}}$. This is Krebes' tangle invariant.
\end{df}

Let $n$ be a positive integer. Then an $n$-punctured ball tangle
$T^n$ with $(H_k)_{0 \leq k \leq n}$ can be regarded as an
$n$-variable function
$T^n:\textbf{\textit{A}}_1\times\cdots\times\textbf{\textit{A}}_n
\rightarrow \textbf{\textit{T}}$ defined as $T^n(X_1,\dots,X_n)$
is a tangle filled up in the $i$-th hole $H_i$ of $T^n$ by $X_i
\in \textbf{\textit{A}}_i$ for each $i \in \{1,\dots,n\}$, where
$\textbf{\textit{A}}_i$ is a class of $t_i$-punctured ball tangles
for each $i \in \{1,\dots,n\}$ and $\textbf{\textit{T}}$ is a
class of tangles. However, this representation of $n$-punctured
ball tangles as $n$-variable functions is not perfect in the sense
that $n$-punctured ball tangles are equivalent only if they induce
the same function. On the other hand, $n$-punctured ball tangles
which induce the same function need not be equivalent. That is, we
can say that tangles are stronger than functions.

\bigskip
\centerline{\epsfxsize=5 in \epsfbox{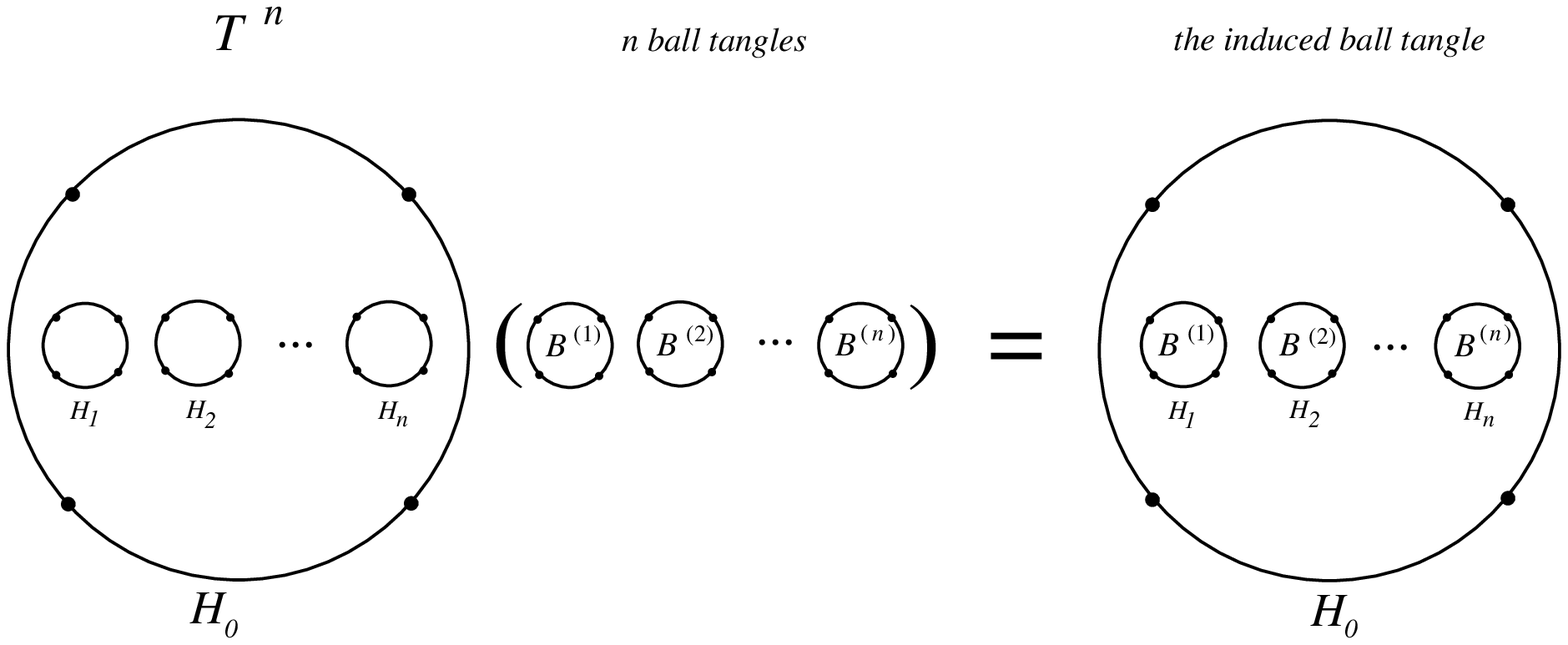}}
\medskip
\centerline{\small Figure 5. The induced ball tangle
$T^n(B^{(1)},\cdots,B^{(n)})$ by $T^n$ and
$B^{(1)},\cdots,B^{(n)}$.}
\bigskip

Roughly speaking, the class of $n$-punctured ball tangles as only
$n$-variable functions gives us an \textit{operad}, a mathematical
device which describes algebraic structure of many varieties and
in various categories. See \cite{MSS}.

First of all, to construct the invariant $F^n$ of $n$-punctured
ball tangle $T^n$, let us regard $T^n$ as a `hole-filling
function', in sense described as above $T^n:\textbf{\textit{BT}}^n
\rightarrow \textbf{\textit{BT}}$, where
$\textbf{\textit{BT}}^n=\textbf{\textit{BT}}_1\times\cdots\times\textbf{\textit{BT}}_n$
with
$\textbf{\textit{BT}}_1=\cdots=\textbf{\textit{BT}}_n=\textbf{\textit{BT}}$
(Figure 5).

Also, to construct our invariant of $n$-punctured ball tangles, we
need to use some quite complicated notations. Let us start with a
gentle introduction to our notations:

(1) For a diagram of 0-punctured ball tangle $T^0$ (a ball
tangle), we can produce 2 links $T_1^0$ and $T_2^0$, which are the
numerator closure and the denominator closure of $T^0$,
respectively.

(2) For a diagram of 1-punctured ball tangle $T^1$ (a spherical
tangle), we can produce $2^{1+1}$ links $T_{1(1)}^1$,
$T_{1(2)}^1$; $T_{2(1)}^1$, $T_{2(2)}^1$, where the subscript 1(1)
means to take the numerator closure of $T$ with its hole filled by
the fundamental tangle 1.

(3) For a diagram of 2-punctured ball tangle $T^2$, we can produce
$2^{2+1}$ links $T_{1(11)}^2$, $T_{1(12)}^2$, $T_{1(21)}^2$,
$T_{1(22)}^2$; $T_{2(11)}^2$, $T_{2(12)}^2$, $T_{2(21)}^2$,
$T_{2(22)}^2$.

If $n$ is a positive integer, $J_1=\cdots=J_n=\{1,2\}$, and
$J(n)=\prod_{k=1}^n J_k$, then $J(n)$ is linearly ordered by a
dictionary order, or lexicographic order, consisting of $2^n$
ordered $n$-tuples each of whose components is either 1 or 2. That
is, if $x,y \in J(n)$ and $x=(x_1,\dots,x_n)$,
$y=(y_1,\dots,y_n)$, then $x<y$ if and only if $x_1<y_1$ or there
is $k \in \{1,\dots,n-1\}$ such that
$x_1=y_1,\dots,x_k=y_k,x_{k+1}<y_{k+1}$.

(4) $J(n)=\{\alpha_i^n|1 \leq i \leq 2^n\}$ and
$\alpha_1^n<\alpha_2^n<\cdots<\alpha_{2^n}^n$, where $<$ is the
dictionary order on $J(n)$. Hence, $\alpha_1^n$ is the least
element $(1,1,\dots,1)$ and $\alpha_{2^n}^n$ is the greatest
element $(2,2,\dots,2)$ of $J(n)$. Let us denote
$\alpha_i^n=(\alpha_{i1}^n,\dots,\alpha_{in}^n)$ for each $i \in
\{1,\dots,2^n\}$.

\bigskip
\centerline{\epsfxsize=5 in \epsfbox{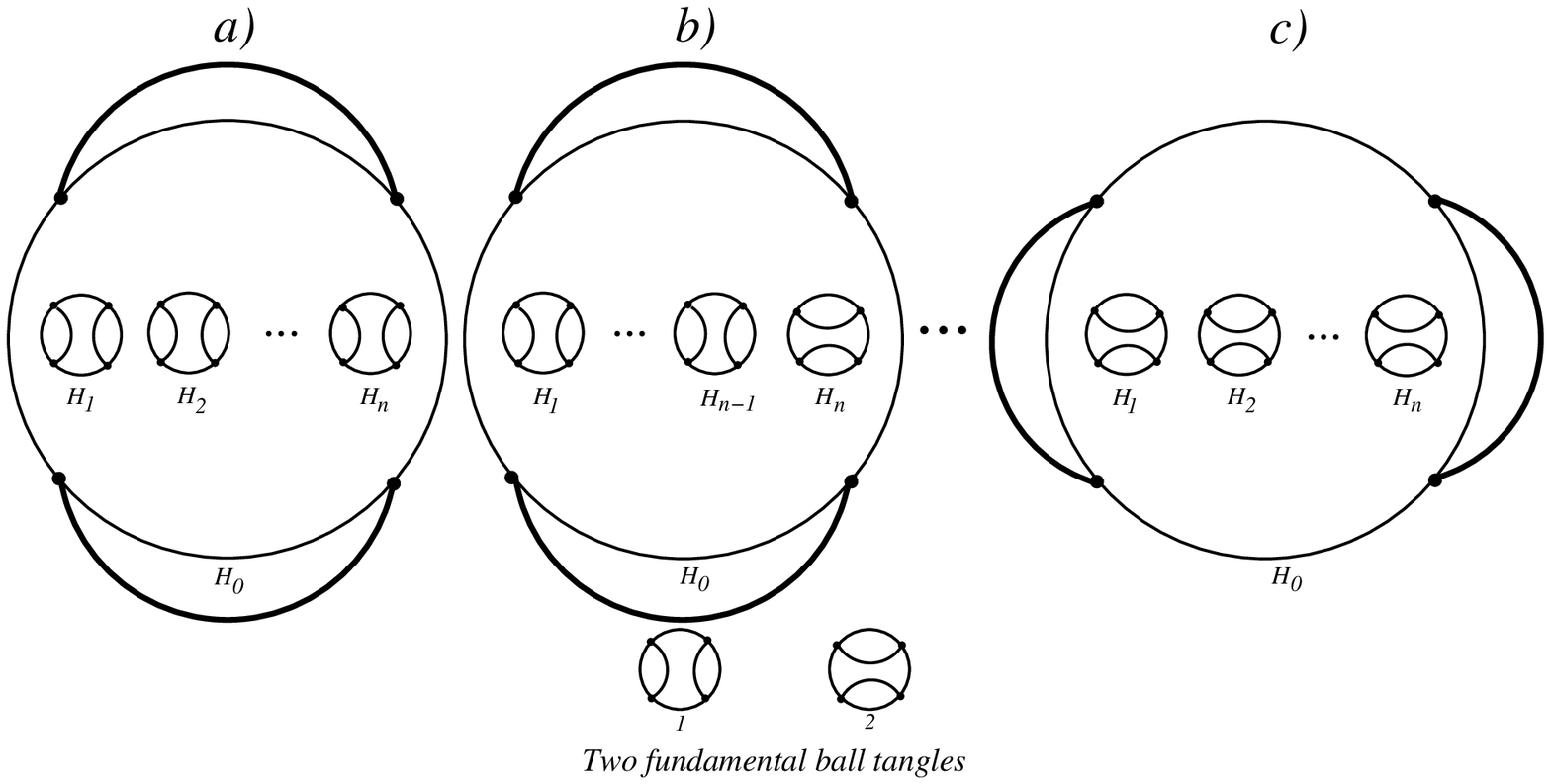}}
\medskip
\centerline{\small Figure 6. The closures of $T^n$, a)
$T_{1\alpha_1^n}^n$, b) $T_{1\alpha_2^n}^n$, c)
$T_{2\alpha_{2^n}^n}^n$.}
\bigskip

(5) For a diagram of $n$-punctured ball tangle $T^n$, we can
produce $2^{n+1}$ links
$T_{1\alpha_1^n}^n,\dots,T_{1\alpha_{2^n}^n}^n$;
$T_{2\alpha_1^n}^n,\dots,T_{2\alpha_{2^n}^n}^n$.

(6) The sequence $(a_n)_{n \geq 0}=((t_k)_{1 \leq k \leq 2^n})_{n
\geq 0}$ is defined recursively as follows:

1) $a_0=(0)$;

2) If $a_{k-1}=(t_1,\dots,t_{2^{k-1}})$, then
$a_k=(t_1,\dots,t_{2^{k-1}},t_1+1,\dots,t_{2^{k-1}}+1)$ for each
$k \in N$. Note that $t_{2^n}=n$ for each $n \in \mathbb N \cup
\{0\}$.

Now, we define our invariant of $n$-punctured ball tangles
inductively.

\begin{thm} For each $n \in \mathbb N$, define $F^n:\textbf{\textit{nPBT}}
\rightarrow PM_{2\times2^n}(\mathbb{Z})$ by
$$F^n(T^n)=\left\{\begin{pmatrix}
(-i)^{t_1}z\langle T_{1\alpha_1^n}^n \rangle &
\cdots & (-i)^{t_{2^n}}z\langle T_{1\alpha_{2^n}^n}^n \rangle \\
(-i)^{t_1}iz\langle T_{2\alpha_1^n}^n \rangle &
\cdots & (-i)^{t_{2^n}}iz\langle T_{2\alpha_{2^n}^n}^n \rangle
\end{pmatrix}\,|\,z \in \Phi\right\}
\cap M_{2\times2^n}(\mathbb Z)$$ for each $T^n \in
\textbf{\textit{nPBT}}$. Then $F^n$ is an isotopy invariant of
$n$-punctured ball tangle diagrams. In particular, $F^0$ is
Krebes' ball tangle invariant $f$.
\end{thm}

\begin{df} For each nonnegative integer $n$, $F^n$ is called the
$n$-punctured ball tangle invariant, simply, the $n$-punctured
tangle invariant.
\end{df}

In order to think of $n$-punctured ball tangle $T^n$ as a
`hole-filling function', we define a function which makes a
dictionary order on complex numbers.

Let $n$ be a positive integer, and let $(k_1,\dots,k_n)$ be an
$n$-tuple of positive integers, and let
$J(n,k_1,\dots,k_n)=\prod_{i=1}^n I_{k_i}$. Then
$J(n,k_1,\dots,k_n)$ is linearly ordered by a dictionary order,
where $I_k=\{1,\dots,k\}$ for each $k \in \mathbb{N}$.

(4$^*$) $J(n,k_1,\dots,k_n)= \{\alpha_i^{n,k_1,\dots,k_n}|1 \leq i
\leq k_1\cdots k_n\}$ and $\alpha_1^{n,k_1,\dots,k_n}<\cdots<
\alpha_{k_1\cdots k_n}^{n,k_1,\dots,k_n}$, where $<$ is the
dictionary order on $J(n,k_1,\dots,k_n)$ and
$\alpha_1^{n,k_1,\dots,k_n}$ is the least element $(1,1,\dots,1)$
and $\alpha_{k_1\cdots k_n}^{n,k_1,\dots,k_n}$ is the greatest
element $(k_1,k_2,\dots,k_n)$ of $J(n,k_1,\dots,k_n)$. Let us
denote $\alpha_i^{n,k_1,\dots,k_n}=(\alpha_{i1}^{n,k_1,\dots,k_n},
\dots,\alpha_{in}^{n,k_1,\dots,k_n})$ for each $i \in
\{1,\dots,k_1\cdots k_n\}$.

\begin{df}
For each $n \in \mathbb{N}$ and $n$-tuple $(k_1,\dots,k_n)$ of
positive integers, define
$$\xi^{n,k_1,\dots,k_n}:\mathbb{C}^{k_1}\times\cdots\times
\mathbb{C}^{k_n} \rightarrow \mathbb{C}^{k_1\cdots k_n}$$ by
$$\xi^{n,k_1,\dots,k_n}((v_1^1,\dots,v_{k_1}^1),
\dots,(v_1^n,\dots,v_{k_n}^n))=\left(\prod_{j=1}^n
v_{\alpha_{1j}^{n,k_1,\dots,k_n}}^j, \dots, \prod_{j=1}^n
v_{\alpha_{k_1\cdots k_nj}^{n,k_1,\dots,k_n}}^j\right)$$ for all
$(v_1^1,\dots,v_{k_1}^1) \in
\mathbb{C}^{k_1},\dots,(v_1^n,\dots,v_{k_n}^n) \in
\mathbb{C}^{k_n}$. Then $\xi^{n,k_1,\dots,k_n}$ is well-defined
and called the dictionary order function on $\mathbb C$ with
respect to $k_1,\dots,k_n$. Also, the $i$-th projection of
$\xi^{n,k_1,\dots,k_n}$ is denoted by $\xi_i^{n,k_1,\dots,k_n}$
for each $i \in \{1,\dots,k_1\cdots k_n\}$. In particular, we
simply denote $\xi^{n,k_1,\dots,k_n}$ by $\xi^n$ when
$k_1=\cdots=k_n=2$.
\end{df}

Denote by $\mathbb{C}^{k\dag}$ the $k$-dimensional column vector
space over $\mathbb{C}$, so the map
$$(v_1,\dots,v_k)\mapsto(v_1,\dots,v_k)^\dag:\mathbb{C}^k\longrightarrow\mathbb{C}^{k\dag}$$
is to transpose row vectors to column vectors. Let
$P\mathbb{C}^{k\dag}=\mathbb{C}^{k\dag}/\pm1$. If
$(v_1,\dots,v_k)^\dag \in \mathbb{C}^{k\dag}$, then we denote by
$$[v_1,\dots,v_k]^\dag=\{(v_1,\dots,v_k)^\dag,(-v_1,\dots,-v_k)^\dag\}$$
the corresponding element in $P\mathbb{C}^{k\dag}$.

Remark that we may extend the above notation to matrices modulo
$\pm1$. Under this extension, matrix multiplication is
well-defined. That is, if $A$ and $B$ are matrices and $AB$ is
defined, then $[A][B]=[A][-B]=[-A][B]=[-A][-B]=[-AB]=[AB]$.

\begin{lm} For each $n \in\mathbb{N}$ and $n$-tuple $(k_1,\dots,k_n)$
of positive integers, define
$$[\xi^{n,k_1,\dots,k_n}]:P\mathbb{C}^{k_1\dag}\times\cdots\times
P\mathbb{C}^{k_n\dag} \longrightarrow P\mathbb{C}^{k_1\cdots
k_n\dag}$$ by $$[\xi^{n,k_1,\dots,k_n}](\left[\begin{matrix} v_1^1
\\ \cdot \\ \cdot \\ \cdot \\ v_{k_1}^1 \end{matrix}\right], \dots,
\left[\begin{matrix} v_1^n \\ \cdot \\ \cdot \\ \cdot \\ v_{k_n}^n
\end{matrix}\right])=\left[\begin{matrix}
\prod_{j=1}^n v_{\alpha_{1j}^{n,k_1,\dots,k_n}}^j
\\ \cdot \\ \cdot \\ \cdot \\
\prod_{j=1}^n v_{\alpha_{k_1\cdots k_nj}^{n,k_1,\dots,k_n}}^j
\end{matrix}\right]$$ for all $(v_1^1,\dots,v_{k_1}^1) \in
\mathbb{C}^{k_1},\dots,(v_1^n,\dots,v_{k_n}^n) \in
\mathbb{C}^{k_n}$. Then $[\xi^{n,k_1,\dots,k_n}]$ is well-defined
and called the dictionary order function induced by
$\xi^{n,k_1,\dots,k_n}$.
\end{lm}

As another notation, if $L$ is a link diagram and $T^n$ is a
diagram of $n$-punctured ball tangle for some $n \in\mathbb{ N}
\cup \{0\}$, then the sets of all crossings of $L$ and $T^n$ are
denoted by $c(L)$ and $c(T^n)$, respectively.

\begin{lm} If $n \in \mathbb{N}$ and $T^n$ is an $n$-punctured ball tangle
diagram and $B^{(1)},\dots,B^{(n)}$ are ball tangle diagrams, then
$$\begin{pmatrix} \langle T^n(B^{(1)},\dots,B^{(n)})_1 \rangle \\
\langle T^n(B^{(1)},\dots,B^{(n)})_2 \rangle
\end{pmatrix}=\begin{pmatrix} \sum_{i=1}^{2^n}
\langle T_{1\alpha_i^n}^n \rangle \langle B_{\alpha_{i1}^n}^{(1)}
\rangle \cdots \langle B_{\alpha_{in}^n}^{(n)} \rangle \\
\sum_{i=1}^{2^n} \langle T_{2\alpha_i^n}^n \rangle \langle
B_{\alpha_{i1}^n}^{(1)} \rangle \cdots \langle
B_{\alpha_{in}^n}^{(n)} \rangle
\end{pmatrix}.$$
\end{lm}

\begin{thm} [J.-W. Chung and X.-S. Lin \cite{C-L}]
For each $n \in \mathbb{N}$, $F^n$ is an $n$-punctured ball tangle
invariant such that
$$F^0(T^n(B^{(1)},\dots,B^{(n)}))=F^n(T^n)[\xi^n](F^0(B^{(1)}),
\dots,F^0(B^{(n)}))$$ for all $B^{(1)},\dots,B^{(n)} \in
\textbf{\textit{BT}}$.
\end{thm}

\begin{proof} Suppose that $T^n$ is an $n$-punctured ball tangle
such that $F^n(T^n)=[zX(T^n)]$ for some $z \in \Phi$ and
$B^{(1)},\dots,B^{(n)}$ are ball tangles such that
$$F^0(B^{(1)})=\left[\begin{matrix} z_1\langle B_1^{(1)} \rangle \\
iz_1\langle B_2^{(1)} \rangle \end{matrix}\right],\dots,
F^0(B^{(n)})=\left[\begin{matrix} z_n\langle B_1^{(n)} \rangle \\
iz_n\langle B_2^{(n)} \rangle \end{matrix}\right]$$ for some
$z_1,\dots,z_n \in \Phi$, where $\langle B_1^{(i)} \rangle$ and
$\langle B_2^{(i)} \rangle$ are the numerator closure and the
denominator closure of $B^{(i)}$, respectively, for each $i \in
\{1,\dots,n\}$. Then
$$
\begin{aligned}
&F^n(T^n)[\xi^n](F^0(B^{(1)}),\dots,F^0(B^{(n)}))\\
&=\left[\begin{matrix} \begin{pmatrix} (-i)^{t_1}z\langle
T_{1\alpha_1^n}^n \rangle & \cdots &
(-i)^{t_{2^n}}z\langle T_{1\alpha_{2^n}^n}^n \rangle \\
(-i)^{t_1}iz\langle T_{2\alpha_1^n}^n \rangle & \cdots &
(-i)^{t_{2^n}}iz\langle T_{2\alpha_{2^n}^n}^n \rangle
\end{pmatrix}
\begin{pmatrix} i^{t_1}z_1\cdots z_n\langle B_{\alpha_{11}^n}^{(1)} \rangle
\cdots \langle B_{\alpha_{1n}^n}^{(n)} \rangle \\ \cdot \\ \cdot \\ \cdot \\
i^{t_{2^n}}z_1\cdots z_n\langle B_{\alpha_{2^n1}^n}^{(1)} \rangle
\cdots \langle B_{\alpha_{2^nn}^n}^{(n)} \rangle
\end{pmatrix} \end{matrix}\right]
\end{aligned}$$
$$
\begin{aligned}
&=\left[\begin{matrix} zz_1\cdots z_n(\langle T_{1\alpha_1^n}^n
\rangle \langle B_{\alpha_{11}^n}^{(1)} \rangle \cdots \langle
B_{\alpha_{1n}^n}^{(n)} \rangle + \cdots\cdots + \langle
T_{1\alpha_{2^n}^n}^n \rangle \langle B_{\alpha_{2^n1}^n}^{(1)}
\rangle \cdots \langle B_{\alpha_{2^nn}^n}^{(n)} \rangle) \\
izz_1\cdots z_n(\langle T_{2\alpha_1^n}^n \rangle \langle
B_{\alpha_{11}^n}^{(1)} \rangle \cdots \langle
B_{\alpha_{1n}^n}^{(n)} \rangle + \cdots\cdots + \langle
T_{2\alpha_{2^n}^n}^n \rangle \langle B_{\alpha_{2^n1}^n}^{(1)}
\rangle \cdots \langle B_{\alpha_{2^nn}^n}^{(n)} \rangle)
\end{matrix}\right]\\
&=\left[\begin{matrix} zz_1\cdots z_n \sum_{i=1}^{2^n} \langle
T_{1\alpha_i^n}^n \rangle \langle B_{\alpha_{i1}^n}^{(1)}
\rangle \cdots \langle B_{\alpha_{in}^n}^{(n)} \rangle \\
izz_1\cdots z_n \sum_{i=1}^{2^n} \langle T_{2\alpha_i^n}^n \rangle
\langle B_{\alpha_{i1}^n}^{(1)} \rangle \cdots \langle
B_{\alpha_{in}^n}^{(n)} \rangle
\end{matrix}\right]=\left[\begin{matrix} zz_1\cdots z_n
\langle T^n(B^{(1)},\dots,B^{(n)})_1 \rangle \\
izz_1\cdots z_n\langle T^n(B^{(1)},\dots,B^{(n)})_2 \rangle
\end{matrix}\right]\\
&=F^0(T^n(B^{(1)},\dots,B^{(n)}))
\end{aligned}$$
by Lemma 2.13.
\end{proof}

\section{Generalized formulas for invariant of $n$-punctured ball tangles}

Notice that an $n$-punctured ball tangle $T^n$ may be regarded as
an $n$ variable function about not only $0$-punctured ball tangles
but also various punctured ball tangles. Given an $n$-punctured
ball tangle diagram $T^n$ and $k_1,\dots,k_n$-punctured ball
tangle diagrams $T^{k_1(1)},\dots,T^{k_n(n)}$, respectively, we
consider the induced $k_1 + \cdots + k_n$-punctured ball tangle
diagram $T^n(T^{k_1(1)},\dots,T^{k_n(n)})$, where $n \in \mathbb
N$ and $k_1,\dots,k_n \in \mathbb N \cup \{0\}$. We show how to
calculate the invariant
$F^{k_1+\cdots+k_n}(T^n(T^{k_1(1)},\dots,T^{k_n(n)}))$ of it if
$F^n(T^n), F^{k_1}(T^{k_1(1)}),\dots,F^{k_n}(T^{k_n(n)})$ are
given (Theorem 3.2). On the other hand, we consider the horizontal
connect sum $T^{k_1(1)} +_h T^{k_2(2)}$ and the vertical connect
sum $T^{k_1(1)} +_v T^{k_2(2)}$ of $k_1$ and $k_2$-punctured ball
tangles $T^{k_1(1)}$ and $T^{k_2(2)}$, respectively, and provide
the formulas for the invariants $F^{k_1+k_2}(T^{k_1(1)} +_h
T^{k_2(2)})$ and $F^{k_1+k_2}(T^{k_1(1)} +_v T^{k_2(2)})$ from
$F^{k_1}(T^{k_1(1)})$ and $F^{k_2}(T^{k_2(2)})$ (Theorem 3.3). To
prove these two generalized formulas, we require a statement from
`Projective Linear Algebra' (Lemma 3.1). Let us start from the
following notations:

Let $n \in \mathbb N$. Then

(1) $e_i^n=\left[\begin{matrix} v_1 \\ \cdot \\ \cdot \\ \cdot \\
v_{2^n} \end{matrix}\right]$ such that $v_i=1$ and $v_j=0$ if $j
\neq i$ for each $i \in \{1,\dots,2^n\}$. In particular,
$e_1^1=\left[\begin{matrix} 1 \\ 0
\end{matrix}\right]$ and $e_2^1=\left[\begin{matrix} 0 \\ 1
\end{matrix}\right]$. Hence, $e_i^n=[\xi^n]
(e_{\alpha_{i1}^n}^1,\dots,e_{\alpha_{in}^n}^1)$ for each $i \in
\{1,\dots,2^n\}$.

(2) $x=\left[\begin{matrix} 1 \\ 1
\end{matrix}\right]$.

(3) $E_j^n$ is the set of all $[\xi^n](y_1,\dots,y_n)$ such that
$j$ components of $(y_1,\dots,y_n)$ are $x$ and each of the others
is $e_1^1$ or $e_2^1$ for each $j \in \{0,1,\dots,n\}$. In
particular,
$$E_0^n=\{[\xi^n](e_{\alpha_{i1}^n}^1,\dots,e_{\alpha_{in}^n}^1)|i
\in \{1,\dots,2^n\}\}$$ and
$$E_n^n=\{[\xi^n](y_1,\dots,y_n)|y_1=\cdots=y_n=x\}.$$
Notice that $\{E_0^n,E_1^n,\dots,E_n^n\}$ is pairwise disjoint and
$|E_j^n|=\,_nC_j\,2^{n-j}$ for each $j \in \{0,1,\dots,n\}$, where
$\,_nC_j\,=\frac{n!}{(n-j)!j\,!}$. Hence,
$$|\coprod_{j=0}^n E_j^n|=\,_nC_0\,2^n + \,_nC_1\,2^{n-1} + \cdots
+ \,_nC_{n-1}\,2^1 + \,_nC_n\,2^0=(2+1)^n=3^n.$$ Note that
$$\coprod_{j=0}^n E_j^n=[\xi^n](\{e_1^1,e_2^1,x\}^n).$$

For example, when $n=3$, we have
$$E_0^3=\{e_1^3,e_2^3,e_3^3,e_4^3,e_5^3,e_6^3,e_7^3,e_8^3\},$$
$$E_1^3=\{\,[\,1\,0\,0\,0\,1\,0\,0\,0\,]^\dag,[\,0\,1\,0\,0\,0\,1\,0\,0\,]^\dag,
[\,0\,0\,1\,0\,0\,0\,1\,0\,]^\dag,[\,0\,0\,0\,1\,0\,0\,0\,1\,]^\dag,$$
$$[\,1\,0\,1\,0\,0\,0\,0\,0\,]^\dag,[\,0\,1\,0\,1\,0\,0\,0\,0\,]^\dag,
[\,0\,0\,0\,0\,1\,0\,1\,0\,]^\dag,[\,0\,0\,0\,0\,0\,1\,0\,1\,]^\dag,$$
$$[\,1\,1\,0\,0\,0\,0\,0\,0\,]^\dag,[\,0\,0\,1\,1\,0\,0\,0\,0\,]^\dag,
[\,0\,0\,0\,0\,1\,1\,0\,0\,]^\dag,[\,0\,0\,0\,0\,0\,0\,1\,1\,]^\dag\},$$
$$E_2^3=\{\,[\,1\,0\,1\,0\,1\,0\,1\,0\,]^\dag,[\,0\,1\,0\,1\,0\,1\,0\,1\,]^\dag,$$
$$[\,1\,1\,0\,0\,1\,1\,0\,0\,]^\dag,[\,0\,0\,1\,1\,0\,0\,1\,1\,]^\dag,$$
$$[\,1\,1\,1\,1\,0\,0\,0\,0\,]^\dag,[\,0\,0\,0\,0\,1\,1\,1\,1\,]^\dag\},$$
$$E_3^3=\{\,[\,1\,1\,1\,1\,1\,1\,1\,1\,]^\dag\}.$$

Now, we have the following lemma which supports our main theorems.

\begin{lm} If $n \in \mathbb N$ and $A,B \in PM_{2\times2^n}(\mathbb Z)$ and
$AX=BX$ for each $X \in [\xi^n](\{e_1^1,e_2^1,x\}^n)$, then $A=B$.
\end{lm}

\begin{proof} We prove the statement by induction on $n \in \mathbb N$.

Step 1. We show that the statement is true for $n=1$.

Let $A=\left[\begin{matrix} a_{11} & a_{12} \\ a_{21} & a_{22}
\end{matrix}\right]$, and let $B=\left[\begin{matrix}
b_{11} & b_{12} \\ b_{21} & b_{22} \end{matrix}\right]$. Since
$E_0^1=\{\left[\begin{matrix} 1 \\ 0
\end{matrix}\right],\left[\begin{matrix} 0 \\ 1
\end{matrix}\right]\}$ and $E_1^1=\{\left[\begin{matrix} 1 \\ 1
\end{matrix}\right]\}$ and $AX=BX$ for each $X \in [\xi^1](\{e_1^1,e_2^1,x\}^1)$,
$\left[\begin{matrix} a_{11} \\ a_{21}
\end{matrix}\right]=\left[\begin{matrix} b_{11} \\ b_{21}
\end{matrix}\right]$ and $\left[\begin{matrix} a_{12} \\ a_{22}
\end{matrix}\right]=\left[\begin{matrix} b_{12} \\ b_{22}
\end{matrix}\right]$ and $\left[\begin{matrix} a_{11}+a_{12} \\
a_{21}+a_{22} \end{matrix}\right]=\left[\begin{matrix} b_{11}+b_{12} \\
b_{21}+b_{22} \end{matrix}\right]$. Hence, $\begin{pmatrix}
a_{11} \\ a_{21} \end{pmatrix}=\epsilon\begin{pmatrix} b_{11} \\
b_{21} \end{pmatrix}$ and $\begin{pmatrix}
a_{12} \\ a_{22} \end{pmatrix}=\epsilon'\begin{pmatrix} b_{12} \\
b_{22} \end{pmatrix}$ and $\begin{pmatrix} a_{11}+a_{12} \\
a_{21}+a_{22} \end{pmatrix}= \epsilon_1\begin{pmatrix}
b_{11}+b_{12} \\ b_{21}+b_{22} \end{pmatrix}$ for some
$\epsilon,\epsilon',\epsilon_1 \in \{1,-1\}$. Suppose that
$\epsilon\epsilon'=-1$.

{\it Case 1}. $\begin{pmatrix} a_{11} \\ a_{21}
\end{pmatrix}=\begin{pmatrix} b_{11} \\ b_{21}
\end{pmatrix}$ and $\begin{pmatrix} a_{12} \\ a_{22}
\end{pmatrix}=\begin{pmatrix} -b_{12} \\ -b_{22} \end{pmatrix}$.
If $\epsilon_1=1$, then
$a_{11}+a_{12}=b_{11}+b_{12}=b_{11}-b_{12}$ and
$a_{21}+a_{22}=b_{21}+b_{22}=b_{21}-b_{22}$, so $b_{12}=b_{22}=0$.
If $\epsilon_1=-1$, then
$a_{11}+a_{12}=-b_{11}-b_{12}=b_{11}-b_{12}$ and
$a_{21}+a_{22}=-b_{21}-b_{22}=b_{21}-b_{22}$, so
$b_{11}=b_{21}=0$. Hence, $A=B$.

{\it Case 2}. $\begin{pmatrix} a_{11} \\ a_{21}
\end{pmatrix}=\begin{pmatrix} -b_{11} \\ -b_{21}
\end{pmatrix}$ and $\begin{pmatrix} a_{12} \\ a_{22}
\end{pmatrix}=\begin{pmatrix} b_{12} \\ b_{22} \end{pmatrix}$.
If $\epsilon_1=1$, then
$a_{11}+a_{12}=b_{11}+b_{12}=-b_{11}+b_{12}$ and
$a_{21}+a_{22}=b_{21}+b_{22}=-b_{21}+b_{22}$, so
$b_{11}=b_{21}=0$. If $\epsilon_1=-1$, then
$a_{11}+a_{12}=-b_{11}-b_{12}=-b_{11}+b_{12}$ and
$a_{21}+a_{22}=-b_{21}-b_{22}=-b_{21}+b_{22}$, so
$b_{12}=b_{22}=0$. Hence, $A=B$.

Step 2. Suppose that the statement is true for $n \in \mathbb N$.
We show that the statement is also true for $n+1$.

Suppose that $$A=\left[\begin{matrix} a_{11} & \cdots & a_{12^n} &
a_{12^n+1} & \cdots & a_{12^{n+1}} \\ a_{21} & \cdots & a_{22^n} &
a_{22^n+1} & \cdots & a_{22^{n+1}} \end{matrix}\right]$$ and
$$B=\left[\begin{matrix} b_{11} & \cdots & b_{12^n} & b_{12^n+1} &
\cdots & b_{12^{n+1}} \\ b_{21} & \cdots & b_{22^n} & b_{22^n+1} &
\cdots & b_{22^{n+1}} \end{matrix}\right]$$ and
$$A_1=\begin{pmatrix} a_{11} & \cdots & a_{12^n} \\
a_{21} & \cdots & a_{22^n} \end{pmatrix},\,
A_2=\begin{pmatrix} a_{12^n+1} & \cdots & a_{12^{n+1}} \\
a_{22^n+1} & \cdots & a_{22^{n+1}} \end{pmatrix},$$
$$A_3=\begin{pmatrix}
a_{11} & \cdots & a_{12^{n-1}} & a_{12^n+1} & \cdots &
a_{12^n+2^{n-1}} \\ a_{21} & \cdots & a_{22^{n-1}} & a_{22^n+1} &
\cdots & a_{22^n+2^{n-1}} \end{pmatrix},$$
$$A_4=\begin{pmatrix}
a_{12^{n-1}+1} & \cdots & a_{12^n} & a_{12^n+2^{n-1}+1} & \cdots &
a_{12^{n+1}} \\ a_{22^{n-1}+1} & \cdots & a_{22^n} &
a_{22^n+2^{n-1}+1} & \cdots & a_{22^{n+1}}
\end{pmatrix}$$ and $$B_1=\begin{pmatrix} b_{11} & \cdots & b_{12^n} \\
b_{21} & \cdots & b_{22^n} \end{pmatrix},\,
B_2=\begin{pmatrix} b_{12^n+1} & \cdots & b_{12^{n+1}} \\
b_{22^n+1} & \cdots & b_{22^{n+1}} \end{pmatrix},$$
$$B_3=\begin{pmatrix}
b_{11} & \cdots & b_{12^{n-1}} & b_{12^n+1} & \cdots &
b_{12^n+2^{n-1}} \\ b_{21} & \cdots & b_{22^{n-1}} & b_{22^n+1} &
\cdots & b_{22^n+2^{n-1}} \end{pmatrix},$$
$$B_4=\begin{pmatrix}
b_{12^{n-1}+1} & \cdots & b_{12^n} & b_{12^n+2^{n-1}+1} & \cdots &
b_{12^{n+1}} \\ b_{22^{n-1}+1} & \cdots & b_{22^n} &
b_{22^n+2^{n-1}+1} & \cdots & b_{22^{n+1}}
\end{pmatrix}.$$

Then $A=\left[\begin{matrix} A_1 & A_2
\end{matrix}\right]$ and $B=\left[\begin{matrix} B_1 & B_2
\end{matrix}\right]$. Notice that $$[\xi^{n+1}](\{e_1^1,e_2^1,x\}^{n+1})=$$
$$[\xi^{n+1}](\{e_1^1\} \times \{e_1^1,e_2^1,x\}^n) \coprod
[\xi^{n+1}](\{e_2^1\} \times \{e_1^1,e_2^1,x\}^n) \coprod
[\xi^{n+1}](\{x\} \times \{e_1^1,e_2^1,x\}^n)$$ and
$[\xi^{n+1}](\{e_1^1\} \times \{e_1^1,e_2^1,x\}^n)$,
$[\xi^{n+1}](\{e_2^1\} \times \{e_1^1,e_2^1,x\}^n)$,
$[\xi^{n+1}](\{x\} \times \{e_1^1,e_2^1,x\}^n)$ have exactly $3^n$
elements, respectively.

Since $AX=BX$ for each $X \in [\xi^{n+1}](\{e_1^1\} \times
\{e_1^1,e_2^1,x\}^n)$, $[A_1]X=[B_1]X$ for each $X \in
[\xi^n](\{e_1^1,e_2^1,x\}^n)$.

Similarly, since $AX=BX$ for each $X \in [\xi^{n+1}](\{e_2^1\}
\times \{e_1^1,e_2^1,x\}^n)$, $[A_2]X=[B_2]X$ for each $X \in
[\xi^n](\{e_1^1,e_2^1,x\}^n)$.

Also, since $AX=BX$ for each $X \in [\xi^{n+1}](\{e_1^1,e_2^1,x\}
\times \{e_1^1\} \times \{e_1^1,e_2^1,x\}^{n-1})$, $[A_3]X=[B_3]X$
for each $X \in [\xi^n](\{e_1^1,e_2^1,x\}^n)$.

Similarly, since $AX=BX$ for each $X \in
[\xi^{n+1}](\{e_1^1,e_2^1,x\} \times \{e_2^1\} \times
\{e_1^1,e_2^1,x\}^{n-1})$, $[A_4]X=[B_4]X$ for each $X \in
[\xi^n](\{e_1^1,e_2^1,x\}^n)$.

By induction hypothesis, we have $$[A_1]=[B_1],\,\,
[A_2]=[B_2],\,\, [A_3]=[B_3],\,\, [A_4]=[B_4].$$ Hence, $A_1 =
\epsilon B_1$ and $A_2 = \epsilon' B_2$ for some
$\epsilon,\epsilon' \in \{1,-1\}$. Now, we claim that, if
$\epsilon\epsilon'=-1$, then $B_1$ or $B_2$ is the $2\times2^n$
zero matrix.

Suppose that $\epsilon\epsilon'=-1$. Without loss of generality,
we may assume that $$A_1=B_1\,\, {\rm and}\,\, A_2=-B_2.$$

Suppose that $B_1$ is not the $2\times2^n$ zero matrix. Then there
is $i \in \{1,\dots,2^n\}$ such that $$\begin{pmatrix} b_{1i} \\
b_{2i} \end{pmatrix} \neq \begin{pmatrix} 0 \\ 0 \end{pmatrix}.$$

{\it Case 1}. If $1 \leq i \leq 2^{n-1}$, then
$\begin{pmatrix} b_{12^n+1} & \cdots & b_{12^n+2^{n-1}} \\
b_{22^n+1} & \cdots & b_{22^n+2^{n-1}} \end{pmatrix}$ is the
$2\times2^{n-1}$ zero matrix since $[A_3]=[B_3]$. We claim that
$\begin{pmatrix} b_{12^n+2^{n-1}+1} & \cdots & b_{12^{n+1}} \\
b_{22^n+2^{n-1}+1} & \cdots & b_{22^{n+1}} \end{pmatrix}$ is also
the $2\times2^{n-1}$ zero matrix.

Suppose that $\begin{pmatrix} b_{12^n+2^{n-1}+1} & \cdots & b_{12^{n+1}} \\
b_{22^n+2^{n-1}+1} & \cdots & b_{22^{n+1}} \end{pmatrix}$ is not
the $2\times2^{n-1}$ zero matrix. Then there is $j \in
\{2^n+2^{n-1}+1,\dots,2^{n+1}\}$ such that $\begin{pmatrix} b_{1j}
\\ b_{2j} \end{pmatrix} \neq \begin{pmatrix} 0 \\ 0 \end{pmatrix}$.

Since $[A_4]=[B_4]$, $\begin{pmatrix} b_{12^{n-1}+1} & \cdots &
b_{12^n} \\ b_{22^{n-1}+1} & \cdots & b_{22^n} \end{pmatrix}$ is
the $2\times2^{n-1}$ zero matrix. In this case, the fact that
$AX=BX$ for each $X \in [\xi^{n+1}](\{e_1^1,e_2^1,x\}^{n+1})$
implies
$$\left[\begin{matrix} a_{11} & \cdots & a_{12^{n-1}} &
a_{12^n+2^{n-1}+1} & \cdots & a_{12^{n+1}} \\ a_{21} & \cdots &
a_{22^{n-1}} & a_{22^n+2^{n-1}+1} & \cdots & a_{22^{n+1}}
\end{matrix}\right]X$$
$$=\left[\begin{matrix} b_{11} & \cdots & b_{12^{n-1}} &
b_{12^n+2^{n-1}+1} & \cdots & b_{12^{n+1}} \\ b_{21} & \cdots &
b_{22^{n-1}} & b_{22^n+2^{n-1}+1} & \cdots & b_{22^{n+1}}
\end{matrix}\right]X$$ for each $X \in
[\xi^n](\{e_1^1,e_2^1,x\}^n)$. Hence, by induction hypothesis, we
have $$\left[\begin{matrix} a_{11} & \cdots & a_{12^{n-1}} &
a_{12^n+2^{n-1}+1} & \cdots & a_{12^{n+1}} \\ a_{21} & \cdots &
a_{22^{n-1}} & a_{22^n+2^{n-1}+1} & \cdots & a_{22^{n+1}}
\end{matrix}\right]$$
$$=\left[\begin{matrix} b_{11} & \cdots & b_{12^{n-1}} &
b_{12^n+2^{n-1}+1} & \cdots & b_{12^{n+1}} \\ b_{21} & \cdots &
b_{22^{n-1}} & b_{22^n+2^{n-1}+1} & \cdots & b_{22^{n+1}}
\end{matrix}\right].$$

Since $A_1=B_1$ and $A_2=-B_2$,
$$\left[\begin{matrix} b_{11} & \cdots & b_{12^{n-1}} &
-b_{12^n+2^{n-1}+1} & \cdots & -b_{12^{n+1}} \\ b_{21} & \cdots &
b_{22^{n-1}} & -b_{22^n+2^{n-1}+1} & \cdots & -b_{22^{n+1}}
\end{matrix}\right]$$
$$=\left[\begin{matrix} b_{11} & \cdots & b_{12^{n-1}} &
b_{12^n+2^{n-1}+1} & \cdots & b_{12^{n+1}} \\ b_{21} & \cdots &
b_{22^{n-1}} & b_{22^n+2^{n-1}+1} & \cdots & b_{22^{n+1}}
\end{matrix}\right].$$

Since $\begin{pmatrix} b_{1i} \\ b_{2i} \end{pmatrix} \neq
\begin{pmatrix} 0 \\ 0 \end{pmatrix}$,
$\begin{pmatrix} b_{12^n+2^{n-1}+1} & \cdots & b_{12^{n+1}} \\
b_{22^n+2^{n-1}+1} & \cdots & b_{22^{n+1}} \end{pmatrix}$ is the
$2\times2^{n-1}$ zero matrix. This is a contradiction. Therefore,
$B_2$ is the $2\times2^n$ zero matrix.

Similarly, we show the other case.

{\it Case 2}. If $2^{n-1}+1 \leq i \leq 2^n$, then
$\begin{pmatrix} b_{12^n+2^{n-1}+1} & \cdots & b_{12^{n+1}} \\
b_{22^n+2^{n-1}+1} & \cdots & b_{22^{n+1}} \end{pmatrix}$ is the
$2\times2^{n-1}$ zero matrix since $[A_4]=[B_4]$. We claim that
$\begin{pmatrix} b_{12^n+1} & \cdots & b_{12^n+2^{n-1}} \\
b_{22^n+1} & \cdots & b_{22^n+2^{n-1}} \end{pmatrix}$ is also the
$2\times2^{n-1}$ zero matrix.

Suppose that $\begin{pmatrix} b_{12^n+1} & \cdots & b_{12^n+2^{n-1}} \\
b_{22^n+1} & \cdots & b_{22^n+2^{n-1}} \end{pmatrix}$ is not the
$2\times2^{n-1}$ zero matrix. Then there is $j \in
\{2^n+1,\dots,2^n+2^{n-1}\}$ such that $\begin{pmatrix} b_{1j}
\\ b_{2j} \end{pmatrix} \neq \begin{pmatrix} 0 \\ 0 \end{pmatrix}$.

Since $[A_3]=[B_3]$, $\begin{pmatrix} b_{11} & \cdots &
b_{12^{n-1}} \\ b_{21} & \cdots & b_{22^{n-1}} \end{pmatrix}$ is
the $2\times2^{n-1}$ zero matrix. In this case, the fact that
$AX=BX$ for each $X \in [\xi^{n+1}](\{e_1^1,e_2^1,x\}^{n+1})$
implies
$$\left[\begin{matrix} a_{12^{n-1}+1} & \cdots & a_{12^n} &
a_{12^n+1} & \cdots & a_{12^n+2^{n-1}} \\ a_{22^{n-1}+1} & \cdots
& a_{22^n} & a_{22^n+1} & \cdots & a_{22^n+2^{n-1}}
\end{matrix}\right]X$$
$$=\left[\begin{matrix} b_{12^{n-1}+1} & \cdots & b_{12^n} &
b_{12^n+1} & \cdots & b_{12^n+2^{n-1}} \\ b_{22^{n-1}+1} & \cdots
& b_{22^n} & b_{22^n+1} & \cdots & b_{22^n+2^{n-1}}
\end{matrix}\right]X$$ for each $X \in [\xi^n](\{e_1^1,e_2^1,x\}^n)$. Hence, by
induction hypothesis and $A_1=B_1$ and $A_2=-B_2$, we have
$$\left[\begin{matrix} b_{12^{n-1}+1} & \cdots & b_{12^n} &
-b_{12^n+1} & \cdots & -b_{12^n+2^{n-1}} \\ b_{22^{n-1}+1} &
\cdots & b_{22^n} & -b_{22^n+1} & \cdots & -b_{22^n+2^{n-1}}
\end{matrix}\right]$$
$$=\left[\begin{matrix} b_{12^{n-1}+1} & \cdots & b_{12^n} &
b_{12^n+1} & \cdots & b_{12^n+2^{n-1}} \\ b_{22^{n-1}+1} & \cdots
& b_{22^n} & b_{22^n+1} & \cdots & b_{22^n+2^{n-1}}
\end{matrix}\right].$$
Since $\begin{pmatrix} b_{1i} \\ b_{2i} \end{pmatrix} \neq
\begin{pmatrix} 0 \\ 0 \end{pmatrix}$,
$\begin{pmatrix} b_{12^n+1} & \cdots & b_{12^n+2^{n-1}} \\
b_{22^n+1} & \cdots & b_{22^n+2^{n-1}} \end{pmatrix}$ is the
$2\times2^{n-1}$ zero matrix. This is a contradiction. Therefore,
$B_2$ is the $2\times2^n$ zero matrix.

Hence, for each case, we have $A=\left[\begin{matrix} A_1 & A_2
\end{matrix}\right]=\left[\begin{matrix} B_1 & B_2
\end{matrix}\right]=B$.
This proves the lemma.
\end{proof}

Remark that the invariant of an $n$-punctured ball tangle is a $2
\times 2^n$ `projective matrix' which means a matrix in
$PM_{2\times2^n}(\mathbb Z)$. To prove Theorem 3.2, we will show
that the projective matrices send each of all possible `projective
column vectors' coming from ball tangle invariants the same value.
Fortunately, there are ball tangle diagrams whose invariants are
$\left[\begin{matrix} 1 \\ 0
\end{matrix}\right],\left[\begin{matrix} 0 \\ 1
\end{matrix}\right],\left[\begin{matrix} 1 \\ 1
\end{matrix}\right]$, respectively (See Figure 8).

From this fact, we can say that $n$-punctured ball tangles have
the same invariant if they are the same function on
$\textbf{\textit{BT}}^n$.

\begin{thm} Let $n \in \mathbb N$, and let $k_1,\dots,k_n \in \mathbb N \cup
\{0\}$, and let $T^n,T^{k_1(1)},\dots,T^{k_n(n)}$ be
$n,k_1,\dots,k_n$-punctured ball tangle diagrams, respectively.
Then

if $F^{k_1}(T^{k_1(1)})=\left[\begin{matrix} b_{11}^1 & \cdots &
b_{12^{k_1}}^1 \\ b_{21}^1 & \cdots & b_{22^{k_1}}^1
\end{matrix}\right],\dots,
F^{k_n}(T^{k_n(n)})=\left[\begin{matrix} b_{11}^n & \cdots &
b_{12^{k_n}}^n \\ b_{21}^n & \cdots & b_{22^{k_n}}^n
\end{matrix}\right]$, then

$F^{k_1+\cdots+k_n}(T^n(T^{k_1(1)},\dots,T^{k_n(n)}))=
F^n(T^n)[\eta^n](F^{k_1}(T^{k_1(1)}),\dots,F^{k_n}(T^{k_n(n)}))$,

where $[\eta^n](F^{k_1}(T^{k_1(1)}),\dots,F^{k_n}(T^{k_n(n)}))$

$=\left[\begin{matrix}
\xi^{n,2^{k_1},\dots,2^{k_n}}((b_{\alpha_{11}^n1}^1,\dots,
b_{\alpha_{11}^n2^{k_1}}^1),\dots,(b_{\alpha_{1n}^n1}^n,
\dots,b_{\alpha_{1n}^n2^{k_n}}^n)) \\
\xi^{n,2^{k_1},\dots,2^{k_n}}((b_{\alpha_{21}^n1}^1,\dots,
b_{\alpha_{21}^n2^{k_1}}^1),\dots,(b_{\alpha_{2n}^n1}^n,
\dots,b_{\alpha_{2n}^n2^{k_n}}^n)) \\ \cdot \\ \cdot \\ \cdot \\
\xi^{n,2^{k_1},\dots,2^{k_n}}((b_{\alpha_{2^n1}^n1}^1,\dots,
b_{\alpha_{2^n1}^n2^{k_1}}^1),\dots,(b_{\alpha_{2^nn}^n1}^n,
\dots,b_{\alpha_{2^nn}^n2^{k_n}}^n))
\end{matrix}\right]$

$=\left[\begin{matrix} \prod_{j=1}^n
b_{\alpha_{1j}^n\alpha_{1j}^{n,2^{k_1},\dots,2^{k_n}}}^j &
\prod_{j=1}^n
b_{\alpha_{1j}^n\alpha_{2j}^{n,2^{k_1},\dots,2^{k_n}}}^j & \cdots
& \prod_{j=1}^n b_{\alpha_{1j}^n\alpha_{2^{k_1+\cdots+
k_n}j}^{n,2^{k_1},\dots,2^{k_n}}}^j
\\ \prod_{j=1}^n
b_{\alpha_{2j}^n\alpha_{1j}^{n,2^{k_1},\dots,2^{k_n}}}^j &
\prod_{j=1}^n
b_{\alpha_{2j}^n\alpha_{2j}^{n,2^{k_1},\dots,2^{k_n}}}^j & \cdots
& \prod_{j=1}^n b_{\alpha_{2j}^n\alpha_{2^{k_1+\cdots+
k_n}j}^{n,2^{k_1},\dots,2^{k_n}}}^j
\\ \cdot & \cdot & \cdot & \cdot \\ \cdot & \cdot & \cdot & \cdot
\\ \cdot & \cdot & \cdot & \cdot \\
\prod_{j=1}^n
b_{\alpha_{2^nj}^n\alpha_{1j}^{n,2^{k_1},\dots,2^{k_n}}}^j &
\prod_{j=1}^n
b_{\alpha_{2^nj}^n\alpha_{2j}^{n,2^{k_1},\dots,2^{k_n}}}^j &
\cdots & \prod_{j=1}^n b_{\alpha_{2^nj}^n\alpha_{2^{k_1+\cdots+
k_n}j}^{n,2^{k_1},\dots,2^{k_n}}}^j
\end{matrix}\right]$.
\end{thm}

\begin{proof} Without loss of generality, we may assume that
$k_1,\dots,k_n \in \mathbb N$.

Let $T=T^n(T^{k_1(1)},\dots,T^{k_n(n)})$, and let
$B^{(11)},\dots,B^{(1k_1)},\dots\dots,B^{(n1)},\dots,B^{(nk_n)}
\in \textbf{\textit{BT}}$ with
$$F^0(B^{(11)})=\left[\begin{matrix} v_1^{11} \\ v_2^{11} \end{matrix}\right],
\dots, F^0(B^{(1k_1)})=\left[\begin{matrix} v_1^{1k_1} \\
v_2^{1k_1} \end{matrix}\right],$$
$$\dots\dots,$$
$$F^0(B^{(n1)})=\left[\begin{matrix} v_1^{n1} \\ v_2^{n1} \end{matrix}\right],
\dots, F^0(B^{(nk_n)})=\left[\begin{matrix} v_1^{nk_n} \\
v_2^{nk_n} \end{matrix}\right].$$

Then $T(B^{(11)},\dots,B^{(1k_1)},\dots\dots,
B^{(n1)},\dots,B^{(nk_n)})$

$=T^n(T^{k_1(1)}(B^{(11)},\dots,B^{(1k_1)}),\dots,
T^{k_n(n)}(B^{(n1)},\dots,B^{(nk_n)}))$ and

$F^0(T(B^{(11)},\dots,B^{(1k_1)},\dots\dots,
B^{(n1)},\dots,B^{(nk_n)}))$

$=F^0(T^n(T^{k_1(1)}(B^{(11)},\dots,B^{(1k_1)}),\dots,
T^{k_n(n)}(B^{(n1)},\dots,B^{(nk_n)})))$

$=F^n(T^n)[\xi^n](F^0(T^{k_1(1)}(B^{(11)},\dots,B^{(1k_1)})),
\dots,F^0(T^{k_n(n)}(B^{(n1)},\dots,B^{(nk_n)})))$

$=F^n(T^n)[\xi^n](F^{k_1}(T^{k_1(1)})[\xi^{k_1}](F^0(B^{(11)}),\dots,F^0(B^{(1k_1)})),
\dots, F^{k_n}(T^{k_n(n)})[\xi^{k_n}]$

$(F^0(B^{(n1)}),\dots,F^0(B^{(nk_n)})))=$

$F^n(T^n)[\xi^n](\left[\begin{matrix} b_{11}^1 & \cdots &
b_{12^{k_1}}^1 \\ b_{21}^1 & \cdots & b_{22^{k_1}}^1
\end{matrix}\right]
\left[\begin{matrix} \prod_{j=1}^{k_1} v_{\alpha_{1j}^{k_1}}^{1j}
\\ \cdot \\ \cdot \\ \cdot \\
\prod_{j=1}^{k_1} v_{\alpha_{2^{k_1}j}^{k_1}}^{1j}
\end{matrix}\right],
\dots, \left[\begin{matrix} b_{11}^n & \cdots & b_{12^{k_n}}^n
\\ b_{21}^n & \cdots & b_{22^{k_n}}^n
\end{matrix}\right]
\left[\begin{matrix} \prod_{j=1}^{k_n} v_{\alpha_{1j}^{k_n}}^{nj}
\\ \cdot \\ \cdot \\ \cdot \\
\prod_{j=1}^{k_n} v_{\alpha_{2^{k_n}j}^{k_n}}^{nj}
\end{matrix}\right])$

$=F^n(T^n)[\xi^n](\left[\begin{matrix} b_{11}^1 \prod_{j=1}^{k_1}
v_{\alpha_{1j}^{k_1}}^{1j} + \cdots +
b_{12^{k_1}}^1 \prod_{j=1}^{k_1} v_{\alpha_{2^{k_1}j}^{k_1}}^{1j} \\
b_{21}^1 \prod_{j=1}^{k_1} v_{\alpha_{1j}^{k_1}}^{1j} + \cdots +
b_{22^{k_1}}^1 \prod_{j=1}^{k_1} v_{\alpha_{2^{k_1}j}^{k_1}}^{1j}
\end{matrix}\right],\dots,$

$\left[\begin{matrix} b_{11}^n \prod_{j=1}^{k_n}
v_{\alpha_{1j}^{k_n}}^{nj} + \cdots +
b_{12^{k_n}}^n \prod_{j=1}^{k_n} v_{\alpha_{2^{k_n}j}^{k_n}}^{nj} \\
b_{21}^n \prod_{j=1}^{k_n} v_{\alpha_{1j}^{k_n}}^{nj} + \cdots +
b_{22^{k_n}}^n \prod_{j=1}^{k_n} v_{\alpha_{2^{k_n}j}^{k_n}}^{nj}
\end{matrix}\right])$

$=F^n(T^n)\left[\begin{matrix}
\xi^{n,2^{k_1},\dots,2^{k_n}}((b_{\alpha_{11}^n1}^1,\dots,
b_{\alpha_{11}^n2^{k_1}}^1),\dots,(b_{\alpha_{1n}^n1}^n,
\dots,b_{\alpha_{1n}^n2^{k_n}}^n)) \\
\xi^{n,2^{k_1},\dots,2^{k_n}}((b_{\alpha_{21}^n1}^1,\dots,
b_{\alpha_{21}^n2^{k_1}}^1),\dots,(b_{\alpha_{2n}^n1}^n,
\dots,b_{\alpha_{2n}^n2^{k_n}}^n)) \\ \cdot \\ \cdot \\ \cdot \\
\xi^{n,2^{k_1},\dots,2^{k_n}}((b_{\alpha_{2^n1}^n1}^1,\dots,
b_{\alpha_{2^n1}^n2^{k_1}}^1),\dots,(b_{\alpha_{2^nn}^n1}^n,
\dots,b_{\alpha_{2^nn}^n2^{k_n}}^n))
\end{matrix}\right] \times$

$\left[\begin{matrix}
\xi^{n,2^{k_1},\dots,2^{k_n}}((\prod_{j=1}^{k_1}
v_{\alpha_{1j}^{k_1}}^{1j},\dots,\prod_{j=1}^{k_1}
v_{\alpha_{2^{k_1}j}^{k_1}}^{1j}),\dots,(\prod_{j=1}^{k_n}
v_{\alpha_{1j}^{k_n}}^{nj},\dots,\prod_{j=1}^{k_n}
v_{\alpha_{2^{k_n}j}^{k_n}}^{nj}))
\end{matrix}\right]^\dag$

$=F^n(T^n)\left[\begin{matrix}
\xi^{n,2^{k_1},\dots,2^{k_n}}((b_{\alpha_{11}^n1}^1,\dots,
b_{\alpha_{11}^n2^{k_1}}^1),\dots,(b_{\alpha_{1n}^n1}^n,
\dots,b_{\alpha_{1n}^n2^{k_n}}^n)) \\
\xi^{n,2^{k_1},\dots,2^{k_n}}((b_{\alpha_{21}^n1}^1,\dots,
b_{\alpha_{21}^n2^{k_1}}^1),\dots,(b_{\alpha_{2n}^n1}^n,
\dots,b_{\alpha_{2n}^n2^{k_n}}^n)) \\ \cdot \\ \cdot \\ \cdot \\
\xi^{n,2^{k_1},\dots,2^{k_n}}((b_{\alpha_{2^n1}^n1}^1,\dots,
b_{\alpha_{2^n1}^n2^{k_1}}^1),\dots,(b_{\alpha_{2^nn}^n1}^n,
\dots,b_{\alpha_{2^nn}^n2^{k_n}}^n))
\end{matrix}\right] \times$

$[\xi^{n,2^{k_1},\dots,2^{k_n}}](\left[\begin{matrix}
\prod_{j=1}^{k_1} v_{\alpha_{1j}^{k_1}}^{1j}
\\ \cdot \\ \cdot \\ \cdot \\
\prod_{j=1}^{k_1} v_{\alpha_{2^{k_1}j}^{k_1}}^{1j}
\end{matrix}\right],\dots,\left[\begin{matrix} \prod_{j=1}^{k_n} v_{\alpha_{1j}^{k_n}}^{nj}
\\ \cdot \\ \cdot \\ \cdot \\
\prod_{j=1}^{k_n} v_{\alpha_{2^{k_n}j}^{k_n}}^{nj}
\end{matrix}\right])$

$=F^n(T^n)\left[\begin{matrix}
\xi^{n,2^{k_1},\dots,2^{k_n}}((b_{\alpha_{11}^n1}^1,\dots,
b_{\alpha_{11}^n2^{k_1}}^1),\dots,(b_{\alpha_{1n}^n1}^n,
\dots,b_{\alpha_{1n}^n2^{k_n}}^n)) \\
\xi^{n,2^{k_1},\dots,2^{k_n}}((b_{\alpha_{21}^n1}^1,\dots,
b_{\alpha_{21}^n2^{k_1}}^1),\dots,(b_{\alpha_{2n}^n1}^n,
\dots,b_{\alpha_{2n}^n2^{k_n}}^n)) \\ \cdot \\ \cdot \\ \cdot \\
\xi^{n,2^{k_1},\dots,2^{k_n}}((b_{\alpha_{2^n1}^n1}^1,\dots,
b_{\alpha_{2^n1}^n2^{k_1}}^1),\dots,(b_{\alpha_{2^nn}^n1}^n,
\dots,b_{\alpha_{2^nn}^n2^{k_n}}^n))
\end{matrix}\right] \times$

$\left[\begin{matrix} \prod_{j=1}^{k_1} v_{\alpha_{1j}^{k_1}}^{1j}
\cdots \prod_{j=1}^{k_{n-2}} v_{\alpha_{1j}^{k_{n-2}}}^{n-2j}
\prod_{j=1}^{k_{n-1}} v_{\alpha_{1j}^{k_{n-1}}}^{n-1j}
\prod_{j=1}^{k_n} v_{\alpha_{1j}^{k_n}}^{nj} \\ \prod_{j=1}^{k_1}
v_{\alpha_{1j}^{k_1}}^{1j} \cdots \prod_{j=1}^{k_{n-2}}
v_{\alpha_{1j}^{k_{n-2}}}^{n-2j} \prod_{j=1}^{k_{n-1}}
v_{\alpha_{1j}^{k_{n-1}}}^{n-1j} \prod_{j=1}^{k_n}
v_{\alpha_{2j}^{k_n}}^{nj} \\ \cdot \\ \cdot \\ \cdot \\
\prod_{j=1}^{k_1} v_{\alpha_{1j}^{k_1}}^{1j} \cdots
\prod_{j=1}^{k_{n-2}} v_{\alpha_{1j}^{k_{n-2}}}^{n-2j}
\prod_{j=1}^{k_{n-1}} v_{\alpha_{1j}^{k_{n-1}}}^{n-1j}
\prod_{j=1}^{k_n} v_{\alpha_{2^{k_n}j}^{k_n}}^{nj} \\
\prod_{j=1}^{k_1} v_{\alpha_{1j}^{k_1}}^{1j} \cdots
\prod_{j=1}^{k_{n-2}} v_{\alpha_{1j}^{k_{n-2}}}^{n-2j}
\prod_{j=1}^{k_{n-1}} v_{\alpha_{2j}^{k_{n-1}}}^{n-1j}
\prod_{j=1}^{k_n} v_{\alpha_{1j}^{k_n}}^{nj} \\
\prod_{j=1}^{k_1} v_{\alpha_{1j}^{k_1}}^{1j} \cdots
\prod_{j=1}^{k_{n-2}} v_{\alpha_{1j}^{k_{n-2}}}^{n-2j}
\prod_{j=1}^{k_{n-1}} v_{\alpha_{2j}^{k_{n-1}}}^{n-1j}
\prod_{j=1}^{k_n} v_{\alpha_{2j}^{k_n}}^{nj}
\\ \cdot \\ \cdot \\ \cdot \\ \cdot \\ \cdot \\ \cdot \\
\prod_{j=1}^{k_1} v_{\alpha_{2^{k_1}j}^{k_1}}^{1j} \cdots
\prod_{j=1}^{k_{n-2}} v_{\alpha_{2^{k_{n-2}}j}^{k_{n-2}}}^{n-2j}
\prod_{j=1}^{k_{n-1}} v_{\alpha_{2^{k_{n-1}}j}^{k_{n-1}}}^{n-1j}
\prod_{j=1}^{k_n} v_{\alpha_{2^{k_n}j}^{k_n}}^{nj}
\end{matrix}\right]=$

$F^{k_1+\cdots+k_n}(T)[\xi^{k_1+\cdots+k_n}](F^0(B^{(11)}),
\dots,F^0(B^{(1k_1)}),\dots\dots,
F^0(B^{(n1)}),\dots,F^0(B^{(nk_n)}))$.

Notice that there are ball tangle diagrams
$B^{(1)},B^{(2)},B^{(3)}$ such that

$F^0(B^{(1)})=\left[\begin{matrix} 1 \\ 0
\end{matrix}\right],F^0(B^{(2)})=\left[\begin{matrix} 0 \\ 1
\end{matrix}\right],F^0(B^{(3)})=\left[\begin{matrix} 1 \\ 1
\end{matrix}\right]$, respectively (See Figure 8).

Therefore, by Lemma 3.1,

$F^{k_1+\cdots+k_n}(T^n(T^{k_1(1)},\dots,T^{k_n(n)}))$

$=F^n(T^n)\left[\begin{matrix}
\xi^{n,2^{k_1},\dots,2^{k_n}}((b_{\alpha_{11}^n1}^1,\dots,
b_{\alpha_{11}^n2^{k_1}}^1),\dots,(b_{\alpha_{1n}^n1}^n,
\dots,b_{\alpha_{1n}^n2^{k_n}}^n)) \\
\xi^{n,2^{k_1},\dots,2^{k_n}}((b_{\alpha_{21}^n1}^1,\dots,
b_{\alpha_{21}^n2^{k_1}}^1),\dots,(b_{\alpha_{2n}^n1}^n,
\dots,b_{\alpha_{2n}^n2^{k_n}}^n)) \\ \cdot \\ \cdot \\ \cdot \\
\xi^{n,2^{k_1},\dots,2^{k_n}}((b_{\alpha_{2^n1}^n1}^1,\dots,
b_{\alpha_{2^n1}^n2^{k_1}}^1),\dots,(b_{\alpha_{2^nn}^n1}^n,
\dots,b_{\alpha_{2^nn}^n2^{k_n}}^n))
\end{matrix}\right]=$

$F^n(T^n)\left[\begin{matrix} \prod_{j=1}^n
b_{\alpha_{1j}^n\alpha_{1j}^{n,2^{k_1},\dots,2^{k_n}}}^j &
\prod_{j=1}^n
b_{\alpha_{1j}^n\alpha_{2j}^{n,2^{k_1},\dots,2^{k_n}}}^j & \cdots
& \prod_{j=1}^n b_{\alpha_{1j}^n\alpha_{2^{k_1+\cdots+
k_n}j}^{n,2^{k_1},\dots,2^{k_n}}}^j
\\ \prod_{j=1}^n
b_{\alpha_{2j}^n\alpha_{1j}^{n,2^{k_1},\dots,2^{k_n}}}^j &
\prod_{j=1}^n
b_{\alpha_{2j}^n\alpha_{2j}^{n,2^{k_1},\dots,2^{k_n}}}^j & \cdots
& \prod_{j=1}^n b_{\alpha_{2j}^n\alpha_{2^{k_1+\cdots+
k_n}j}^{n,2^{k_1},\dots,2^{k_n}}}^j
\\ \cdot & \cdot & \cdot & \cdot \\ \cdot & \cdot & \cdot & \cdot
\\ \cdot & \cdot & \cdot & \cdot \\
\prod_{j=1}^n
b_{\alpha_{2^nj}^n\alpha_{1j}^{n,2^{k_1},\dots,2^{k_n}}}^j &
\prod_{j=1}^n
b_{\alpha_{2^nj}^n\alpha_{2j}^{n,2^{k_1},\dots,2^{k_n}}}^j &
\cdots & \prod_{j=1}^n b_{\alpha_{2^nj}^n\alpha_{2^{k_1+\cdots+
k_n}j}^{n,2^{k_1},\dots,2^{k_n}}}^j
\end{matrix}\right]$

$=F^n(T^n)[\eta^n](F^{k_1}(T^{k_1(1)}),\dots,F^{k_n}(T^{k_n(n)}))$.

This proves the theorem.
\end{proof}

Let us give the following example.

Suppose that $T^3,T^{1(1)},T^{1(2)},T^{1(3)}$ are
$3,1,1,1$-punctured ball tangle diagrams such that
$$F^3(T^3)=\left[\begin{matrix} a_{11} & a_{12} & a_{13} & a_{14}
& a_{15} & a_{16} & a_{17} & a_{18} \\ a_{21} & a_{22} & a_{23} &
a_{24} & a_{25} & a_{26} & a_{27} & a_{28} \end{matrix}\right],$$
$$F^1(T^{1(1)})=\left[\begin{matrix} b_{11}^1 & b_{12}^1 \\
b_{21}^1 & b_{22}^1 \end{matrix}\right],\,\, F^1(T^{1(2)})=\left[\begin{matrix} b_{11}^2 & b_{12}^2 \\
b_{21}^2 & b_{22}^2 \end{matrix}\right],\,\, F^1(T^{1(3)})=\left[\begin{matrix} b_{11}^3 & b_{12}^3 \\
b_{21}^3 & b_{22}^3 \end{matrix}\right],$$ respectively. Then
$T^3(T^{1(1)},T^{1(2)},T^{1(3)})$ is a $3$-punctured ball tangle
diagram and
$$F^3(T^3(T^{1(1)},T^{1(2)},T^{1(3)}))=F^3(T^3)[\eta^3](F^1(T^{1(1)}),F^1(T^{1(2)}),F^1(T^{1(3)}))$$
$$=\left[\begin{matrix} a_{11}
& a_{12} & a_{13} & a_{14} & a_{15} & a_{16} & a_{17} & a_{18} \\
a_{21} & a_{22} & a_{23} & a_{24} & a_{25} & a_{26} & a_{27} &
a_{28} \end{matrix}\right] \times$$
$$\left[\begin{matrix}
b_{11}^1b_{11}^2b_{11}^3 & b_{11}^1b_{11}^2b_{12}^3 &
b_{11}^1b_{12}^2b_{11}^3 & b_{11}^1b_{12}^2b_{12}^3 &
b_{12}^1b_{11}^2b_{11}^3 & b_{12}^1b_{11}^2b_{12}^3 & b_{12}^1b_{12}^2b_{11}^3 & b_{12}^1b_{12}^2b_{12}^3 \\
b_{11}^1b_{11}^2b_{21}^3 & b_{11}^1b_{11}^2b_{22}^3 &
b_{11}^1b_{12}^2b_{21}^3 & b_{11}^1b_{12}^2b_{22}^3 &
b_{12}^1b_{11}^2b_{21}^3 & b_{12}^1b_{11}^2b_{22}^3 & b_{12}^1b_{12}^2b_{21}^3 & b_{12}^1b_{12}^2b_{22}^3 \\
b_{11}^1b_{21}^2b_{11}^3 & b_{11}^1b_{21}^2b_{12}^3 &
b_{11}^1b_{22}^2b_{11}^3 & b_{11}^1b_{22}^2b_{12}^3 &
b_{12}^1b_{21}^2b_{11}^3 & b_{12}^1b_{21}^2b_{12}^3 & b_{12}^1b_{22}^2b_{11}^3 & b_{12}^1b_{22}^2b_{12}^3 \\
b_{11}^1b_{21}^2b_{21}^3 & b_{11}^1b_{21}^2b_{22}^3 &
b_{11}^1b_{22}^2b_{21}^3 & b_{11}^1b_{22}^2b_{22}^3 &
b_{12}^1b_{21}^2b_{21}^3 & b_{12}^1b_{21}^2b_{22}^3 & b_{12}^1b_{22}^2b_{21}^3 & b_{12}^1b_{22}^2b_{22}^3 \\
b_{21}^1b_{11}^2b_{11}^3 & b_{21}^1b_{11}^2b_{12}^3 &
b_{21}^1b_{12}^2b_{11}^3 & b_{21}^1b_{12}^2b_{12}^3 &
b_{22}^1b_{11}^2b_{11}^3 & b_{22}^1b_{11}^2b_{12}^3 & b_{22}^1b_{12}^2b_{11}^3 & b_{22}^1b_{12}^2b_{12}^3 \\
b_{21}^1b_{11}^2b_{21}^3 & b_{21}^1b_{11}^2b_{22}^3 &
b_{21}^1b_{12}^2b_{21}^3 & b_{21}^1b_{12}^2b_{22}^3 &
b_{22}^1b_{11}^2b_{21}^3 & b_{22}^1b_{11}^2b_{22}^3 & b_{22}^1b_{12}^2b_{21}^3 & b_{22}^1b_{12}^2b_{22}^3 \\
b_{21}^1b_{21}^2b_{11}^3 & b_{21}^1b_{21}^2b_{12}^3 &
b_{21}^1b_{22}^2b_{11}^3 & b_{21}^1b_{22}^2b_{12}^3 &
b_{22}^1b_{21}^2b_{11}^3 & b_{22}^1b_{21}^2b_{12}^3 & b_{22}^1b_{22}^2b_{11}^3 & b_{22}^1b_{22}^2b_{12}^3 \\
b_{21}^1b_{21}^2b_{21}^3 & b_{21}^1b_{21}^2b_{22}^3 &
b_{21}^1b_{22}^2b_{21}^3 & b_{21}^1b_{22}^2b_{22}^3 &
b_{22}^1b_{21}^2b_{21}^3 & b_{22}^1b_{21}^2b_{22}^3 &
b_{22}^1b_{22}^2b_{21}^3 & b_{22}^1b_{22}^2b_{22}^3
\end{matrix}\right].$$

Next, let us consider `(outer) connect sums' of various
$n$-punctured ball tangle diagrams and their invariants. They will
be also very useful when we compute invariants of complicated
tangles. Given $k_1$ and $k_2$-punctured ball tangle diagrams
$T^{k_1(1)}$ and $T^{k_2(2)}$, we denote the `horizontal' and the
`vertical' connect sums of them by $T^{k_1(1)} +_h T^{k_2(2)}$ and
$T^{k_1(1)} +_v T^{k_2(2)}$, respectively (See Figure 7).

\bigskip
\centerline{\epsfxsize=5 in \epsfbox{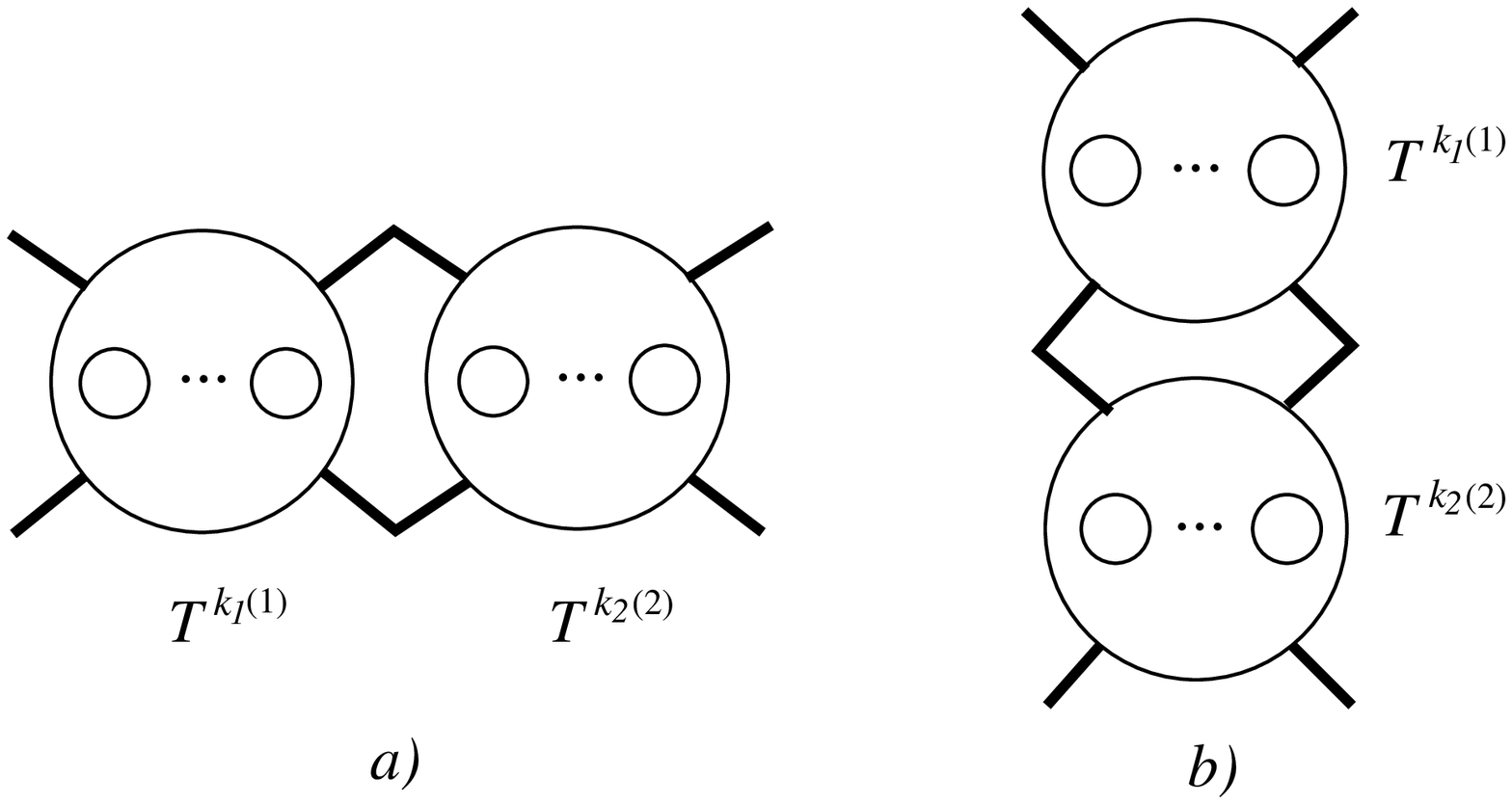}}
\medskip
\centerline{\small Figure 7. Connect sums of punctured ball
tangles.} \centerline{\small a) $T^{k_1(1)} +_h T^{k_2(2)}$, b)
$T^{k_1(1)} +_v T^{k_2(2)}$.}
\bigskip

\begin{thm} Let $k_1, k_2 \in \mathbb N \cup
\{0\}$, and let $T^{k_1(1)}, T^{k_2(2)}$ be $k_1, k_2$-punctured
ball tangle diagrams, respectively. Then

if $F^{k_1}(T^{k_1(1)})=\left[\begin{matrix} a_{11} & a_{12} &
\cdots & a_{12^{k_1}} \\ a_{21} & a_{22} & \cdots & a_{22^{k_1}}
\end{matrix}\right]$ and $F^{k_2}(T^{k_2(2)})=\left[\begin{matrix}
b_{11} & b_{12} & \cdots & b_{12^{k_2}} \\
b_{21} & b_{22} & \cdots & b_{22^{k_2}}
\end{matrix}\right]$, then

{\rm (1)} $F^{k_1+k_2}(T^{k_1(1)} +_h
T^{k_2(2)})=\left[\begin{matrix}\left(\begin{pmatrix}
a_{1i}b_{2j}+a_{2i}b_{1j} \\ a_{2i}b_{2j}
\end{pmatrix}_{j=1,\dots,2^{k_2}}\right)_{i=1,\dots,2^{k_1}}\end{matrix}\right]$,

{\rm (2)} $F^{k_1+k_2}(T^{k_1(1)} +_v
T^{k_2(2)})=\left[\begin{matrix}\left(\begin{pmatrix} a_{1i}b_{1j}
\\ a_{2i}b_{1j}+a_{1i}b_{2j}
\end{pmatrix}_{j=1,\dots,2^{k_2}}\right)_{i=1,\dots,2^{k_1}}\end{matrix}\right]$.
\end{thm}

\begin{proof} We denote $F^0(B^{(1)} +_h B^{(2)})$ by
$F^0(B^{(1)}) +_h F^0(B^{(2)})$ and $F^0(B^{(1)} +_v B^{(2)})$ by
$F^0(B^{(1)}) +_v F^0(B^{(2)})$ if $B^{(1)},B^{(2)} \in
\textbf{\textit{BT}}$.

(1) Let $T=T^{k_1(1)} +_h T^{k_2(2)}$, and let
$B^{(11)},\dots,B^{(1k_1)},B^{(21)},\dots,B^{(2k_2)} \in
\textbf{\textit{BT}}$ with
$$F^0(B^{(11)})=\left[\begin{matrix} v_1^{11} \\ v_2^{11} \end{matrix}\right],
\dots, F^0(B^{(1k_1)})=\left[\begin{matrix} v_1^{1k_1} \\
v_2^{1k_1} \end{matrix}\right],$$
$$F^0(B^{(21)})=\left[\begin{matrix} v_1^{21} \\ v_2^{21} \end{matrix}\right],
\dots, F^0(B^{(2k_2)})=\left[\begin{matrix} v_1^{2k_2} \\
v_2^{2k_2} \end{matrix}\right].$$

Then $T(B^{(11)},\dots,B^{(1k_1)},B^{(21)},\dots,B^{(2k_2)})$

$=T^{k_1(1)}(B^{(11)},\dots,B^{(1k_1)}) +_h
T^{k_2(2)}(B^{(21)},\dots,B^{(2k_2)})$ and

$F^0(T(B^{(11)},\dots,B^{(1k_1)},B^{(21)},\dots,B^{(2k_2)}))$

$=F^0(T^{k_1(1)}(B^{(11)},\dots,B^{(1k_1)}) +_h
T^{k_2(2)}(B^{(21)},\dots,B^{(2k_2)}))$

$=F^0(T^{k_1(1)}(B^{(11)},\dots,B^{(1k_1)})) +_h
F^0(T^{k_2(2)}(B^{(21)},\dots,B^{(2k_2)}))$

$=F^{k_1}(T^{k_1(1)})[\xi^{k_1}](F^0(B^{(11)}),\dots,F^0(B^{(1k_1)}))$

$+_h
F^{k_2}(T^{k_2(2)})[\xi^{k_2}](F^0(B^{(21)}),\dots,F^0(B^{(2k_2)}))$

$=\left[\begin{matrix} a_{11} & \cdots & a_{12^{k_1}} \\
a_{21} & \cdots & a_{22^{k_1}}
\end{matrix}\right]
\left[\begin{matrix} \prod_{j=1}^{k_1} v_{\alpha_{1j}^{k_1}}^{1j}
\\ \cdot \\ \cdot \\ \cdot \\
\prod_{j=1}^{k_1} v_{\alpha_{2^{k_1}j}^{k_1}}^{1j}
\end{matrix}\right] +_h
\left[\begin{matrix} b_{11} & \cdots & b_{12^{k_2}} \\
b_{21} & \cdots & b_{22^{k_2}}
\end{matrix}\right]
\left[\begin{matrix} \prod_{j=1}^{k_2} v_{\alpha_{1j}^{k_2}}^{2j}
\\ \cdot \\ \cdot \\ \cdot \\
\prod_{j=1}^{k_2} v_{\alpha_{2^{k_2}j}^{k_2}}^{2j}
\end{matrix}\right]=$

$$\left[\begin{matrix} a_{11} \prod_{j=1}^{k_1}
v_{\alpha_{1j}^{k_1}}^{1j} + \cdots +
a_{12^{k_1}} \prod_{j=1}^{k_1} v_{\alpha_{2^{k_1}j}^{k_1}}^{1j} \\
a_{21} \prod_{j=1}^{k_1} v_{\alpha_{1j}^{k_1}}^{1j} + \cdots +
a_{22^{k_1}} \prod_{j=1}^{k_1} v_{\alpha_{2^{k_1}j}^{k_1}}^{1j}
\end{matrix}\right] +_h
\left[\begin{matrix} b_{11} \prod_{j=1}^{k_2}
v_{\alpha_{1j}^{k_2}}^{2j} + \cdots +
b_{12^{k_2}} \prod_{j=1}^{k_2} v_{\alpha_{2^{k_2}j}^{k_2}}^{2j} \\
b_{21} \prod_{j=1}^{k_2} v_{\alpha_{1j}^{k_2}}^{2j} + \cdots +
b_{22^{k_2}} \prod_{j=1}^{k_2} v_{\alpha_{2^{k_2}j}^{k_2}}^{2j}
\end{matrix}\right]$$

$=\left[\begin{matrix} a_{11}b_{21}+a_{21}b_{11} & a_{21}b_{21} \\
\cdot & \cdot \\ \cdot & \cdot \\ \cdot & \cdot \\
a_{11}b_{22^{k_2}}+a_{21}b_{12^{k_2}} & a_{21}b_{22^{k_2}} \\
\cdot & \cdot \\ \cdot & \cdot \\ \cdot & \cdot \\
\cdot & \cdot \\ \cdot & \cdot \\ \cdot & \cdot \\
a_{12^{k_1}}b_{21}+a_{22^{k_1}}b_{11} & a_{22^{k_1}}b_{21} \\
\cdot & \cdot \\ \cdot & \cdot \\ \cdot & \cdot \\
a_{12^{k_1}}b_{22^{k_2}}+a_{22^{k_1}}b_{12^{k_2}} &
a_{22^{k_1}}b_{22^{k_2}} \end{matrix}\right]^\dag
\left[\begin{matrix} \prod_{j=1}^{k_1} v_{\alpha_{1j}^{k_1}}^{1j}
\prod_{j=1}^{k_2} v_{\alpha_{1j}^{k_2}}^{2j} \\
\cdot \\ \cdot \\ \cdot \\
\prod_{j=1}^{k_1} v_{\alpha_{1j}^{k_1}}^{1j}
\prod_{j=1}^{k_2} v_{\alpha_{2^{k_2}j}^{k_2}}^{2j} \\
\cdot \\ \cdot \\ \cdot \\ \cdot \\ \cdot \\ \cdot \\
\prod_{j=1}^{k_1} v_{\alpha_{1j}^{k_1}}^{1j}
\prod_{j=1}^{k_2} v_{\alpha_{1j}^{k_2}}^{2j} \\
\cdot \\ \cdot \\ \cdot \\
\prod_{j=1}^{k_1} v_{\alpha_{2^{k_1}j}^{k_1}}^{1j}
\prod_{j=1}^{k_2} v_{\alpha_{2^{k_2}j}^{k_2}}^{2j}
\end{matrix}\right]$

$=F^{k_1+k_2}(T)[\xi^{k_1+k_2}](F^0(B^{(11)}),
\dots,F^0(B^{(1k_1)}),F^0(B^{(21)}),\dots,F^0(B^{(2k_2)}))$.

Notice that there are ball tangle diagrams
$B^{(1)},B^{(2)},B^{(3)}$ such that

$F^0(B^{(1)})=\left[\begin{matrix} 1 \\ 0
\end{matrix}\right],F^0(B^{(2)})=\left[\begin{matrix} 0 \\ 1
\end{matrix}\right],F^0(B^{(3)})=\left[\begin{matrix} 1 \\ 1
\end{matrix}\right]$, respectively (See Figure 8).

Therefore, by Lemma 3.1,

$F^{k_1+k_2}(T^{k_1(1)} +_h
T^{k_2(2)})=\left[\begin{matrix}\left(\begin{pmatrix}
a_{1i}b_{2j}+a_{2i}b_{1j} \\ a_{2i}b_{2j}
\end{pmatrix}_{j=1,\dots,2^{k_2}}\right)_{i=1,\dots,2^{k_1}}\end{matrix}\right]$.

(2) Similarly, we can show that

$F^{k_1+k_2}(T^{k_1(1)} +_v
T^{k_2(2)})=\left[\begin{matrix}\left(\begin{pmatrix} a_{1i}b_{1j}
\\ a_{2i}b_{1j}+a_{1i}b_{2j}
\end{pmatrix}_{j=1,\dots,2^{k_2}}\right)_{i=1,\dots,2^{k_1}}\end{matrix}\right]$.

This proves the theorem.
\end{proof}

Notice that each of the horizontal connect sum and the vertical
connect sum of punctured ball tangles is associative but not
commutative. However, their invariants are not changed.

\begin{cor} Let $k_1, k_2 \in \mathbb N \cup
\{0\}$, and let $T^{k_1(1)}, T^{k_2(2)}$ be $k_1, k_2$-punctured
ball tangle diagrams, respectively. Then

{\rm (1)} $F^{k_1+k_2}(T^{k_1(1)} +_h
T^{k_2(2)})=F^{k_2+k_1}(T^{k_2(2)} +_h T^{k_1(1)})$,

{\rm (2)} $F^{k_1+k_2}(T^{k_1(1)} +_v
T^{k_2(2)})=F^{k_2+k_1}(T^{k_2(2)} +_v T^{k_1(1)})$.
\end{cor}

From now on, we denote simply by $F$ and $f$ for $F^1$ and $F^0$,
respectively. The following corollaries of our main theorems are
for the invariants of ball tangles and spherical tangles.

\begin{cor} If $A, B \in PM_{2\times2}(\mathbb Z)$ and

$A\left[\begin{matrix} 1 \\ 0 \end{matrix}\right]=
B\left[\begin{matrix} 1 \\ 0 \end{matrix}\right],
A\left[\begin{matrix} 0 \\ 1 \end{matrix}\right]=
B\left[\begin{matrix} 0 \\ 1 \end{matrix}\right],
A\left[\begin{matrix} 1 \\ 1 \end{matrix}\right]=
B\left[\begin{matrix} 1 \\ 1 \end{matrix}\right]$, then $A=B$.
\end{cor}

\begin{cor} If $S^{(1)}, S^{(2)} \in \textbf{\textit{ST}}$, then
$F(S^{(2)}(S^{(1)}))=F(S^{(2)})F(S^{(1)})$.
\end{cor}

\begin{cor} If $B^{(1)}, B^{(2)} \in \textbf{\textit{BT}}$ with
$f(B^{(1)})=\left[\begin{matrix} p \\ q \end{matrix}\right]$ and
$f(B^{(2)})=\left[\begin{matrix} r \\ s \end{matrix}\right]$, then

{\rm (1)} $f(B^{(1)} +_h B^{(2)})=\left[\begin{matrix} ps+qr \\
qs \end{matrix}\right]$ \rm{(Krebes \cite{K})},
{\rm (2)} $f(B^{(1)} +_v B^{(2)})=\left[\begin{matrix} pr \\
qr+ps \end{matrix}\right]$.
\end{cor}

\begin{cor} If $B \in \textbf{\textit{BT}}$ with $f(B)=\left[\begin{matrix} p \\
q \end{matrix}\right]$ and $S \in \textbf{\textit{ST}}$ with
$F(S)=\left[\begin{matrix} \alpha & \gamma \\ \beta & \delta
\end{matrix}\right]$, then

{\rm (1)} $F(B +_h S)=\left[\begin{matrix} p\beta + q\alpha &
p\delta + q\gamma \\ q\beta & q\delta
\end{matrix}\right]$, {\rm (2)} $F(B +_v S)=\left[\begin{matrix} p\alpha & p\gamma
\\ q\alpha + p\beta & q\gamma + p\delta
\end{matrix}\right]$.
\end{cor}

A connect sum of two spherical tangles is a $2$-punctured ball
tangle, so it has a $2\times2^2$ matrix in
$PM_{2\times2^2}(\mathbb Z)$. As a corollary of Theorem 3.3, we
give one more statement as follows.

\begin{cor} If $S^{(1)}, S^{(2)} \in \textbf{\textit{ST}}$ with
$F(S^{(1)})=\left[\begin{matrix} p & r \\ q & s
\end{matrix}\right]$ and
$F(S^{(2)})=\left[\begin{matrix} \alpha & \gamma \\ \beta & \delta
\end{matrix}\right]$, then

{\rm (1)} $F^2(S^{(1)} +_h S^{(2)})=\left[\begin{matrix} p\beta +
q\alpha & p\delta + q\gamma & r\beta + s\alpha & r\delta + s\gamma
\\ q\beta & q\delta & s\beta & s\delta
\end{matrix}\right]$,

{\rm (2)} $F^2(S^{(1)} +_v S^{(2)})=\left[\begin{matrix} p\alpha &
p\gamma & r\alpha & r\gamma \\ q\alpha + p\beta & q\gamma +
p\delta & s\alpha + r\beta & s\gamma + r\delta
\end{matrix}\right]$.
\end{cor}

Let us calculate the invariant for each of the ball tangles and
the spherical tangles in Figure 8.

1. The fundamental ball tangles \textbf{\textit{a}} and
\textbf{\textit{b}} have invariants $\left[\begin{matrix} 1
\\ 0 \end{matrix}\right]$ and $\left[\begin{matrix} 0 \\ 1
\end{matrix}\right]$, respectively.

2. The ball tangle \textbf{\textit{c}} has invariant
$\left[\begin{matrix} 1 \\ 1 \end{matrix}\right]$.

3. The spherical tangle \textbf{\textit{d}} is
$\textbf{\textit{I}}$ and has invariant $\left[\begin{matrix} 1 & 0 \\
0 & 1 \end{matrix}\right]$.

4. The spherical tangle \textbf{\textit{e}} has invariant
$\left[\begin{matrix} 1 \\ 1 \end{matrix}\right] +_v
\left[\begin{matrix} 1 & 0 \\ 0 & 1 \end{matrix}\right] =
\left[\begin{matrix} 1 & 0 \\ 1 & 1
\end{matrix}\right]$.

5. The spherical tangle \textbf{\textit{f}} has invariant
$\left[\begin{matrix} 1 \\ 0 \end{matrix}\right] +_h
\left[\begin{matrix} 1 & 0 \\ 1 & 1 \end{matrix}\right] =
\left[\begin{matrix} 1 & 1 \\ 0 & 0
\end{matrix}\right]$.

\vskip 0.2in

\bigskip
\centerline{\epsfxsize=5 in \epsfbox{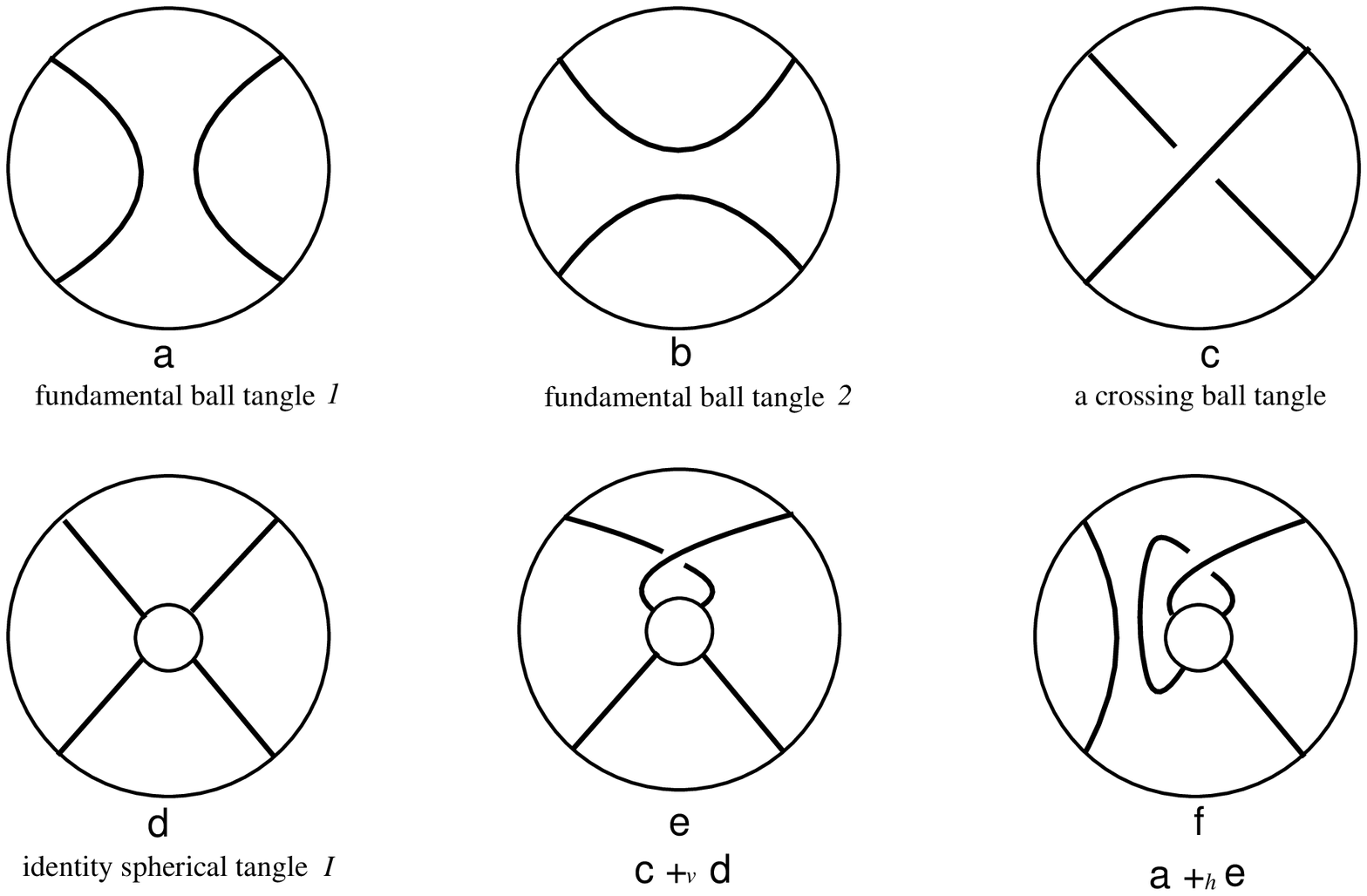}}
\medskip
\centerline{\small Figure 8. Ball tangle diagrams and spherical
tangle diagrams.}
\bigskip

When we denote the statement that $n$-punctured ball tangles
$T^{n(1)}$ and $T^{n(2)}$ induces the same function from
$\textbf{\textit{BT}}^n$ to $\textbf{\textit{BT}}$ by $T^{n(1)}
\simeq T^{n(2)}$ and $F^n(T^{n(1)})=F^n(T^{n(2)})$ by $T^{n(1)}
\sim T^{n(2)}$, $\simeq$ and $\sim$ are clearly equivalence
relations on $\textbf{\textit{nPBT}}$ and we have
$$T^{n(1)} \cong T^{n(2)} \Longrightarrow T^{n(1)} \simeq T^{n(2)}
\Longrightarrow T^{n(1)} \sim T^{n(2)}.$$ The first implication
comes from the definition of $\cong$ and the second implication is
proved by Theorem 2.14 and Lemma 3.1 immediately.

Notice that neither the converse of the first implication nor that
of the second implication is true (See Figure 9). In particular,
the spherical tangles $C$ and $D$ in Figure 9 have the matrix
$\left[\begin{matrix} 3 & 0 \\ 0 & 3 \end{matrix}\right]$ as
invariant. For another nonzero matrix invariant, we can take the
spherical tangle $A$ in Figure 9 and a spherical tangle $B'$
obtained from a single twist of the hole of $A$. We easily know
that $A$ and $B'$ are different functions. However, $A$ and $B'$
have the same invariant. By these reasons, we may consider the
equivalence relation $\simeq$ instead of $\cong$ for our
$n$-punctured ball tangle invariant.

This aspect is quite similar to that in Algebraic Topology in the
sense as follows:

If $X$ and $Y$ are pathconnected topological spaces, then
$$X \cong Y \Longrightarrow X \simeq Y \Longrightarrow X \sim Y,$$
where $X \cong Y$, $X \simeq Y$, and $X \sim Y$ mean the
statements that $X$ and $Y$ are topologically equivalent, $X$ and
$Y$ are homotopically equivalent, and $\pi_1(X)$ and $\pi_1(Y)$
are isomorphic, respectively.

\vskip 0.2in

\bigskip
\centerline{\epsfxsize=5.7 in \epsfbox{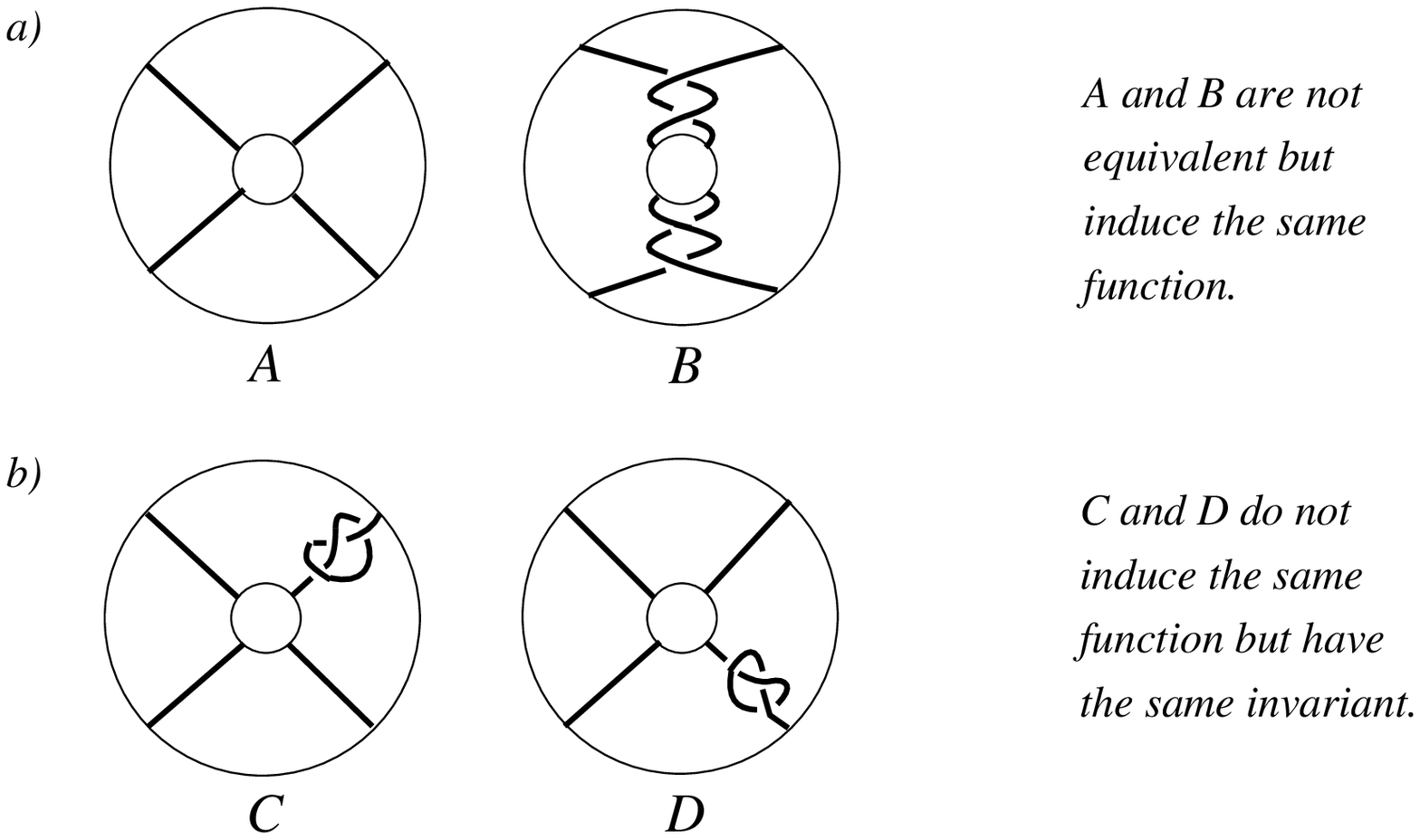}}
\medskip
\centerline{\small Figure 9. Tangles and functions.}
\bigskip

\begin{lm} [J.-W. Chung and X.-S. Lin \cite{C-L}] Let $\textbf{\textit{J}}$ be
the spherical tangle shown in Figure 10. Let $p_1,p_2,p_3,p_4$ be
the number of half twists inside of the balls marked by 1,2,3,4,
respectively. Then $$F(\textbf{\textit{J}})=\left[
\begin{matrix}p_1p_2p_3+p_1p_2p_4+p_1p_3p_4+p_2p_3p_4 &
-p_1p_3-p_1p_4-p_2p_3-p_2p_4 \\
p_1p_2+p_1p_4+p_3p_2+p_3p_4 & -p_1-p_2-p_3-p_4
\end{matrix}\right].$$ Therefore,
$${\rm det}\,F(\textbf{\textit{J}})=(p_1p_4-p_2p_3)^2.$$
\end{lm}

This is by a direct calculation.

Now, let us indicate a direct way to compute the invariant
$F(\textbf{\textit{J}})$ of the spherical tangle
$\textbf{\textit{J}}$ in Figure 10, in the special case of
$p_1=p_2=p_4=-4$ and $p_3=2$. Check with the formula in Lemma
3.10. Suppose that $T^5,B^{(1)},B^{(2)},B^{(3)},B^{(4)}$ are the
$5,0,0,0,0$-punctured ball tangle diagrams in Figure 10,
respectively. Then
$$\textbf{\textit{J}}=T^5(B^{(1)},B^{(2)},B^{(3)},B^{(4)},\textbf{\textit{I}}).$$
By Theorem 3.2, $F(\textbf{\textit{J}})=F^5(T^5)[\eta^5]
(f(B^{(1)}),f(B^{(2)}),f(B^{(3)}),f(B^{(4)}),F(\textbf{\textit{I}}))$.

We have
$f(B^{(1)})=f(B^{(2)})=f(B^{(4)})=\left[\begin{matrix} -4 \\
1 \end{matrix}\right]$ and $f(B^{(3)})=\left[\begin{matrix} 2 \\
1 \end{matrix}\right]$.

\vskip 0.2in

\bigskip
\centerline{\epsfxsize=3 in \epsfbox{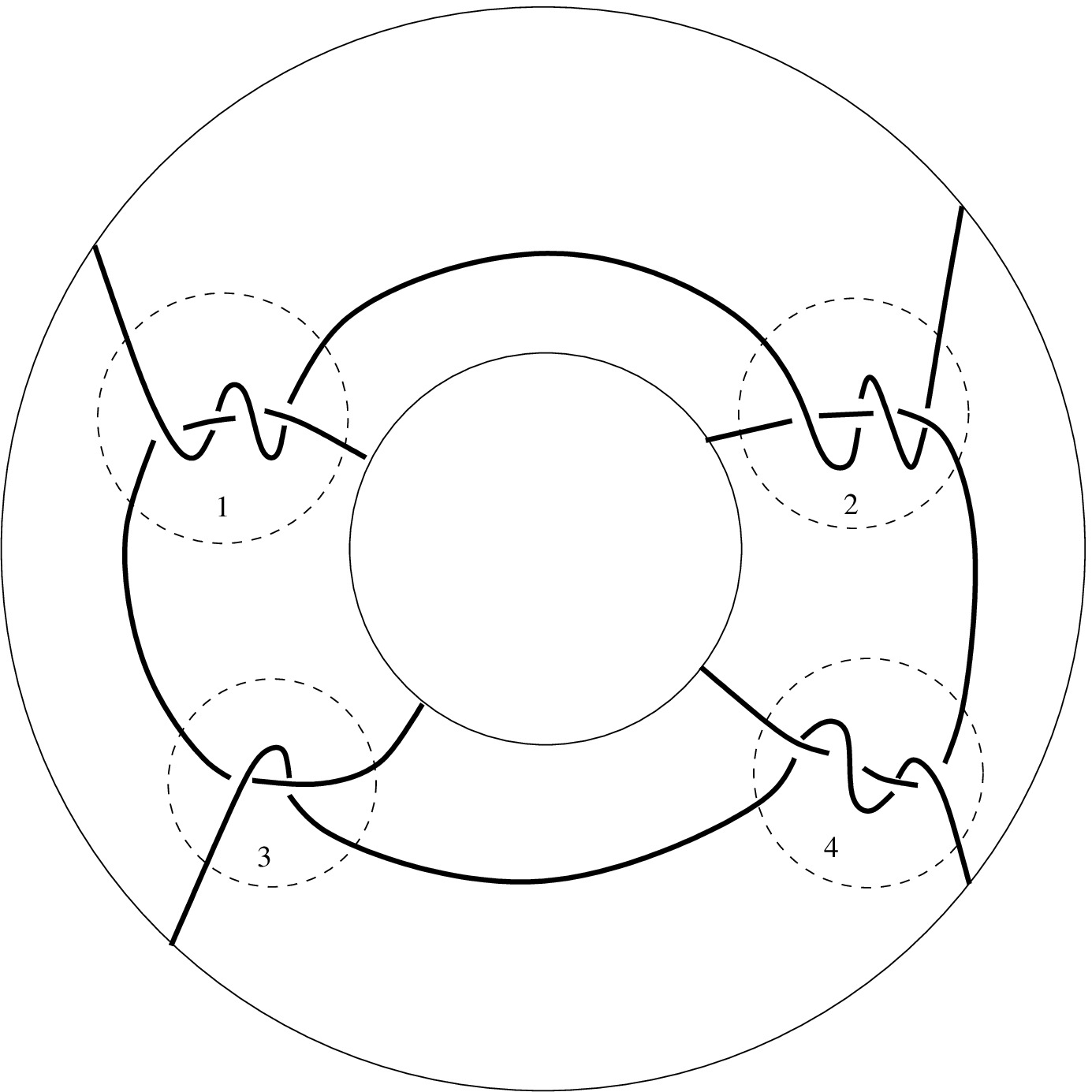}}
\medskip
\centerline{\small Figure 10. The spherical tangle
$\textbf{\textit{J}}$.}
\bigskip

\vskip 0.2in

\bigskip
\centerline{\epsfxsize=4.4 in \epsfbox{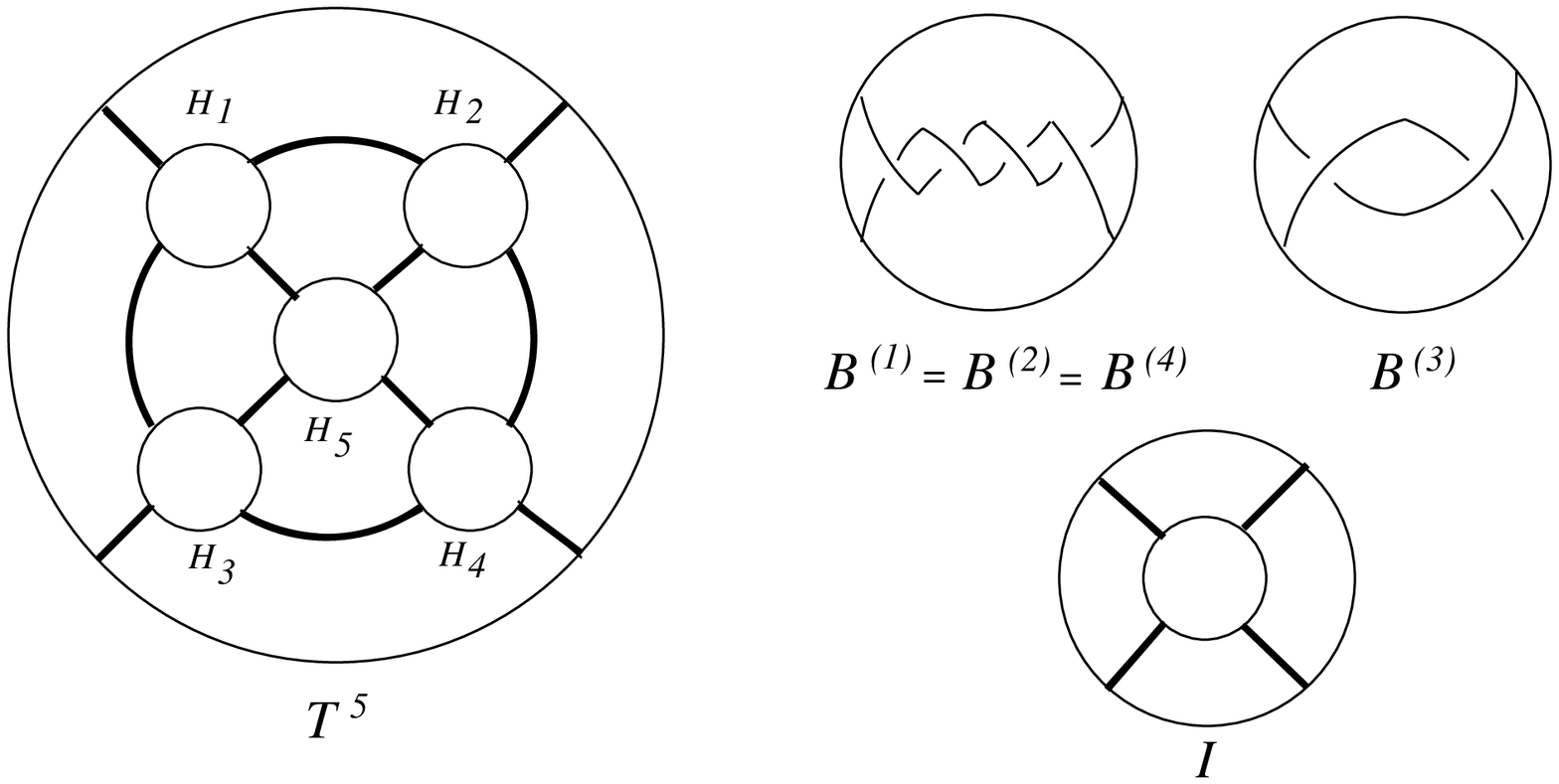}}
\medskip
\centerline{\small Figure 11. A decomposition of the spherical
tangle $\textbf{\textit{J}}$.}
\bigskip

First, let us compute $F^5(T^5)$ as the following steps:

1) The matrix $$\begin{pmatrix} \langle T_{1\alpha_1^5}^5 \rangle
& \cdots & \langle T_{1\alpha_{2^5}^5}^5 \rangle \\ \langle
T_{2\alpha_1^5}^5 \rangle & \cdots & \langle T_{2\alpha_{2^5}^5}^5
\rangle \end{pmatrix}$$ is
$$\begin{pmatrix} 0010 & 1000 & 1001 & 0100 & 1001 & 0100 & 0000 & 0000 \\
0000 & 0010 & 0000 & 1001 & 0010 & 0001 & 1001 & 0100
\end{pmatrix}.$$

2) Let $$F^5(T^5)=\left\{\begin{pmatrix} (-i)^{t_1}z\langle
T_{1\alpha_1^5}^5 \rangle & \cdots & (-i)^{t_{2^5}}z\langle T_{1\alpha_{2^5}^5}^5 \rangle \\
(-i)^{t_1}iz\langle T_{2\alpha_1^5}^5 \rangle & \cdots &
(-i)^{t_{2^5}}iz\langle T_{2\alpha_{2^5}^5}^5 \rangle
\end{pmatrix}\,|\,z \in \Phi\right\}
\cap M_{2\times2^5}(\mathbb Z).$$ Then the sequence $(t_k)_{1 \leq
k \leq 2^5}$ of exponents of $-i$ is $$\begin{matrix} 0112 & 1223
& 1223 & 2334 & 1223 & 2334 & 2334 & 3445 \end{matrix}.$$
Therefore, by taking $z=\pm i$, we have the invariant $F^5(T^5)$
as follows.
$$F^5(T^5)=\left[\begin{matrix}
0010 & 1000 & 100-1 & 0-100 & 100-1 & 0-100 & 0000 & 0000 \\
0000 & 0010 & 0000 & 100-1 & 0010 & 000-1 & 100-1 & 0-100
\end{matrix}\right].$$

Second, we compute $[\eta^5]
(f(B^{(1)}),f(B^{(2)}),f(B^{(3)}),f(B^{(4)}),F(\textbf{\textit{I}}))$
and describe it row-by-row as follows. That is, each pair of the
following means a row of the matrix $[\eta^5]
(f(B^{(1)}),f(B^{(2)}),f(B^{(3)}),f(B^{(4)}),F(\textbf{\textit{I}}))$.
$$128\,\,\,0;\, 0\,\,\,128;\, -32\,\,\,0;\, 0\,\,\,-32;\, 64\,\,\,0;\, 0\,\,\,64;\, -16\,\,\,0;\, 0\,\,\,-16;$$
$$-32\,\,\,0;\, 0\,\,\,-32;\, 8\,\,\,0;\, 0\,\,\,8;\, -16\,\,\,0;\, 0\,\,\,-16;\, 4\,\,\,0;\, 0\,\,\,4;$$
$$-32\,\,\,0;\, 0\,\,\,-32;\, 8\,\,\,0;\, 0\,\,\,8;\, -16\,\,\,0;\, 0\,\,\,-16;\, 4\,\,\,0;\, 0\,\,\,4;$$
$$8\,\,\,0;\, 0\,\,\,8;\, -2\,\,\,0;\, 0\,\,\,-2;\, 4\,\,\,0;\, 0\,\,\,4;\, -1\,\,\,0;\, 0\,\,\,-1.$$

Therefore, $$F(\textbf{\textit{J}})=F^5(T^5)[\eta^5]
(f(B^{(1)}),f(B^{(2)}),f(B^{(3)}),f(B^{(4)}),F(\textbf{\textit{I}}))$$
$$=\left[\begin{matrix} -32+64-32+0+0-32+0+0\ \ & 0+0+0-8+16+0-8+16 \\
-16-16+0+8+0+8+0+0\ \ & 0+0-4+0-4+0+2-4 \end{matrix}\right]$$
$$=\left[\begin{matrix} -32 & 16 \\ -16 & -10 \end{matrix}\right].$$

Also, we have
$${\rm det}\,F(\textbf{\textit{J}})=(-32)(-10)-16(-16)=576=24^2.$$
Thus, ${\rm det}\,F(\textbf{\textit{J}})$ is a square of integer.

We generalize lemma 3.10 as follows.

\begin{thm} Let $T^5$ be the $5$-punctured ball tangle shown in Figure
11, and let $B^{(i)}$ be ball tangles with $f(B^{(i)})=\left[\begin{matrix} p_i \\
q_i \end{matrix}\right]$ for each $i \in \{1,2,3,4\}$. Let
$\textbf{\textit{X}}=T^5(B^{(1)},B^{(2)},B^{(3)},B^{(4)},\textbf{\textit{I}}\,)$.
Then $$F(\textbf{\textit{X}})=\left[
\begin{smallmatrix}p_1p_2p_3q_4+p_1p_2q_3p_4+p_1q_2p_3p_4+q_1p_2p_3p_4\ \ &
-p_1q_2p_3q_4-p_1q_2q_3p_4-q_1p_2p_3q_4-q_1p_2q_3p_4 \\
p_1p_2q_3q_4+p_1q_2q_3p_4+q_1p_2p_3q_4+q_1q_2p_3p_4\ \ &
-p_1q_2q_3q_4-q_1p_2q_3q_4-q_1q_2p_3q_4-q_1q_2q_3p_4
\end{smallmatrix}\right].$$ Also, we have
$${\rm det}\,F(\textbf{\textit{X}})=(p_1q_2q_3p_4-q_1p_2p_3q_4)^2.$$
\end{thm}

The proof of Theorem 3.11 is quite similar as above for Lemma
3.10. Its proof is left to the reader. Notice that the determinant
of $F(\textbf{\textit{X}})$ is also a square of integer.

\vskip 0.1in

\bigskip
\centerline{\epsfxsize=4.5 in \epsfbox{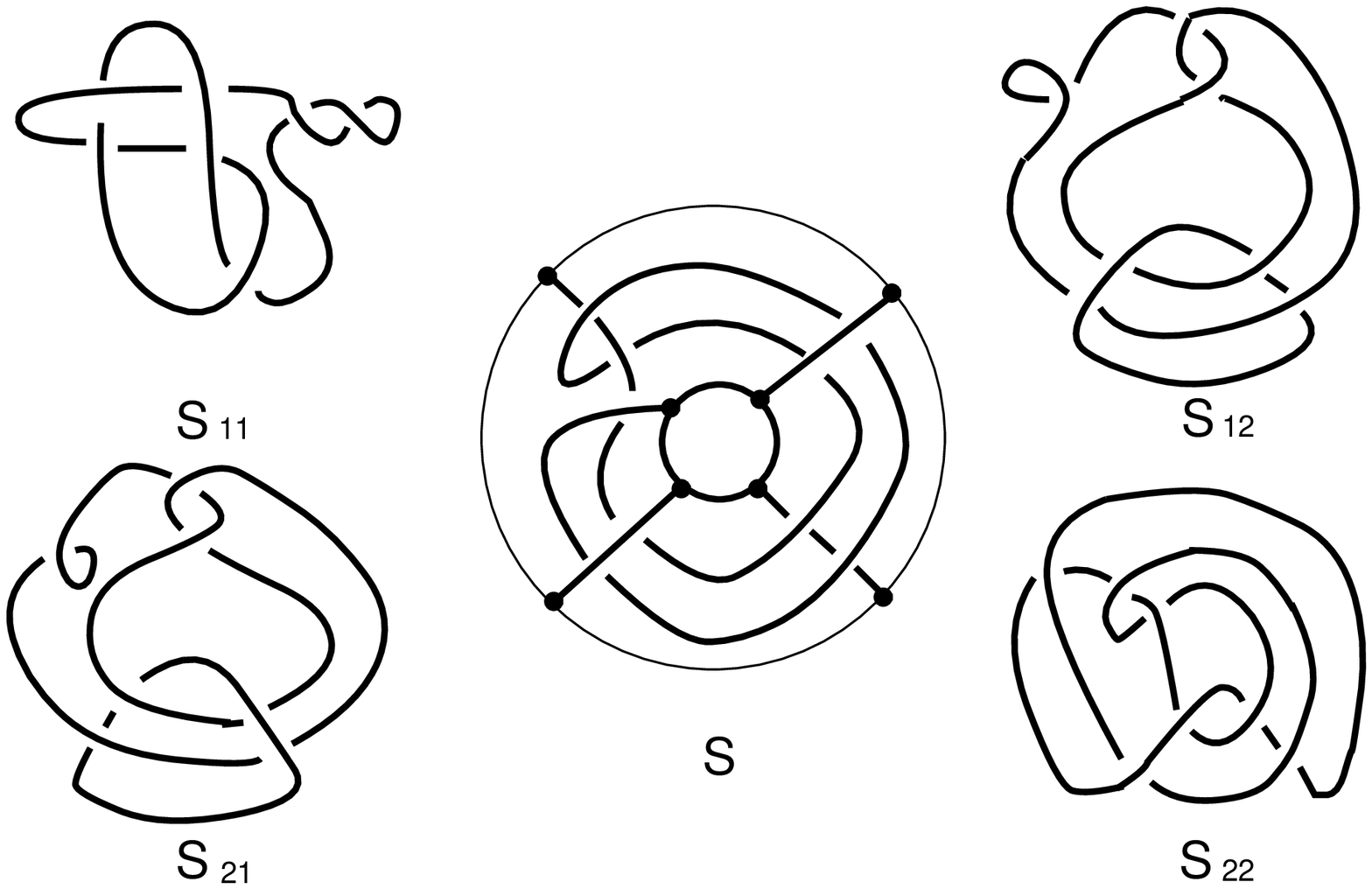}}
\medskip
\centerline{\small Figure 12. A spherical tangle $S$ cannot be
decomposed in terms of connect sums.}
\bigskip

As another example, let us consider the spherical tangle diagram
$S$ in Figure 12. Remark that the Kauffman bracket is a regular
isotopy invariant of link diagrams. By the definition of
invariant, we have
$$F(S)=\{\begin{pmatrix} z5A^3 & -iz8A \\ iz8A & z11A^{-1}
\end{pmatrix}|z \in \Phi\} \cap M_{2\times2}(Z)=\left[\begin{matrix}
5 & -8 \\ 8 & -11 \end{matrix}\right].$$ Hence, ${\rm
det}\,F(S)=-55-(-64)=9$. That is, ${\rm det}\,F(S)$ is a square of
integer. However, it seems that $S$ can not be decomposed in terms
of connect sums although we are not able to prove this fact.

For convenience, we use the following notation throughout the next
section:

(1) The subscripts 1,2 of ball tangles or spherical tangles will
no longer used to denote different kinds of closures. They will be
used simply to distinguish different ball tangles or spherical
tangles.

(2) $PM_2=PM_{2\times1}(\mathbb{Z})$ and
$PM_{2\times2}=PM_{2\times2}(\mathbb{Z})$.

\section{The Elementary operations on $PM_{2\times2}$ and Coxeter groups}

In this section, we introduce the group structure generated by the
elementary operations on $PM_{2\times2}$ induced by the elementary
operations on $\textbf{\textit{ST}}$.

Let us introduce the elementary operations on
$\textbf{\textit{ST}}$. Remark that a spherical tangle has exactly
2 holes which are inside and outside.

\begin{df} \cite{C-L}. Let $S$ be a spherical tangle diagram. Then

(1) $S^*$ is the mirror image of $S$,

(2) $S^-$ is the spherical tangle diagram obtained by
interchanging the inside hole with the outside hole of $S$,

(3) $S^{r_1}$ is the spherical tangle diagram obtained by only
rotating inside hole of $S$ $90^{\circ}$ counterclockwise on the
projection plane,

(4) $S^{r_2}$ is the spherical tangle diagram obtained by only
rotating outside hole of $S$ $90^{\circ}$ counterclockwise on the
projection plane,

(5) $S^R$ is the spherical tangle diagram obtained by the
$90^{\circ}$ rotation of $S$ itself counterclockwise on the
projection plane.
\end{df}

Note that $S^{r_2}=S^{-r_1-}$, $S^{r_1}=S^{-r_2-}$, and
$S^R=S^{r_1r_2}=S^{r_2r_1}$ for each $S \in \textbf{\textit{ST}}$.

\begin{lm} \cite{C-L}. If $S \in \textbf{\textit{ST}}$ with the invariant
$F(S)=\left[\begin{matrix} \alpha & \gamma \\ \beta & \delta
\end{matrix}\right]$, then

{\rm (1)} $F(S^*)=\left[\begin{matrix} \alpha & -\gamma \\
-\beta & \delta \end{matrix}\right]$,
{\rm (2)} $F(S^-)=\left[\begin{matrix} \delta & \gamma \\
\beta & \alpha \end{matrix}\right]$,
{\rm (3)} $F(S^{r_1})=\left[\begin{matrix} -\gamma & \alpha \\
-\delta & \beta \end{matrix}\right]$,

{\rm (4)} $F(S^{r_2})=\left[\begin{matrix} -\beta & -\delta \\
\alpha & \gamma \end{matrix}\right]$,
{\rm (5)} $F(S^R)=\left[\begin{matrix} \delta & -\beta \\
-\gamma & \alpha \end{matrix}\right]$.
\end{lm}

\bigskip
\centerline{\epsfxsize=4.5 in \epsfbox{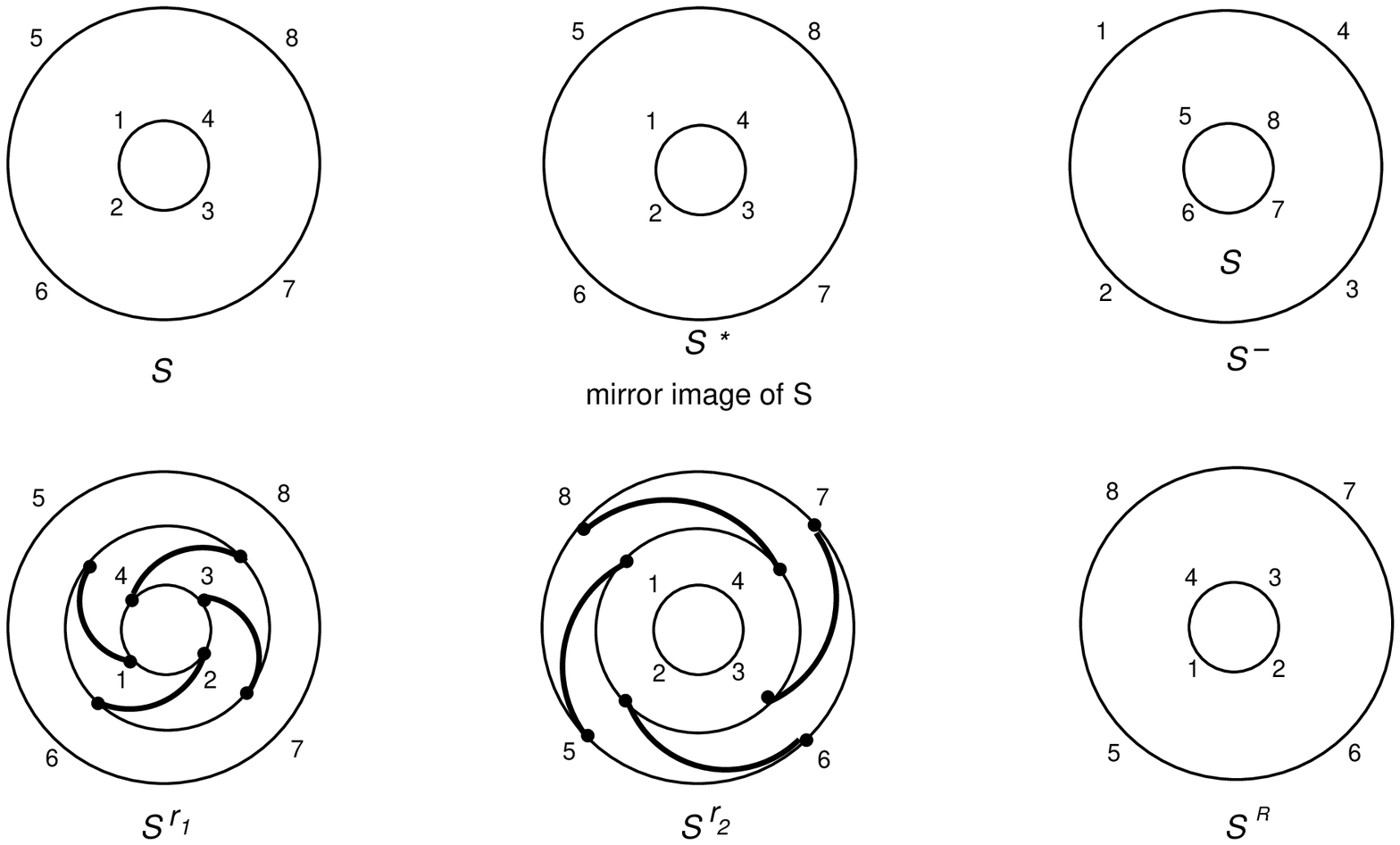}}
\medskip
\centerline{\small Figure 13. Elementary operations on
$\textbf{\textit{ST}}$.}
\bigskip

\begin{proof} Let $S \in \textbf{\textit{ST}}$ with
$F(S)=\left[\begin{matrix} \alpha & \gamma \\ \beta & \delta
\end{matrix}\right]$. Then there is $u \in \Phi$ such that
$\langle S_{11} \rangle=\alpha u$, $\langle S_{12} \rangle=\gamma
iu$, $\langle S_{21} \rangle=\beta (-i)u$, $\langle S_{22}
\rangle=\delta u$. Here the link $S_{ij}$, $i,j\in\{1,2\}$, is
obtained by taking the numerator closure ($i=1$) or the
denominator closure ($i=2$) of $S$ with its hole filled by the
fundamental tangle $j$. Therefore, $$\left[\begin{matrix} \alpha &
\gamma
\\ \beta & \delta
\end{matrix}\right]=\left[\begin{matrix} u^{-1}\alpha u &
u^{-1}(-i)\gamma iu
\\ u^{-1}i\beta (-i)u & u^{-1}\delta u
\end{matrix}\right].$$

Now we have

(1) $\langle S_{11}^* \rangle=\alpha u^{-1}, \langle S_{12}^*
\rangle=\gamma (iu)^{-1}, \langle S_{21}^* \rangle=\beta
(-iu)^{-1}, \langle S_{22}^* \rangle=\delta u^{-1}$,

(2) $\langle S_{11}^- \rangle=\langle S_{22} \rangle, \langle
S_{12}^- \rangle=\langle S_{12} \rangle, \langle S_{21}^-
\rangle=\langle S_{21} \rangle, \langle S_{22}^- \rangle=\langle
S_{11} \rangle$,

(3) $\langle S_{11}^{r_1} \rangle=\langle S_{12} \rangle, \langle
S_{12}^{r_1} \rangle=\langle S_{11} \rangle, \langle S_{21}^{r_1}
\rangle=\langle S_{22} \rangle, \langle S_{22}^{r_1}
\rangle=\langle S_{21} \rangle$.

Hence, $F(S^*)=\left[\begin{matrix} \alpha & -\gamma \\
-\beta & \delta \end{matrix}\right]$,
$F(S^-)=\left[\begin{matrix} \delta & \gamma \\
\beta & \alpha \end{matrix}\right]$,
$F(S^{r_1})=\left[\begin{matrix} \gamma & -\alpha \\
\delta & -\beta \end{matrix}\right]=
\left[\begin{matrix} -\gamma & \alpha \\
-\delta & \beta \end{matrix}\right]$.

Since $S^{r_2}=S^{-r_1-}$ and $S^R=S^{r_1r_2}$, (4) and (5) are
easily proved by (2) and (3).
\end{proof}

Like the case of ball tangle operations and invariants, it is
convenient to use the following notations.

Notation: Let $\left[\begin{matrix} \alpha & \gamma
\\ \beta & \delta
\end{matrix}\right] \in PM_{2\times2}$. Then

(1) $\left[\begin{matrix} \alpha & \gamma
\\ \beta & \delta
\end{matrix}\right]^*=\left[\begin{matrix} \alpha & -\gamma
\\ -\beta & \delta
\end{matrix}\right]$,
(2) $\left[\begin{matrix} \alpha & \gamma
\\ \beta & \delta
\end{matrix}\right]^-=\left[\begin{matrix} \delta & \gamma
\\ \beta & \alpha
\end{matrix}\right]$,
(3) $\left[\begin{matrix} \alpha & \gamma
\\ \beta & \delta
\end{matrix}\right]^{r_1}=\left[\begin{matrix} -\gamma & \alpha \\
-\delta & \beta \end{matrix}\right]$,

(4) $\left[\begin{matrix} \alpha & \gamma
\\ \beta & \delta
\end{matrix}\right]^{r_2}=\left[\begin{matrix} -\beta & -\delta \\
\alpha & \gamma \end{matrix}\right]$, (5) $\left[\begin{matrix}
\alpha & \gamma
\\ \beta & \delta
\end{matrix}\right]^R=\left[\begin{matrix} \delta & -\beta \\
-\gamma & \alpha \end{matrix}\right]$.

With these notations, we can write: $F(S^*)=F(S)^*, F(S^-)=F(S)^-,
F(S^{r_1})=F(S)^{r_1}, F(S^{r_2})=F(S)^{r_2}, F(S^R)=F(S)^R$ if $S
\in \textbf{\textit{ST}}$.

The determinant function {\rm det} is well-defined on
$PM_{2\times2}$ since ${\rm det}\,(-A)=(-1)^2\,{\rm det}\,A$ for
each $A \in PM_{2\times2}$.

Notice that the 5 elementary operations on $\textbf{\textit{ST}}$
do not change the determinant of invariants of spherical tangles.

Recall that $F(S_2 \circ S_1)=F(S_2)F(S_1)$ if $S_1, S_2 \in
\textbf{\textit{ST}}$ (Corollary 3.6).

\bigskip
\centerline{\epsfxsize=5 in \epsfbox{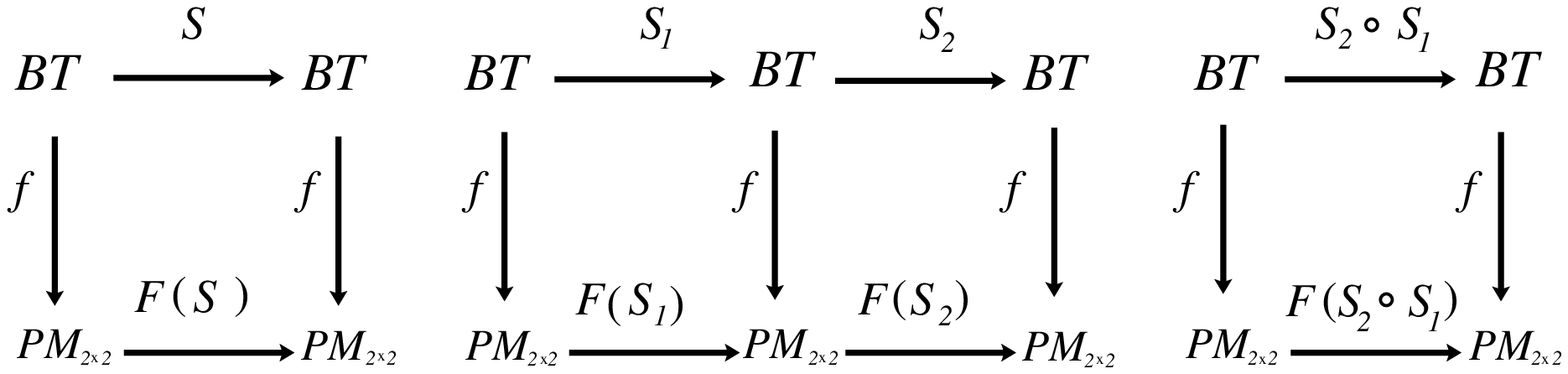}}
\medskip
\centerline{\small Figure 14. Commutative diagrams of invariants.}
\bigskip

The following lemma shows the elementary operations on the
composed spherical tangle.

\begin{lm} \cite{C-L}. If $S_1, S_2 \in \textbf{\textit{ST}}$, then

{\rm (1)} $(S_1 \circ S_2)^*=S_1^* \circ S_2^*$, {\rm (2)} $(S_1
\circ S_2)^-=S_2^- \circ S_1^-$, {\rm (3)} $(S_1 \circ
S_2)^{r_1}=S_1 \circ S_2^{r_1}$,

{\rm (4)} $(S_1 \circ S_2)^{r_2}=S_1^{r_2} \circ S_2$, {\rm (5)}
$(S_1 \circ S_2)^R=S_1^R \circ S_2^R$.
\end{lm}

\begin{df} An $n \times n$ matrix $M$ is called a Coxeter matrix
if $M_{ii}=1$ and $M_{ij}=M_{ji}>1$ for all $i,j \in
\{1,\dots,n\}$ with $i \neq j$, where $M_{ij}$ is the
$(i,j)$-entry of $M$.
\end{df}

\begin{df} Let $M$ be an $n \times n$ Coxeter matrix. Then a group
presented by $$\langle\, x_1,\dots,x_n\, |\,
(x_ix_j)^{M_{ij}}=1\,\, {\rm for}\,\, {\rm all}\,\, i,j \in
\{1,\dots,n\}\, \rangle,$$ denoted by $C_M$, is called the Coxeter
group with the Coxeter matrix $M$.
\end{df}

Let us think of the $5$ elementary operations $*, -, r_1, r_2, R$
on $PM_{2\times2}$ induced by the elementary operations on
$\textbf{\textit{ST}}$ as functions from $PM_{2\times2}$ to
$PM_{2\times2}$, respectively. For convenience, we use the
opposite composition of functions for the binary operation. For
instance, $-r_1$ means the composition $r_1 \circ -$. Recall that
$S^{r_2}=S^{-r_1-}$, $S^{r_1}=S^{-r_2-}$, and
$S^R=S^{r_1r_2}=S^{r_2r_1}$ for each $S \in \textbf{\textit{ST}}$
and observe the followings:

Suppose that $\left[\begin{matrix} \alpha & \gamma
\\ \beta & \delta \end{matrix}\right] \in PM_{2\times2}$. Then
$$\left[\begin{matrix} \alpha & \gamma
\\ \beta & \delta \end{matrix}\right] \underrightarrow{-}
\left[\begin{matrix} \delta & \gamma
\\ \beta & \alpha \end{matrix}\right] \underrightarrow{r_1}
\left[\begin{matrix} -\gamma & \delta
\\ -\alpha & \beta \end{matrix}\right] \underrightarrow{-}
\left[\begin{matrix} \beta & \delta
\\ -\alpha & -\gamma \end{matrix}\right] \underrightarrow{r_1}
\left[\begin{matrix} -\delta & \beta
\\ \gamma & -\alpha \end{matrix}\right],$$
$$\left[\begin{matrix} \alpha & \gamma
\\ \beta & \delta \end{matrix}\right] \underrightarrow{r_1}
\left[\begin{matrix} -\gamma & \alpha
\\ -\delta & \beta \end{matrix}\right] \underrightarrow{-}
\left[\begin{matrix} \beta & \alpha
\\ -\delta & -\gamma \end{matrix}\right] \underrightarrow{r_1}
\left[\begin{matrix} -\alpha & \beta
\\ \gamma & -\delta \end{matrix}\right] \underrightarrow{-}
\left[\begin{matrix} -\delta & \beta
\\ \gamma & -\alpha \end{matrix}\right]$$ and
$$\left[\begin{matrix} \alpha & \gamma
\\ \beta & \delta \end{matrix}\right] \underrightarrow{-}
\left[\begin{matrix} \delta & \gamma
\\ \beta & \alpha \end{matrix}\right] \underrightarrow{*}
\left[\begin{matrix} \delta & -\gamma
\\ -\beta & \alpha \end{matrix}\right],$$
$$\left[\begin{matrix} \alpha & \gamma
\\ \beta & \delta \end{matrix}\right] \underrightarrow{*}
\left[\begin{matrix} \alpha & -\gamma
\\ -\beta & \delta \end{matrix}\right] \underrightarrow{-}
\left[\begin{matrix} \delta & -\gamma
\\ -\beta & \alpha \end{matrix}\right]$$ and
$$\left[\begin{matrix} \alpha & \gamma
\\ \beta & \delta \end{matrix}\right] \underrightarrow{r_1}
\left[\begin{matrix} -\gamma & \alpha
\\ -\delta & \beta \end{matrix}\right] \underrightarrow{*}
\left[\begin{matrix} -\gamma & -\alpha
\\ \delta & \beta \end{matrix}\right],$$
$$\left[\begin{matrix} \alpha & \gamma
\\ \beta & \delta \end{matrix}\right] \underrightarrow{*}
\left[\begin{matrix} \alpha & -\gamma
\\ -\beta & \delta \end{matrix}\right] \underrightarrow{r_1}
\left[\begin{matrix} \gamma & \alpha
\\ -\delta & -\beta \end{matrix}\right].$$
Hence, we have $-r_1-r_1=r_1-r_1-$ and $-*=*-$ and $r_1*=*r_1$.
Also, $--$ and $r_1r_1$ and $**$ are the identity function from
$PM_{2\times2}$ to $PM_{2\times2}$. We show that the group
generated by the elementary operations on $PM_{2\times2}$ induced
by those on $\textbf{\textit{ST}}$ has the group presentation
$\langle\, x,y,z\, |\, x^2=y^2=z^2=1,\, xyxy=yxyx,\, xz=zx,\,
yz=zy\, \rangle$ which is a Coxeter group.

\begin{thm} The group $G(F)$ generated by the
elementary operations on $PM_{2\times2}$ induced by those on
$\textbf{\textit{ST}}$ has the group presentation $$\langle\,
x,y,z\, |\, x^2=y^2=z^2=1,\, xyxy=yxyx,\, xz=zx,\, yz=zy\,
\rangle.$$ Furthermore, $G(F)$ is isomorphic to the Coxeter group
$C_M$ with the Coxeter matrix $$M=\begin{pmatrix} 1 & 4 & 2 \\
4 & 1 & 2 \\ 2 & 2 & 1 \end{pmatrix}.$$ That is, $$G(F)=\langle\,
x,y,z\, |\,
x^2=y^2=z^2=(xy)^4=(yx)^4=(xz)^2=(zx)^2=(yz)^2=(zy)^2=1\,
\rangle.$$
\end{thm}

\begin{proof} Let $G=\langle\, x,y,z\, |\, x^2=y^2=z^2=1,\,
xyxy=yxyx,\, xz=zx,\, yz=zy\, \rangle$. Suppose that $\phi:G
\rightarrow G(F)$ is the epimorphism such that $\phi(x)=-$,
$\phi(y)=r_1$, $\phi(z)=*$. We claim that ${\rm Ker}\,\phi=\{1\}$.
Let $W(x,y,z)$ be a word in ${\rm Ker}\,\phi$. Then
$\phi(W(x,y,z))=W(-,r_1,*)=Id_{PM_{2\times2}}$. Since
$x^2=y^2=z^2=1$, we may assume that $W(x,y,z)$ has no consecutive
letters and no inverses of letters. Since $xz=zx$ and $yz=zy$ and
$z^2=1$, we have either $W(x,y,z)=W_1(x,y)z$ or
$W(x,y,z)=W_1(x,y)$ for some word $W_1(x,y)$ in $\{x,y\}$. We may
also assume that $W_1(x,y)$ has no consecutive letters and no
inverses of letters. We show that $W(x,y,z) \neq W_1(x,y)z$. If
$W(x,y,z)=W_1(x,y)z$, then $W_1(-,r_1)*=Id_{PM_{2\times2}}$. That
is, $W_1(-,r_1)=*$.

Observe that $$\begin{matrix} - \neq *, & -r_1 \neq *, & -r_1-
\neq *, & -r_1-r_1 \neq *, \\ r_1 \neq *, & r_1- \neq *, & r_1-r_1
\neq *, & r_1-r_1- \neq *. \end{matrix}$$ By $-r_1-r_1=r_1-r_1-$
and $-^2=r_1^2=Id_{PM_{2\times2}}$, we have $W_1(-,r_1) \neq *$.
This is a contradiction. Hence, $W(x,y,z) \neq W_1(x,y)z$.
Therefore, $W(x,y,z)=W_1(x,y)$ and the number of $z$ in $W(x,y,z)$
must be even.

Since $W_1(x,y)$ has no consecutive letters and no inverses of
letters, we have either

there are $k \in \mathbb N \cup \{0\}$ and $R \in
\{1,x,xy,xyx,xyxy,xyxyx,xyxyxy,xyxyxyx\}$

such that $W_1(x,y)=(xy)^{4k}R$ or

there are $k' \in \mathbb N \cup \{0\}$ and $R' \in
\{1,y,yx,yxy,yxyx,yxyxy,yxyxyx,yxyxyxy\}$

such that $W_1(x,y)=(yx)^{4k'}R'$.

Also, since $W_1(x,y)=W(x,y,z) \in {\rm Ker}\,\phi$, we have
either

$Id_{PM_{2\times2}}=W_1(-,r_1)=(-r_1)^{4k}\phi(R)$ or
$Id_{PM_{2\times2}}=W_1(-,r_1)=(r_1-)^{4k'}\phi(R')$.

Similarly, as above, observe that
$$\begin{matrix} - \neq
Id_{PM_{2\times2}}, & -r_1 \neq Id_{PM_{2\times2}}, & -r_1-
\neq Id_{PM_{2\times2}}, & -r_1-r_1 \neq Id_{PM_{2\times2}}, \\
r_1 \neq Id_{PM_{2\times2}}, & r_1- \neq Id_{PM_{2\times2}}, &
r_1-r_1 \neq Id_{PM_{2\times2}}, & r_1-r_1- \neq
Id_{PM_{2\times2}}.
\end{matrix}$$

Also, notice that $$\begin{matrix} -r_1-r_1-=r_1-r_1, &
-r_1-r_1-r_1=r_1-, \\ -r_1-r_1-r_1-=r_1, &
-r_1-r_1-r_1-r_1=Id_{PM_{2\times2}} \end{matrix}$$ and
$$\begin{matrix} r_1-r_1-r_1=-r_1-, & r_1-r_1-r_1-=-r_1, \\
r_1-r_1-r_1-r_1=-, & r_1-r_1-r_1-r_1-=Id_{PM_{2\times2}}.
\end{matrix}$$ Hence, we know that $\phi(R)=Id_{PM_{2\times2}}$ if and only if $R=1$
and $\phi(R')=Id_{PM_{2\times2}}$ if and only if $R'=1$.

Since $(-r_1)^{4k}=Id_{PM_{2\times2}}$ and
$(r_1-)^{4k'}=Id_{PM_{2\times2}}$, we have
$\phi(R)=Id_{PM_{2\times2}}$ and $\phi(R')=Id_{PM_{2\times2}}$.
Hence, $R=1$ and $R'=1$.

Therefore, $W_1(x,y)=(xy)^{4k}$ or $W_1(x,y)=(yx)^{4k}$ for some
$k \in \mathbb N \cup \{0\}$. Since $(xy)^4=1$ and $(yx)^4=1$,
$W_1(x,y)=1$. That is, $W(x,y,z)=1$. We have proved ${\rm
Ker}\,\phi=\{1\}$. Hence, $\phi:G \rightarrow G(F)$ is a group
isomorphism and $G(F)$ has the group presentation $\langle\,
x,y,z\, |\, x^2=y^2=z^2=1,\, xyxy=yxyx,\, xz=zx,\, yz=zy\,
\rangle$.

Now, we show that $G(F)$ is isomorphic to $C_M$. Since
$(xy)^2=(yx)^2$, $(xy)^2(xy)^2=(yx)^2(xy)^2$ and
$(xy)^2(yx)^2=(yx)^2(yx)^2$. Since $x^2=y^2=1$, $(xy)^4=(yx)^4=1$.
Also, since $xz=zx$, $(xz)(xz)=(zx)(xz)$ and $(xz)(zx)=(zx)(zx)$.
Since $x^2=z^2=1$, $(xz)^2=(zx)^2=1$. Similarly, since $yz=zy$,
$(yz)(yz)=(zy)(yz)$ and $(yz)(zy)=(zy)(zy)$. Since $y^2=z^2=1$,
$(yz)^2=(zy)^2=1$. Hence, the consequence of relators of $C_M$ is
contained in that of $G(F)$. Conversely, Since $(xy)^4=1$,
$(xy)^4(yx)^2=(yx)^2$. Since $x^2=y^2=1$, $(xy)^2=(yx)^2$. Also,
since $(xz)^2=1$, $(xz)^2(zx)=zx$. Since $x^2=z^2=1$, $xz=zx$.
Similarly, since $(yz)^2=1$, $(yz)^2(zy)=zy$. Since $y^2=z^2=1$,
$yz=zy$. Hence, the consequence of relators of $G(F)$ is contained
in that of $C_M$. Thus, $G(F)$ is isomorphic to $C_M$.
\end{proof}

We have just shown that the group $G(F)$ is a Coxeter group.
However, the group generated by the elementary operations on
$\textbf{\textit{ST}}$ is not a Coxeter group because $r_1$ on
$\textbf{\textit{ST}}$ has infinite order.

On the other hand, we showed the determinant of invariant of a
spherical tangle is a square of integer modulo 4 in \cite{C-L}.
However, it seems that the determinant is a square of integer even
though we don't know how to prove it yet.

\bigskip
\bigskip
\centerline {\textbf{APPENDIX: A guide to the nature of the
calculations}}

We have used so complicated notations to prove Theorem 3.2 which
is our first main theorem that most readers would probably feel
difficult to read the proof. However, to prove it precisely, we
could not help using such notations. Here, as this appendix, we
try to explain such complicated notations by concrete examples
with motivations to help to understand our proof of it. Also, we
introduce examples for the calculation of invariant of connect
sums looked like addition of fractions.

To explain the calculation process, we use elementary well-known
facts, in particular, expansion of product of several polynomials
by dictionary order, and finite sequences on the set $\{1,2\}$
which are combinations of our binary digits $1$ and $2$.

1. Examples of finite sequences on $\{1,2\}$:

For elements of linearly ordered set $J(n)$ by dictionary order,
we write as follows.

$$\alpha_1^1=(1), \,\,\,\,\, \alpha_{2^1}^1=(2).$$

$$\alpha_1^2=(1 1), \,\,\,\,\, \alpha_2^2=(1 2), \,\,\,\,\, \alpha_3^2=(2 1), \,\,\,\,\,
\alpha_{2^2}^2=(2 2).$$

$$\alpha_1^3=(1 1 1), \,\,\,\,\, \alpha_2^3=(1 1 2), \,\,\,\,\, \alpha_3^3=(1 2 1), \,\,\,\,\,
\alpha_4^3=(1 2 2),$$
$$\alpha_5^3=(2 1 1), \,\,\,\,\, \alpha_6^3=(2 1 2), \,\,\,\,\, \alpha_7^3=(2 2 1), \,\,\,\,\,
\alpha_{2^3}^3=(2 2 2).$$

$$\alpha_1^4=(1 1 1 1), \,\,\,\,\, \alpha_2^4=(1 1 1 2), \,\,\,\,\, \alpha_3^4=(1 1 2 1), \,\,\,\,\,
\alpha_4^4=(1 1 2 2),$$
$$\alpha_5^4=(1 2 1 1), \,\,\,\,\, \alpha_6^4=(1 2 1 2), \,\,\,\,\, \alpha_7^4=(1 2 2 1), \,\,\,\,\,
\alpha_8^4=(1 2 2 2),$$
$$\alpha_9^4=(2 1 1 1), \,\,\,\,\, \alpha_{10}^4=(2 1 1 2), \,\,\,\,\, \alpha_{11}^4=(2 1 2 1), \,\,\,\,\,
\alpha_{12}^4=(2 1 2 2),$$
$$\alpha_{13}^4=(2 2 1 1), \,\,\,\,\, \alpha_{14}^4=(2 2 1 2), \,\,\,\,\, \alpha_{15}^4=(2 2 2 1), \,\,\,\,\,
\alpha_{2^4}^4=(2 2 2 2).$$

Also, some examples of coordinates of above sequences are as
follows.
$$\alpha_{72}^3=2, \,\,\,\,\, \alpha_{32}^2=1, \,\,\,\,\, \alpha_{2^32}^3=2, \,\,\,\,\,
\alpha_{74}^4=1, \,\,\,\,\, \alpha_{2^43}^4=2.$$

2. Examples to key idea (motivation to dictionary order):

When we expand a product of several polynomials, we can use the
dictionary order as described. One of very complicated functions
$[\eta^n]$ which is the key for the proof of Theorem 3.2 is based
on the dictionary orders by which we expand the products of
several polynomials.

Let us explain the following two examples which involve our idea
for the main theorem.

(1) When $n=2$, $k_1=2^1$, $k_2=2^1$,
$$(a_1x_1+a_{2^1}x_{2^1})(b_1y_1+b_{2^1}y_{2^1})=a_1b_1x_1y_1+a_1b_{2^1}x_1y_{2^1}+
a_{2^1}b_1x_{2^1}y_1+a_{2^1}b_{2^1}x_{2^1}y_{2^1}$$
$$=\begin{pmatrix} a_1b_1 & a_1b_{2^1} & a_{2^1}b_1 & a_{2^1}b_{2^1} \end{pmatrix}
\begin{pmatrix} x_1y_1 \\ x_1y_{2^1} \\ x_{2^1}y_1 \\ x_{2^1}y_{2^1} \end{pmatrix}$$
$$=\xi^{2,2^1,2^1}((a_1, a_{2^1}), \,(b_1, b_{2^1})) \,\, \xi^{2,2^1,2^1}((x_1, x_{2^1}), \,(y_1, y_{2^1}))^\dag.$$

(2) When $n=3$, $k_1=2^1$, $k_2=2^2$, $k_3=2^2$,
$$(a_1x_1+a_{2^1}x_{2^1})(b_1y_1+b_2y_2+b_3y_3+b_{2^2}y_{2^2})(c_1z_1+c_2z_2+c_3z_3+c_{2^2}z_{2^2})$$
$$=a_1b_1c_1x_1y_1z_1+a_1b_1c_2x_1y_1z_2+a_1b_1c_3x_1y_1z_3+a_1b_1c_{2^2}x_1y_1z_{2^2}$$
$$+a_1b_2c_1x_1y_2z_1+a_1b_2c_2x_1y_2z_2+a_1b_2c_3x_1y_2z_3+a_1b_2c_{2^2}x_1y_2z_{2^2}$$
$$+a_1b_3c_1x_1y_3z_1+a_1b_3c_2x_1y_3z_2+a_1b_3c_3x_1y_3z_3+a_1b_3c_{2^2}x_1y_3z_{2^2}$$
$$+a_1b_{2^2}c_1x_1y_{2^2}z_1+a_1b_{2^2}c_2x_1y_{2^2}z_2+a_1b_{2^2}c_3x_1y_{2^2}z_3+a_1b_{2^2}c_{2^2}x_1y_{2^2}z_{2^2}$$
$$+a_{2^1}b_1c_1x_{2^1}y_1z_1+a_{2^1}b_1c_2x_{2^1}y_1z_2+a_{2^1}b_1c_3x_{2^1}y_1z_3+a_{2^1}b_1c_{2^2}x_{2^1}y_1z_{2^2}$$
$$+a_{2^1}b_2c_1x_{2^1}y_2z_1+a_{2^1}b_2c_2x_{2^1}y_2z_2+a_{2^1}b_2c_3x_{2^1}y_2z_3+a_{2^1}b_2c_{2^2}x_{2^1}y_2z_{2^2}$$
$$+a_{2^1}b_3c_1x_{2^1}y_3z_1+a_{2^1}b_3c_2x_{2^1}y_3z_2+a_{2^1}b_3c_3x_{2^1}y_3z_3+a_{2^1}b_3c_{2^2}x_{2^1}y_3z_{2^2}$$
$$+a_{2^1}b_{2^2}c_1x_{2^1}y_{2^2}z_1+a_{2^1}b_{2^2}c_2x_{2^1}y_{2^2}z_2+a_{2^1}b_{2^2}c_3x_{2^1}y_{2^2}z_3
+a_{2^1}b_{2^2}c_{2^2}x_{2^1}y_{2^2}z_{2^2}$$
$$=\xi^{3,2^1,2^2,2^2}((a_1, a_{2^1}), \,(b_1, b_2, b_3, b_{2^2}), \,(c_1, c_2, c_3, c_{2^2}))\times$$
$$\xi^{3,2^1,2^2,2^2}((x_1, x_{2^1}), \,(y_1, y_2, y_3, y_{2^2}), \,(z_1, z_2, z_3, z_{2^2}))^\dag.$$

3. An explanation of the proof of Theorem 3.2 by an example:

Let us consider the following example.

Suppose that $T^2,T^{2(1)},T^{1(2)}$ are $2,2,1$-punctured ball
tangle diagrams such that
$$F^2(T^2)=\left[\begin{matrix} a_{11} & a_{12} & a_{13} & a_{12^2} \\ a_{21} & a_{22} & a_{23} & a_{22^2}
\end{matrix}\right],$$
$$F^2(T^{2(1)})=\left[\begin{matrix} b_{11}^1 & b_{12}^1 & b_{13}^1 & b_{12^2}^1 \\
b_{21}^1 & b_{22}^1 & b_{23}^1 & b_{22^2}^1
\end{matrix}\right],\,\,
F^1(T^{1(2)})=\left[\begin{matrix} b_{11}^2 & b_{12^1}^2 \\
b_{21}^2 & b_{22^1}^2 \end{matrix}\right],$$ respectively. Notice
that $n=2$, $k_1=2$, $k_2=1$.

Let $T=T^2(T^{2(1)},T^{1(2)})$, and let
$B^{(11)},B^{(12)},B^{(21)} \in \textbf{\textit{BT}}$ with
$$F^0(B^{(11)})=\left[\begin{matrix} v_1^{11} \\ v_2^{11} \end{matrix}\right],
F^0(B^{(12)})=\left[\begin{matrix} v_1^{12} \\ v_2^{12}
\end{matrix}\right], F^0(B^{(21)})=\left[\begin{matrix} v_1^{21}
\\ v_2^{21} \end{matrix}\right].$$

Then
$T(B^{(11)},B^{(12)},B^{(21)})=T^2(T^{2(1)}(B^{(11)},B^{(12)}),T^{1(2)}(B^{(21)}))$
and

$F^0(T(B^{(11)},B^{(12)},B^{(21)}))=F^0(T^2(T^{2(1)}(B^{(11)},B^{(12)}),T^{1(2)}(B^{(21)})))$

$=F^2(T^2)[\xi^2](F^0(T^{2(1)}(B^{(11)},B^{(12)})),F^0(T^{1(2)}(B^{(21)})))$

$=F^2(T^2)[\xi^2](F^2(T^{2(1)})[\xi^2](F^0(B^{(11)}),F^0(B^{(12)})),
F^1(T^{1(2)})[\xi^1](F^0(B^{(21)})))$

$=F^2(T^2)[\xi^2](\left[\begin{matrix} b_{11}^1 & b_{12}^1 & b_{13}^1 & b_{12^2}^1 \\
b_{21}^1 & b_{22}^1 & b_{23}^1 & b_{22^2}^1 \end{matrix}\right]
\left[\begin{matrix}
\prod_{j=1}^2 v_{\alpha_{1j}^2}^{1j} \\ \prod_{j=1}^2 v_{\alpha_{2j}^2}^{1j} \\
\prod_{j=1}^2 v_{\alpha_{3j}^2}^{1j} \\ \prod_{j=1}^2
v_{\alpha_{2^2j}^2}^{1j}
\end{matrix}\right], \left[\begin{matrix} b_{11}^2 & b_{12^1}^2 \\ b_{21}^2 & b_{22^1}^2
\end{matrix}\right] \left[\begin{matrix} \prod_{j=1}^1 v_{\alpha_{1j}^1}^{2j}
\\ \prod_{j=1}^1 v_{\alpha_{2^1j}^1}^{2j}
\end{matrix}\right])$

$=F^2(T^2)[\xi^2](\left[\begin{matrix} b_{11}^1 \prod_{j=1}^2
v_{\alpha_{1j}^2}^{1j} + b_{12}^1 \prod_{j=1}^2
v_{\alpha_{2j}^2}^{1j} + b_{13}^1 \prod_{j=1}^2
v_{\alpha_{3j}^2}^{1j} + b_{12^2}^1 \prod_{j=1}^2
v_{\alpha_{2^2j}^2}^{1j} \\ b_{21}^1 \prod_{j=1}^2
v_{\alpha_{1j}^2}^{1j} + b_{22}^1 \prod_{j=1}^2
v_{\alpha_{2j}^2}^{1j} + b_{23}^1 \prod_{j=1}^2
v_{\alpha_{3j}^2}^{1j} + b_{22^2}^1 \prod_{j=1}^2
v_{\alpha_{2^2j}^2}^{1j}
\end{matrix}\right],$

$\left[\begin{matrix} b_{11}^2 \prod_{j=1}^1
v_{\alpha_{1j}^1}^{2j} +
b_{12^1}^2 \prod_{j=1}^1 v_{\alpha_{2^1j}^1}^{2j} \\
b_{21}^2 \prod_{j=1}^1 v_{\alpha_{1j}^1}^{2j} + b_{22^1}^2
\prod_{j=1}^1 v_{\alpha_{2^1j}^1}^{2j}
\end{matrix}\right])$

$=F^2(T^2)\left[\begin{matrix}
\xi^{2,2^2,2^1}((b_{\alpha_{11}^21}^1,b_{\alpha_{11}^22}^1,b_{\alpha_{11}^23}^1,b_{\alpha_{11}^22^2}^1),
(b_{\alpha_{12}^21}^2,b_{\alpha_{12}^22^1}^2)) \\
\xi^{2,2^2,2^1}((b_{\alpha_{21}^21}^1,b_{\alpha_{21}^22}^1,b_{\alpha_{21}^23}^1,b_{\alpha_{21}^22^2}^1),
(b_{\alpha_{22}^21}^2,b_{\alpha_{22}^22^1}^2)) \\
\xi^{2,2^2,2^1}((b_{\alpha_{31}^21}^1,b_{\alpha_{31}^22}^1,b_{\alpha_{31}^23}^1,b_{\alpha_{31}^22^2}^1),
(b_{\alpha_{32}^21}^2,b_{\alpha_{32}^22^1}^2)) \\
\xi^{2,2^2,2^1}((b_{\alpha_{2^21}^21}^1,b_{\alpha_{2^21}^22}^1,b_{\alpha_{2^21}^23}^1,
b_{\alpha_{2^21}^22^2}^1),(b_{\alpha_{2^22}^21}^2,
b_{\alpha_{2^22}^22^1}^2))
\end{matrix}\right]\times$

$\left[\begin{matrix} \xi^{2,2^2,2^1}((\prod_{j=1}^2
v_{\alpha_{1j}^2}^{1j},\prod_{j=1}^2 v_{\alpha_{2j}^2}^{1j},
\prod_{j=1}^2 v_{\alpha_{3j}^2}^{1j},\prod_{j=1}^2
v_{\alpha_{2^2j}^2}^{1j}), (\prod_{j=1}^1
v_{\alpha_{1j}^1}^{2j},\prod_{j=1}^1 v_{\alpha_{2^1j}^1}^{2j}))
\end{matrix}\right]^\dag$

$=F^2(T^2)\left[\begin{matrix}
\xi^{2,2^2,2^1}((b_{11}^1,b_{12}^1,b_{13}^1,b_{12^2}^1),(b_{11}^2,b_{12^1}^2)) \\
\xi^{2,2^2,2^1}((b_{11}^1,b_{12}^1,b_{13}^1,b_{12^2}^1),(b_{21}^2,b_{22^1}^2)) \\
\xi^{2,2^2,2^1}((b_{21}^1,b_{22}^1,b_{23}^1,b_{22^2}^1),(b_{11}^2,b_{12^1}^2)) \\
\xi^{2,2^2,2^1}((b_{21}^1,b_{22}^1,b_{23}^1,b_{22^2}^1),(b_{21}^2,b_{22^1}^2))
\end{matrix}\right]$
$[\xi^{2,2^2,2^1}](\left[\begin{matrix} v_1^{11}v_1^{12} \\
v_1^{11}v_2^{12} \\ v_2^{11}v_1^{12} \\ v_2^{11}v_2^{12}
\end{matrix}\right], \left[\begin{matrix} v_1^{21} \\ v_2^{21}
\end{matrix}\right])$

$=F^2(T^2)\left[\begin{matrix} b_{11}^1b_{11}^2 & b_{11}^1b_{12}^2
& b_{12}^1b_{11}^2 & b_{12}^1b_{12}^2 &
b_{13}^1b_{11}^2 & b_{13}^1b_{12}^2 & b_{14}^1b_{11}^2 & b_{14}^1b_{12}^2 \\
b_{11}^1b_{21}^2 & b_{11}^1b_{22}^2 & b_{12}^1b_{21}^2 &
b_{12}^1b_{22}^2 &
b_{13}^1b_{21}^2 & b_{13}^1b_{22}^2 & b_{14}^1b_{21}^2 & b_{14}^1b_{22}^2 \\
b_{21}^1b_{11}^2 & b_{21}^1b_{12}^2 & b_{22}^1b_{11}^2 &
b_{22}^1b_{12}^2 &
b_{23}^1b_{11}^2 & b_{23}^1b_{12}^2 & b_{24}^1b_{11}^2 & b_{24}^1b_{12}^2 \\
b_{21}^1b_{21}^2 & b_{21}^1b_{22}^2 & b_{22}^1b_{21}^2 &
b_{22}^1b_{22}^2 & b_{23}^1b_{21}^2 & b_{23}^1b_{22}^2 &
b_{24}^1b_{21}^2 & b_{24}^1b_{22}^2
\end{matrix}\right]$
$\left[\begin{matrix}
v_1^{11}v_1^{12}v_1^{21} \\ v_1^{11}v_1^{12}v_2^{21} \\ v_1^{11}v_2^{12}v_1^{21} \\
v_1^{11}v_2^{12}v_2^{21} \\ v_2^{11}v_1^{12}v_1^{21} \\ v_2^{11}v_1^{12}v_2^{21} \\
v_2^{11}v_2^{12}v_1^{21} \\ v_2^{11}v_2^{12}v_2^{21}
\end{matrix}\right]$

$=F^2(T^2)[\eta^2](F^2(T^{2(1)}),F^1(T^{1(2)}))[\xi^{2+1}](F^0(B^{(11)}),F^0(B^{(12)}),F^0(B^{(21)}))$

$=F^{2+1}(T)[\xi^{2+1}](F^0(B^{(11)}),F^0(B^{(12)}),F^0(B^{(21)}))$.

By Lemma 3.1, we conclude that

$F^3(T)=F^3(T^2(T^{2(1)},T^{1(2)}))=F^2(T^2)[\eta^2](F^2(T^{2(1)}),F^1(T^{1(2)}))=$

$\left[\begin{matrix} a_{11} & a_{12} & a_{13} & a_{14} \\ a_{21}
& a_{22} & a_{23} & a_{24}
\end{matrix}\right]\left[\begin{matrix} b_{11}^1b_{11}^2 & b_{11}^1b_{12}^2
& b_{12}^1b_{11}^2 & b_{12}^1b_{12}^2 &
b_{13}^1b_{11}^2 & b_{13}^1b_{12}^2 & b_{14}^1b_{11}^2 & b_{14}^1b_{12}^2 \\
b_{11}^1b_{21}^2 & b_{11}^1b_{22}^2 & b_{12}^1b_{21}^2 &
b_{12}^1b_{22}^2 &
b_{13}^1b_{21}^2 & b_{13}^1b_{22}^2 & b_{14}^1b_{21}^2 & b_{14}^1b_{22}^2 \\
b_{21}^1b_{11}^2 & b_{21}^1b_{12}^2 & b_{22}^1b_{11}^2 &
b_{22}^1b_{12}^2 &
b_{23}^1b_{11}^2 & b_{23}^1b_{12}^2 & b_{24}^1b_{11}^2 & b_{24}^1b_{12}^2 \\
b_{21}^1b_{21}^2 & b_{21}^1b_{22}^2 & b_{22}^1b_{21}^2 &
b_{22}^1b_{22}^2 & b_{23}^1b_{21}^2 & b_{23}^1b_{22}^2 &
b_{24}^1b_{21}^2 & b_{24}^1b_{22}^2
\end{matrix}\right]$.

Now, let us explain why Lemma 3.1 is required to complete this
example.

Suppose that $A=F^3(T^2(T^{2(1)},T^{1(2)}))$ and $B=F^2(T^2)[\eta^2](F^2(T^{2(1)}),F^1(T^{1(2)}))$.
Then $A$ and $B$ are matrices in $PM_{2\times2^3}(\mathbb Z)$.
In order to show $A=B$, we have shown that $AX=BX$ for each
$$X \in \{[\xi^3](F^0(B^{(11)}),F^0(B^{(12)}),F^0(B^{(21)}))|B^{(11)},B^{(12)},B^{(21)} \in \textbf{\textit{BT}}\}.$$
In \cite{C-L}, we proved that the 0-punctured ball tangle
invariant $F^0:\textbf{\textit{BT}} \rightarrow PM_{2\times1}(\mathbb Z)$ is surjective.
Note that $PM_{2\times1}(\mathbb Z)=P\mathbb{Z}^{2\dag}$. So we have
$$\{[\xi^3](F^0(B^{(11)}),F^0(B^{(12)}),F^0(B^{(21)}))|B^{(11)},B^{(12)},B^{(21)} \in \textbf{\textit{BT}}\}$$
$$=[\xi^3](P\mathbb{Z}^{2\dag} \times P\mathbb{Z}^{2\dag} \times P\mathbb{Z}^{2\dag}).$$
Fortunately, we have ball tangles $B^{(1)}$, $B^{(2)}$, $B^{(3)}$ whose invariants
are $\left[\begin{matrix} 1 \\ 0 \end{matrix}\right]$,
$\left[\begin{matrix} 0 \\ 1 \end{matrix}\right]$,
$\left[\begin{matrix} 1 \\ 1 \end{matrix}\right]$, respectively (See Figure 8).
Note that $[\xi^3](P\mathbb{Z}^{2\dag} \times P\mathbb{Z}^{2\dag} \times P\mathbb{Z}^{2\dag})
\subsetneqq  P\mathbb{Z}^{8\dag}$.
For example, we easily know that $\left[\begin{matrix} 0 & 1 & 1 & 1 & 1 & 1 & 1 & 1 \end{matrix}\right]^\dag \notin
[\xi^3](P\mathbb{Z}^{2\dag} \times P\mathbb{Z}^{2\dag} \times P\mathbb{Z}^{2\dag})$.

Lemma 3.1 says that we have only to show that $AX=BX$ for each
$$X \in [\xi^3](\{\left[\begin{matrix} 1 \\ 0 \end{matrix}\right],
\left[\begin{matrix} 0 \\ 1 \end{matrix}\right], \left[\begin{matrix} 1 \\ 1 \end{matrix}\right]\} \times
\{\left[\begin{matrix} 1 \\ 0 \end{matrix}\right],
\left[\begin{matrix} 0 \\ 1 \end{matrix}\right], \left[\begin{matrix} 1 \\ 1 \end{matrix}\right]\} \times
\{\left[\begin{matrix} 1 \\ 0 \end{matrix}\right],
\left[\begin{matrix} 0 \\ 1 \end{matrix}\right], \left[\begin{matrix} 1 \\ 1 \end{matrix}\right]\}).$$
That is, we have only to check the following 27 column vectors in $P\mathbb{Z}^{8\dag}$:
$$[\,1\,0\,0\,0\,0\,0\,0\,0\,]^\dag,[\,0\,1\,0\,0\,0\,0\,0\,0\,]^\dag,
[\,1\,1\,0\,0\,0\,0\,0\,0\,]^\dag,[\,0\,0\,1\,0\,0\,0\,0\,0\,]^\dag,[\,0\,0\,0\,1\,0\,0\,0\,0\,]^\dag,$$
$$[\,0\,0\,1\,1\,0\,0\,0\,0\,]^\dag,[\,1\,0\,1\,0\,0\,0\,0\,0\,]^\dag,
[\,0\,1\,0\,1\,0\,0\,0\,0\,]^\dag,[\,1\,1\,1\,1\,0\,0\,0\,0\,]^\dag,[\,0\,0\,0\,0\,1\,0\,0\,0\,]^\dag,$$
$$[\,0\,0\,0\,0\,0\,1\,0\,0\,]^\dag,[\,0\,0\,0\,0\,1\,1\,0\,0\,]^\dag,
[\,0\,0\,0\,0\,0\,0\,1\,0\,]^\dag,[\,0\,0\,0\,0\,0\,0\,0\,1\,]^\dag,[\,0\,0\,0\,0\,0\,0\,1\,1\,]^\dag,$$
$$[\,0\,0\,0\,0\,1\,0\,1\,0\,]^\dag,[\,0\,0\,0\,0\,0\,1\,0\,1\,]^\dag,
[\,0\,0\,0\,0\,1\,1\,1\,1\,]^\dag,[\,1\,0\,0\,0\,1\,0\,0\,0\,]^\dag,[\,0\,1\,0\,0\,0\,1\,0\,0\,]^\dag,$$
$$[\,1\,1\,0\,0\,1\,1\,0\,0\,]^\dag,[\,0\,0\,1\,0\,0\,0\,1\,0\,]^\dag,
[\,0\,0\,0\,1\,0\,0\,0\,1\,]^\dag,[\,0\,0\,1\,1\,0\,0\,1\,1\,]^\dag,[\,1\,0\,1\,0\,1\,0\,1\,0\,]^\dag,$$
$$[\,0\,1\,0\,1\,0\,1\,0\,1\,]^\dag,[\,1\,1\,1\,1\,1\,1\,1\,1\,]^\dag.$$

4. The invariant of connect sums looked like addition of
fractions:

Recall Corollary 3.7 to explain the calculation process by Theorem 3.3 which is our second main Theorem.
If $B^{(1)}, B^{(2)} \in \textbf{\textit{BT}}$ with
$F^0(B^{(1)})=\left[\begin{matrix} p \\ q \end{matrix}\right]$ and
$F^0(B^{(2)})=\left[\begin{matrix} r \\ s \end{matrix}\right]$, then
{\rm (1)} $F^0(B^{(1)} +_h B^{(2)})=\left[\begin{matrix} ps+qr \\
qs \end{matrix}\right]$ \rm{(Krebes \cite{K})},
{\rm (2)} $F^0(B^{(1)} +_v B^{(2)})=\left[\begin{matrix} pr \\
qr+ps \end{matrix}\right]$.

Consider the addition of fractions:
$$\frac {p}{q} + \frac {r}{s} = \frac {ps+qr}{qs},\,\,\, \frac {1}{{\frac {1}{\frac p q}}+{\frac {1}{\frac r s}}}
=\frac {pr}{ps+qr}.$$ They look like the invariant of connect sums of ball tangles.

Let us consider the following example.

Suppose that $T^{2(1)},T^{1(2)}$ are $2,1$-punctured ball
tangle diagrams such that
$$F^2(T^{2(1)})=\left[\begin{matrix} a_{11} & a_{12} & a_{13} & a_{12^2} \\
a_{21} & a_{22} & a_{23} & a_{22^2} \end{matrix}\right],\,\,
F^1(T^{1(2)})=\left[\begin{matrix} b_{11} & b_{12^1} \\ b_{21} & b_{22^1} \end{matrix}\right],$$ respectively.
Notice that $k_1=2$, $k_2=1$.

Let $T=T^{2(1)} +_h T^{1(2)}$, and let $B^{(11)},B^{(12)},B^{(21)} \in \textbf{\textit{BT}}$ with
$$F^0(B^{(11)})=\left[\begin{matrix} v_1^{11} \\ v_2^{11} \end{matrix}\right],
F^0(B^{(12)})=\left[\begin{matrix} v_1^{12} \\ v_2^{12} \end{matrix}\right],
F^0(B^{(21)})=\left[\begin{matrix} v_1^{21} \\ v_2^{21} \end{matrix}\right].$$
Then $T(B^{(11)},B^{(12)},B^{(21)})=T^{2(1)}(B^{(11)},B^{(12)}) +_h T^{1(2)}(B^{(21)})$ and

$F^0(T(B^{(11)},B^{(12)},B^{(21)}))=F^0(T^{2(1)}(B^{(11)},B^{(12)}) +_h T^{1(2)}(B^{(21)}))$

$=F^0(T^{2(1)}(B^{(11)},B^{(12)})) +_h F^0(T^{1(2)}(B^{(21)}))$

$=F^2(T^{2(1)})[\xi^2](F^0(B^{(11)}),F^0(B^{(12)})) +_h F^1(T^{1(2)})[\xi^1](F^0(B^{(21)}))$

$=\left[\begin{matrix} a_{11} & a_{12} & a_{13} & a_{12^2} \\ a_{21} & a_{22} & a_{23} & a_{22^2}
\end{matrix}\right]
\left[\begin{matrix} \prod_{j=1}^2 v_{\alpha_{1j}^2}^{1j} \\ \prod_{j=1}^2 v_{\alpha_{2j}^2}^{1j} \\
\prod_{j=1}^2 v_{\alpha_{3j}^2}^{1j} \\ \prod_{j=1}^2 v_{\alpha_{2^2j}^2}^{1j} \end{matrix}\right]
+_h \left[\begin{matrix} b_{11} & b_{12^1} \\ b_{21} & b_{22^1} \end{matrix}\right]
\left[\begin{matrix} \prod_{j=1}^1 v_{\alpha_{1j}^1}^{2j} \\ \prod_{j=1}^1 v_{\alpha_{2^1j}^1}^{2j}
\end{matrix}\right]$

$=\left[\begin{matrix} a_{11} \prod_{j=1}^2 v_{\alpha_{1j}^2}^{1j} +
a_{12} \prod_{j=1}^2 v_{\alpha_{2j}^2}^{1j} + a_{13} \prod_{j=1}^2 v_{\alpha_{3j}^2}^{1j} +
a_{12^2} \prod_{j=1}^2 v_{\alpha_{2^2j}^2}^{1j} \\
a_{21} \prod_{j=1}^2 v_{\alpha_{1j}^2}^{1j} +
a_{22} \prod_{j=1}^2 v_{\alpha_{2j}^2}^{1j} + a_{23} \prod_{j=1}^2 v_{\alpha_{3j}^2}^{1j} +
a_{22^2} \prod_{j=1}^2 v_{\alpha_{2^2j}^2}^{1j}
\end{matrix}\right] +_h$

$\left[\begin{matrix} b_{11} \prod_{j=1}^1 v_{\alpha_{1j}^1}^{2j} +
b_{12^1} \prod_{j=1}^1 v_{\alpha_{2^1j}^1}^{2j} \\
b_{21} \prod_{j=1}^1 v_{\alpha_{1j}^1}^{2j} +
b_{22^1} \prod_{j=1}^1 v_{\alpha_{2^1j}^1}^{2j}
\end{matrix}\right]$

$=\left[\begin{matrix}
\xi^{2,2^2,2^1}((a_{11},a_{12},a_{13},a_{12^2}),(b_{21},b_{22^1}))+
\xi^{2,2^2,2^1}((a_{21},a_{22},a_{23},a_{22^2}),(b_{11},b_{12^1})) \\
\xi^{2,2^2,2^1}((a_{21},a_{22},a_{23},a_{22^2}),(b_{21},b_{22^1}))
\end{matrix}\right] \times$

$\left[\begin{matrix} \xi^{2,2^2,2^1}((\prod_{j=1}^2
v_{\alpha_{1j}^2}^{1j},\prod_{j=1}^2 v_{\alpha_{2j}^2}^{1j},
\prod_{j=1}^2 v_{\alpha_{3j}^2}^{1j},\prod_{j=1}^2
v_{\alpha_{2^2j}^2}^{1j}), (\prod_{j=1}^1
v_{\alpha_{1j}^1}^{2j},\prod_{j=1}^1 v_{\alpha_{2^1j}^1}^{2j}))
\end{matrix}\right]^\dag$

$=[\begin{smallmatrix} a_{11}b_{21}+a_{21}b_{11} &
a_{11}b_{22^1}+a_{21}b_{12^1} & a_{12}b_{21}+a_{22}b_{11} &
a_{12}b_{22^1}+a_{22}b_{12^1} \\ a_{21}b_{21} & a_{21}b_{22^1} &
a_{22}b_{21} & a_{22}b_{22^1} \end{smallmatrix}$
$$\begin{smallmatrix} a_{13}b_{21}+a_{23}b_{11} &
a_{13}b_{22^1}+a_{23}b_{12^1} & a_{12^2}b_{21}+a_{22^2}b_{11} &
a_{12^2}b_{22^1}+a_{22^2}b_{12^1} \\ a_{23}b_{21} & a_{23}b_{22^1}
& a_{22^2}b_{21} & a_{22^2}b_{22^1}
\end{smallmatrix}]
[\xi^{2,2^2,2^1}](\left[\begin{matrix} v_1^{11}v_1^{12} \\
v_1^{11}v_2^{12} \\ v_2^{11}v_1^{12} \\ v_2^{11}v_2^{12}
\end{matrix}\right], \left[\begin{matrix} v_1^{21} \\ v_2^{21}
\end{matrix}\right])$$

$=[\begin{smallmatrix} a_{11}b_{21}+a_{21}b_{11} &
a_{11}b_{22^1}+a_{21}b_{12^1} & a_{12}b_{21}+a_{22}b_{11} &
a_{12}b_{22^1}+a_{22}b_{12^1} \\ a_{21}b_{21} & a_{21}b_{22^1} &
a_{22}b_{21} & a_{22}b_{22^1} \end{smallmatrix}$
$$\begin{smallmatrix} a_{13}b_{21}+a_{23}b_{11} &
a_{13}b_{22^1}+a_{23}b_{12^1} & a_{12^2}b_{21}+a_{22^2}b_{11} &
a_{12^2}b_{22^1}+a_{22^2}b_{12^1} \\ a_{23}b_{21} & a_{23}b_{22^1}
& a_{22^2}b_{21} & a_{22^2}b_{22^1}
\end{smallmatrix}]
\left[\begin{matrix}
v_1^{11}v_1^{12}v_1^{21} \\ v_1^{11}v_1^{12}v_2^{21} \\ v_1^{11}v_2^{12}v_1^{21} \\
v_1^{11}v_2^{12}v_2^{21} \\ v_2^{11}v_1^{12}v_1^{21} \\ v_2^{11}v_1^{12}v_2^{21} \\
v_2^{11}v_2^{12}v_1^{21} \\ v_2^{11}v_2^{12}v_2^{21}
\end{matrix}\right]$$

$=F^{2+1}(T)[\xi^{2+1}](F^0(B^{(11)}),F^0(B^{(12)}),F^0(B^{(21)}))$.

Therefore, by Lemma 3.1, we have

$F^{2+1}(T^{2(1)} +_h T^{1(2)})=[\begin{smallmatrix}
a_{11}b_{21}+a_{21}b_{11} & a_{11}b_{22^1}+a_{21}b_{12^1} &
a_{12}b_{21}+a_{22}b_{11} & a_{12}b_{22^1}+a_{22}b_{12^1} \\
a_{21}b_{21} & a_{21}b_{22^1} & a_{22}b_{21} & a_{22}b_{22^1}
\end{smallmatrix}$
$$\begin{smallmatrix} a_{13}b_{21}+a_{23}b_{11} &
a_{13}b_{22^1}+a_{23}b_{12^1} & a_{12^2}b_{21}+a_{22^2}b_{11} &
a_{12^2}b_{22^1}+a_{22^2}b_{12^1} \\ a_{23}b_{21} & a_{23}b_{22^1}
& a_{22^2}b_{21} & a_{22^2}b_{22^1}
\end{smallmatrix}].$$
Also, we can write $$F^{2+1}(T^{2(1)} +_h
T^{1(2)})=\left[\begin{matrix}\left(\begin{pmatrix}
a_{1i}b_{2j}+a_{2i}b_{1j} \\ a_{2i}b_{2j}
\end{pmatrix}_{j=1,2^1}\right)_{i=1,2,3,2^2}\end{matrix}\right].$$
Similarly, we can show the following formula for the vertical
connect sum.
$$F^{2+1}(T^{2(1)} +_v
T^{1(2)})=\left[\begin{matrix}\left(\begin{pmatrix} a_{1i}b_{1j}
\\ a_{2i}b_{1j}+a_{1i}b_{2j}
\end{pmatrix}_{j=1,2^1}\right)_{i=1,2,3,2^2}\end{matrix}\right].$$
Notice that the addition of fractions still plays an important
role in the calculation process of the invariant of connect sums
of punctured ball tangles.

We have tried to make our main theorems easier by concrete
examples. Even though we have used very complicated notations, we
think of our method as a kind of primitive applications of
dictionary orders.

\bigskip
\bigskip
\centerline {\textbf{Acknowledgement}}

I am a student of the late Professor Xiao-Song Lin of the
University of California, Riverside. This work had been done with
his invaluable advice and careful suggestions. Most part of the
paper is in my Ph.D. thesis. In this chance, I fixed a mistake in
the figure on page 43 of my dissertation. I made a wrong
counterexample with wrong pictures for a statement but I fixed it
now. Also, I gave a much easier counterexample than original one
on page 44 which is a supplementary explanation. In addition, I
have made the proof of Theorem 3.8 in my thesis more detailed
which is the first main theorem in this paper. Also, I would like
to thank Professor Scott Carter for encouraging me to add the
appendix which helps the readers.

\bigskip
\bigskip
MSC: 57M27

Keyword: Punctured ball tangle; Kauffman bracket; tangles embedded
in links; spherical tangle; group presentation; Coxeter group

\end{document}